\documentclass[appendixprefix=Appendix, 11pt]{article}
\usepackage[utf8]{inputenc}
\usepackage[margin=1in]{geometry}

\usepackage{amsthm}

\newcommand{\Halmos}{\ensuremath{\Box}}

\theoremstyle{definition}
\newtheorem{definition}{Definition}

\theoremstyle{plain}
\newtheorem{lemma}{Lemma}
\theoremstyle{plain}
\newtheorem{theorem}{Theorem}
\theoremstyle{plain}
\newtheorem{corollary}{Corollary}
\theoremstyle{plain}
\newtheorem{assumption}{Assumption}
\theoremstyle{plain}
\newtheorem{example}{Example}
\theoremstyle{plain}
\newtheorem{remark}{Remark}
\theoremstyle{plain}
\newtheorem{claim}{Claim}

\newcommand{\EMAIL}[1]{\texttt{#1}}

\title{Minimum Cut Representability of Stable Matching Problems}

\date{}
\author{
  Yuri Faenza\\
  {\footnotesize IEOR, Columbia University, New York, USA, \EMAIL{yf2414@columbia.edu}}
  \and
  Ayoub Foussoul\\
  {\footnotesize IEOR, Columbia University, New York, USA, \EMAIL{af3209@columbia.edu}}
  \and
  Chengyue He\\
  {\footnotesize IEOR, Columbia University, New York, USA, \EMAIL{ch3480@columbia.edu}}
}

\usepackage{graphicx}
\usepackage{multirow}
\usepackage{hhline}
\usepackage{cancel}
\usepackage{epstopdf}
\usepackage{wrapfig}

\usepackage[small, margin=1cm]{caption}
\usepackage[normalem]{ulem}
\usepackage{subcaption}

\usepackage{appendix}

\usepackage{color}
\definecolor{strcolor}{rgb}{0.6, 0.2, 0.6}
\definecolor{commentcolor}{rgb}{0.3125, 0.5, 0.3125}
\definecolor{keycol}{rgb}{0, 0, 1}

\usepackage{amsmath, amssymb}
\usepackage{bbm}

\usepackage{url}
\usepackage{hyperref}
\hypersetup{colorlinks=true,breaklinks=true,urlcolor=blue,citecolor=blue,linkcolor=blue,bookmarksopen=false,draft=false}

\DeclareMathOperator*{\cI}{\mathcal{I}}
\DeclareMathOperator*{\cD}{\mathcal{D}}

\DeclareMathOperator*{\cK}{\mathcal{K}}
\DeclareMathOperator*{\cP}{\mathcal{P}}
\DeclareMathOperator*{\cS}{\mathcal{S}}
\DeclareMathOperator*{\cA}{\mathcal{A}}
\DeclareMathOperator*{\cF}{\mathcal{F}}
\DeclareMathOperator*{\cR}{\mathcal{R}}
\DeclareMathOperator*{\cU}{\mathcal{U}}
\DeclareMathOperator*{\cL}{\mathcal{L}}
\DeclareMathOperator*{\cC}{\mathcal{C}}

\DeclareMathOperator*{\cval}{\mathsf{c}}

\newcommand{\AgentGreater}[2]{\succ_{#1}^{#2}}

\newcommand{\degree}[1]{\mathsf{deg}(#1)}

\newcommand{\val}[2]{{\mathsf{val}}_{#1}(#2)}

\usepackage{color-edits}
\addauthor{yf}{purple}
\addauthor{af}{blue}
\addauthor{ch}{orange}
\usepackage{enumerate}
\usepackage{todonotes}
\usepackage{tikz}
\usetikzlibrary{cd}
\tikzcdset{every label/.append style = {font = \small}}

\newcommand\vertarrowbox[3][1ex]{%
  \begin{array}[t]{@{}c@{}} #2 \vspace{-1mm} \\
  \left\uparrow\vcenter{\hrule height #1}\right.\kern-0.5ex\vspace{-1mm}\\
  \makebox[0mm]{\scriptsize#3-th}
  \vspace{-3.5mm}
  \vspace{2mm}
  \end{array}%
}

\usepackage[natbib=true, style=authoryear]{biblatex}
\begin{document}

\maketitle

\begin{abstract}%
We introduce and study {\em Minimum Cut Representability}, a framework to solve optimization and feasibility problems over stable matchings by representing them as minimum $s$-$t$ cut problems on digraphs over rotations. We provide necessary and sufficient conditions on objective functions and feasibility sets for problems to be minimum cut representable. In particular, we define the concepts of first and second order differentials of a function over stable matchings and show that a problem is minimum cut representable if and only if, roughly speaking, the objective function can be expressed solely using these differentials, and the feasibility set is a sublattice of the stable matching lattice. To demonstrate the practical relevance of our framework, we study a range of real-world applications, including problems involving school choice with siblings and a two-stage stochastic stable matching problem.
We show how our framework can be used to help solving these problems.
\end{abstract}

\section{Introduction}

Stability is a fundamental fairness concept in matching markets problems when we do not only wish to optimize a global objective function, but we also care that the output solution is fair at the level of individual agents. Since its introduction in the seminal work by~\cite{galeshapley}, stability has been employed in many real-world contexts, including matching medical residents to hospitals~\citep{roth2003origins}, students to schools~\citep{abdulkadirouglu2005new}, organ donors to recipients~\citep{huang2010circular}, drivers to riders in car-sharing~\citep{wang2018stable}, relocating refugees~\citep{delacretaz2023matching} and allocating resources in wireless networks~\citep{gu2015matching}, among others.

To compute a stable matching, \cite{galeshapley} proposed an algorithm that, given in input a two-sided market, outputs a stable matching that is optimal for one side of the market. This algorithm, known as \emph{Deferred Acceptance}, is now widely used in many applications, including most school matching markets (see, e.g.,~\cite{abdulkadirouglu2009strategy}). However, there are many situations where it might be desirable to compute a stable matching other than the one-sided optimal one. 

First, some applications require finding stable matchings that optimize specific objective functions. These functions may, for instance, model additional fairness criteria as is the case for \emph{minimum egalitarian} and \emph{sex-equal} stable marriages~\citep{gusfield1987three}, which aim at finding a stable matching that, roughly speaking, does not favor one side of the market over the other. We refer the reader to the literature review (Section~\ref{sec:additional-work}) for additional examples. 

Second, application-specific requirements often necessitate that the output matching meets additional constraints beyond stability. For instance, in refugee resettlement problems, it is important to find stable matchings that satisfy resource allocation constraints~\citep{delacretaz2023matching}. Similarly, in the school choice problem, it may be  desirable to compute diversity-aware stable matchings, which can be achieved by adding proportionality constraints~\citep{nguyen2019stable}. Moreover, in many situations, one wants to compute a stable matching that is robust to changes in the market~\citep{onlinestablematching,feigenbaum2020dynamic,fixedfirststage}. Therefore, developing algorithms that allow a central planner to output stable matchings other than the one-sided optimal is a significant and active area of research.

The goal of this paper is to introduce and investigate \emph{Minimum Cut Representability}, a framework that allows a central planner to solve a wide range of optimization and feasibility problems over stable matchings. We illustrate the usefulness of the framework by applying it to scenarios that involve sibling-related constraints and to situations where the market changes stochastically over time.

In the following, a {\em school matching instance} $(A,B,\AgentGreater{}{}, q)$ consists of a two-sided market with a set of students $A$, a set of schools $B$, and for every $a\in A$ (resp., $b \in B$), a strict order $\AgentGreater{a}{}$ (resp., $\AgentGreater{b}{}$) over $B^+ = B\cup\{\emptyset\}$ (resp., $A^+ = A\cup \{\emptyset\}$), where $\emptyset$ denotes the outside option. Each school $b \in B$ has a quota $q_b \leq |A|$ which gives the maximum number of students that can be assigned to the school. An element of $ A\cup B$ is called an \emph{agent}. A matching of students to schools is given by a set of pairs $M \subset A^+ \times B^+$ where every student $a \in A$ is paired with either a unique school or with the outside option, and every school $b \in B$ is paired with either exactly $q_b$ students or at most $q_b-1$ students and the outside option. Given a matching $M$, we denote by $M(a) \in B^+$ the unique partner of student $a$ and denote by $M(b) \subset A^+$ the partners of a school $b$. A {\em stable matching} of a school matching instance is a matching where no {\em blocking pairs} or {\em blocking agents} exist. A pair $ab \in A \times B$ is called a blocking pair if $b \succ_a M(a)$ and there exists $a' \in M(b)$ such that $a \succ_b a'$; a student $a \in A$ is called a blocking agent if they prefer the outside option to their assigned school; and a school $b \in B$ is called a blocking agent if they prefer the outside option to at least one of their assigned students. In words, a matching is stable if and only if no pair of agents would prefer each other to their partner (or one of their partners) in the matching, and no agent would prefer the outside option to their partner (or one of their partners) in the matching.

\smallskip

Given a school matching instance $(A,B,\AgentGreater{}{}, q)$, an optimization problem over the stable matchings of $(A,B,\AgentGreater{}{}, q)$ is a problem of the form 
% \vspace{-.4cm}
\begin{align*}
\label{problem}
 \tag{$\Pi(f, \cF)$}
     \min_{M \in \mathcal{F}} \quad & f(M)
\end{align*}

\noindent where $\cF$ is a subset of stable matchings of $(A,B,\AgentGreater{}{}, q)$ and $f: \cF \rightarrow \mathbb{R}$ is a real valued function.

\smallskip

\subsection{Contributions}

{\noindent \bf Necessary and Sufficient Conditions for Minimum Cut Representability.} In this work, we propose a framework to solve a large class of optimization problems over stable matchings. Our approach consists of transforming these problems into minimum $s$-$t$ cut problems over a directed graph, where the vertices are special cycles in the instance called {\em rotations} (see Definition~\ref{def:rotation-onemany}), along with a source vertex $s$ and a sink vertex $t$. More formally, consider a school matching instance $I$. Let $\cF$ be a subset of stable matchings of $I$ and $f: \cF \rightarrow \mathbb{R}$ a real-valued function over $\cF$. We say that the optimization problem \eqref{problem} is \emph{minimum cut representable} if and only if there exists a digraph $\cD$, referred to as a \emph{witness digraph}, whose vertices are the rotations of $I$ plus a source vertex $s$ and a sink vertex $t$, with non-negative (possibly infinite) capacities on arcs, such that finite capacitated $s$-$t$ cuts correspond exactly to the feasible stable matchings $\cF$, and each such cut has capacity, within an additive constant, equal to the value of $f$ on the corresponding matching (see Definition~\ref{def:mcs}). Provided that the witness digraph $\cD$ can be constructed efficiently, a minimum cut representable problem can therefore be solved efficiently using standard minimum $s$-$t$ cut algorithms.

We would like to note that while our work is not the first to reduce an optimization problem over stable matchings to a minimum $s$-$t$ cut problems over rotations, it introduces a much more general framework. In particular, a similar transformation has been used to solve the Minimum Weight Stable Matching (MWSM) problem, a special case of \eqref{problem} where ${\cal F}$ coincides with the whole family of stable matchings and the objective function $f$ is linear, with the value of a stable matching $M$ equal to the sum of the weights of the pairs in $M$ \citep{irving1987efficient}. Our framework, however, addresses more general feasibility sets and non-linear objectives. Most importantly, we provide a complete characterization of the set of problems that are minimum cut representable. 

In Theorem~\ref{thm:main1}, we establish necessary and sufficient conditions on the objective function $f$ and the feasibility set $\cF$ for the problem \eqref{problem} to be minimum cut representable. In particular, we introduce the notions of \emph{first} and \emph{second} order differentials of $f$ and show that the problem \eqref{problem} is minimum cut representable if and only if $\cF$ is a sublattice of the lattice of stable matchings, the second-order differentials of $f$ are non-negative, and $f = f^{\sf aprx}$ for a suitable function $f^{\sf aprx}: \cF \rightarrow \mathbb{R}$ defined in terms of the first and second-order differentials of $f$. Moreover, when a problem is minimum cut representable, we give an explicit construction of a witness digraph. Our construction is efficient for a wide range of applications. We also show that, the class of Minimum Weight Stable Matching (MWSM) problems is precisely the subclass of problems \eqref{problem} where the second-order differentials of the objective function $f$ are all zero, see Theorem~\ref{thm:linear-to-second-order}. This implies in particular that the class of minimum cut representable problems is significantly larger than the class of linear objective problems (MWSM).

\vspace{3mm}

{\noindent\bf Applications.} Our results can be employed for applications as follows. First, they give a tool to recognize when a problem is minimum cut representable and %\textemdash when it is difficult to construct a corresponding witness digraph by hand\textemdash 
a systematic way to construct a witness digraph. Second, when a problem is not minimum cut representable, our results give an explicit certificate. Third, our results can be used to reveal when an objective function can be represented as a linear function on the pairs of students and schools. To illustrate these use cases, we introduce and study the following problems motivated by real-world applications.

\vspace{2mm}
{\noindent\em Matching Siblings.} Our first set of applications addresses school choice when siblings wish to enroll in the same, or in related, programs. %The goal is to compute a stable matching of students to schools ensuring that specific pairs of students are matched to the same or related schools. 
In particular, all our problems start from a school matching instance $(A,B,\AgentGreater{}{}, q)$  where some pairs of students $A$ are marked as siblings.

In the first scenario, we address the problem of deciding whether there exists a stable matching of $(A,B,\AgentGreater{}{}, q)$ where each pair of siblings are matched to the same school (MSSS). Using our framework, specifically Theorem~\ref{thm:main1} and Theorem~\ref{thm:linear-to-second-order}, we show that (MSSS) can be solved in polynomial time by casting it as an instance of the Minimum Weight Stable Matching problem (MWSM). See Theorem~\ref{thm:MSSS}. 

In the second scenario, we assume that each pair of siblings wish to be enrolled in the same activity (e.g., school or after-school classes such soccer, tennis, music, ...). However, various constraints (e.g. age, skill-level, scheduling restrictions, ...) may require each sibling to register in a different section for the same activity. Formally, we partition schools in $B$ into activities, which are further divided into classes. Each student prefers at most one class from each activity over the outside option. The problem is then to decide whether there exists a stable matching of $(A,B,\AgentGreater{}{},q)$ where each pair of siblings is enrolled in classes corresponding to the same activity (MSDP). We show that, in general, (MSDP) is NP-Hard (see Theorem~\ref{thm:MSDP-Np-Complete}) and, in some well-defined sense, not minimum cut representable (see Lemma~\ref{lem:counter-ex-not-MSDP}). However,  when the preference lists of each pair of siblings induce the same order of activities, which is the third scenario we consider denoted by (MSSP), the problem becomes minimum cut representable and can be solved in polynomial time (see Theorem~\ref{thm:after-school-siblings}). Unlike (MSSS), we show that (MSSP) cannot be linearized and is not an instance of (MWSM) (see Lemma~\ref{lem:counter-ex-not-MWSS}).

\vspace{3mm}
{\noindent\em Two-Stage Stochastic Stable Matchings (2STO)\footnote{A preliminary version of this paper, in which we introduce and study the two-stage stochastic stable matching problem (2STO), was published in Integer Programming and Combinatorial Optimization (IPCO) 2024~\citep{faenza2024two}. This full version significantly extends the IPCO 2024 extended abstract to include: the general theory of minimum cut representability, its applications to the matching siblings problem, and missing proofs on the two-stage problem (2STO).}.}
While the classical stable matching problem assumes static and fully known input, in many applications the input changes over time as new agents enter or leave the market. For instance, each year, and after an assignment of the students to San Francisco public middle schools has been decided, around 20\% of the students who were allotted a seat choose not to use it, mostly to join private schools instead~\citep{Irene}. Also, new students move to the city and need to be allotted a seat. On the other side of the market, the schools may have unforeseen budget cuts or expansions, leading to a change in the number of available seats.  In all of these scenarios, a second-stage reallocation of at least part of the seats is required. For instance, while a subset of students might commit to their assigned seats in the first round of matching (and hence leave the market along with their allotted seats), the original assignment needs to be adjusted for the rest of the students after new agents arrive (or old agents depart) in order to maintain stability. It is therefore important that the initial stable matching, besides being of good quality, is adaptable to a changing environment with small adjustments so that it leads to a good quality second round matching and to small dissatisfaction of the students present in both rounds from being downgraded to a less preferred school. It is easy to find examples (see Figure~\ref{fig:running-intro}) where, in order to achieve such goal, one might want to go beyond the student optimal stable matching\footnote{The student optimal matching is a stable matching that all students prefer to every other stable matching, see Definition~\ref{def:student-optimal-pessimal}.}, which is the obvious solution that is almost always employed (see, e.g.,~\cite{abdulkadirouglu2005new,roth1999redesign}).

\begin{figure}[t!]
\begin{center}
\begin{tabular}{c  l c | c c l }
$a^1$: & \fbox{$b^1$}$\AgentGreater{a^1}{} b^2  \AgentGreater{a^1}{} b^3  \AgentGreater{a^1}{} \underline{b^4}   \AgentGreater{a^1}{} b^5 \AgentGreater{a^1}{} \emptyset$ & \hspace{.15cm} &  & $b^1$: & $a^4 \AgentGreater{b^1}{} a^5 \AgentGreater{b^1}{} a^3  \AgentGreater{b^1}{} a^2  \AgentGreater{b^1}{} a^1 \AgentGreater{b^1}{} \emptyset$ \\ 
$a^2$: &  \fbox{$b^2$}$\AgentGreater{a^2}{} b^1  \AgentGreater{a^2}{} b^4  \AgentGreater{a^2}{} \underline{b^3} \AgentGreater{a^2}{} b^5 \AgentGreater{a^2}{} \emptyset$ & \hspace{.15cm} & & $b^2$: & $a^3 \AgentGreater{b^2}{} a^5 \AgentGreater{b^2}{} a^4  \AgentGreater{b^2}{} a^1  \AgentGreater{b^2}{} a^2 \AgentGreater{b^2}{} \emptyset$ \\
$a^3$: &  \fbox{$b^3$}$\AgentGreater{a^3}{} b^4  \AgentGreater{a^3}{} b^1  \AgentGreater{a^3}{} \underline{b^2} \AgentGreater{a^3}{} b^5 \AgentGreater{a^3}{} \emptyset$  & \hspace{.15cm} & & $b^3$: & $a^2 \AgentGreater{b^3}{} a^5 \AgentGreater{b^3}{} a^1  \AgentGreater{b^3}{} a^4  \AgentGreater{b^2}{} a^3 \AgentGreater{b^3}{} \emptyset$ \\
$a^4$: &  \fbox{$b^4$}$\AgentGreater{a^4}{} b^3  \AgentGreater{a^4}{} b^2  \AgentGreater{a^4}{} \underline{b^1} \AgentGreater{a^3}{} b^5 \AgentGreater{a^4}{} \emptyset$  & \hspace{.15cm} & & $b^4$: & $a^1 \AgentGreater{b^4}{} a^5 \AgentGreater{b^4}{} a^2  \AgentGreater{b^4}{} a^3  \AgentGreater{b^4}{} a^4 \AgentGreater{b^4}{} \emptyset$ \\
$a^5$: & \fbox{$\underline{b^5}$}$\AgentGreater{a^5}{} b^3  \AgentGreater{a^5}{} b^2  \AgentGreater{a^5}{} b^1 \AgentGreater{a^5}{} b^4 \AgentGreater{a^5}{} \emptyset$  & \hspace{.15cm} & & $b^5$: & $a^5 \AgentGreater{b^5}{} a^1 \AgentGreater{b^5}{} a^2  \AgentGreater{b^5}{} a^3  \AgentGreater{b^5}{} a^4 \AgentGreater{b^5}{} \emptyset$
\vspace{3mm}
\end{tabular}
\end{center}
\caption{Example adapted from~\cite[Example 8.1]{faenza2022legal}. Consider the first-stage instance $I_1$ given above, where $a^1,\dots, a^5$ are the students and all schools have a quota of $1$. Then $M_0= \{a^1b^1, a^2b^2, a^3b^3, a^4b^4, a^5b^5\}$ and $\underline M = \{a^1b^4, a^2b^3, a^3b^2, a^4b^1, a^5b^5\}$ are two stable matchings of $I_1$, with $M_0$ being student-optimal (the partner of each student $a$ in $M_0$ and $\underline M$ is, respectively, boxed and underlined in $a$'s preference list). If $b^5$ leaves the market in the second-stage, the only stable matching is $M'= \underline M\setminus \{a^5b^5\}$, hence, $\underline M$ minimizes the number of students downgraded in the second-stage.
}\label{fig:running-intro}
\end{figure}

Motivated by these considerations, we introduce and study a two-stage stochastic stable matching problem, where agents randomly enter and leave the market between the two stages. The goal of the decision maker is, roughly speaking, to maximize the expected quality of the matchings across the two stages and minimize the expected students' discontent from being downgraded to a less preferred school when going from the first to the second-stage. See Section~\ref{sec:2-stage} for formal definitions.

We use our framework to show that the two-stage stochastic stable matching problem can be solved in polynomial time when the distribution of the entering/leaving agents is given explicitly by the list of all possible scenarios and their probabilities (see Theorem~\ref{thm:two-stage-main1}). We then consider the more general case when the distribution is only given by a sampling oracle. We show that the problem becomes NP-Hard in this case (see Theorem~\ref{thm:imp-hard}), and give a sampling based pseudo-polynomial algorithm (see Theorem~\ref{thm:imp-fpras}). 

We evaluate the performance of our method on randomly generated instances. In particular, we compare the numerical performance of four first-stage matchings: The solution $M^*$ output by our algorithm, the one-side optimal matchings, and the best stable matching in hindsight denoted by $M^{\sf off}$, which is the first-stage stable matching that would have been optimal to pick if one knew the realization of the second-stage distribution. One of the key insights we observe is that for a wide range of instances, the stable matching computed by our algorithm significantly outperforms the one-sided optimal matchings commonly used in practice, indicating that when minimizing the rank change of students is a priority, our solution provides a more effective approach.

\subsection{Additional Related Work}\label{sec:additional-work}

There is substantial body of literature focusing on finding stable matchings that optimize specific objectives and/or satisfy additional constraints beyond stability. Key optimization objectives explored in the literature include regret (\cite{mandal2021matchings,knuth1976marriages,eirinakis2013finding}), social welfare (\cite{jain2024maximizing,anshelevich2009anarchy}), egalitarian measures (\cite{irving1987efficient,mcvitie1971stable}), fairness (\cite{teo1998geometry,sethuraman2006many}), and balance (\cite{feder1995stable,gupta2021balanced,narang2020achieving}). As for constraints, the literature has considered various types, including forced or forbidden pairs (\cite{boehmer2024adapting,mandal2021matchings}), lower quotas (\cite{biro2010college,boehmer2022fine}), regional constraints (\cite{kamada2017stability,kamada2015efficient,goto2016strategyproof}), diversity constraints (\cite{ehlers2014school,echenique2015control,kojima2018designing}), and budgetary constraints (\cite{ismaili2019weighted}). Our work develops a framework to efficiently solve a wide range of constrained optimization problems over stable matchings via reduction to minimum cut problems over rotations. A different approach to stable matching problems with additional, complex constraints is to relax the stability requirement in a controlled manner, as to output a solution that satisfies the additional constraints and is, in a well-defined sense, close to be stable (see, e.g.,~\cite{nguyen2021stability,nguyen2018near,nguyen2019stable}). It is however often not known if the corresponding solutions can be found efficiently.

A widely studied setting close to our first set of applications is when some applicant students already have siblings enrolled in the schools, in which case classical stable matchings are computed where schools prioritize students with siblings already attending, see, e.g., ~\cite{correa2022school,DOE}. Also closely related is the the classical problem of stable matching with couples (\cite{biro2011stable,klaus2005stable,klaus2007some,kojima2013matching,mcdermid2010keeping}). Unlike our applications where siblings have individual preferences over schools, couples in the stable matching with couples problem have joint preferences making the problem significantly more challenging. More broadly, our applications are related to the literature on stable matchings with externalities, where agents are concerned not only with their own partners but also with the partners of others (\cite{bodine2011peer,dutta1997stability,hafalir2008stability,pycia2007many,revilla2007many}).

As for our second application, there is a growing body of literature addressing two-stage models for the school matching problem. Some of these studies assume a fixed first-stage stable matching and aim to compute a second-stage stable matching that adapts to environmental changes, such as agents entering or leaving the market (\cite{fixedfirststage}) and shifts in preferences (\cite{compte2008voluntary,bredereck2020adapting}). Other research focuses on computing both first- and second-stage stable matchings with desirable properties, including scenarios involving market entry/exit (\cite{bampis2023online}) or changing preferences (\cite{khuller1994line}). The work by \cite{mentzer2023student} and \cite{feigenbaum2020dynamic} explores the effects of overbooking in the first-stage on the adjustments needed in the second-stage, while \cite{genc2017finding} and \cite{chen2021matchings} investigate computing stable matchings that are robust to changes.

Also related to our work is the work on representability of set functions over the Boolean lattice as cut functions (e.g.,~\cite{cunningham1985minimum,fujishige2001realization}). These works give necessary and sufficient conditions on functions defined over the Boolean lattice to be cut functions of directed capacitated graphs. Another relevant stand of literature is two-stage stochastic optimization. This framework has been applied to a variety of combinatorial problems including bipartite matching (\cite{feng2021two,katriel2008commitment}), vehicle routing (\cite{verweij2003sample}), scheduling (\cite{shmoys2007approximation}), Steiner tree problems (\cite{bomze2010solving}), vertex cover (\cite{gupta2004boosted}), and facility location (\cite{koca2023two}), among others. In this context, the Sample Average Approximation (SAA) technique is a widely used method to approach these problems. The method involves approximating the expected value in the second-stage by a sample average over a set of i.i.d. samples, then solving the resulting sample average problem (\cite{charikar2005sampling,ravi2006hedging,swamy2004sample,swamy2005sampling}).

\section{Preliminaries}\label{sec:prelim}

Let us begin by recalling some known notions and facts about partially ordered sets and lattices.

\vspace{3mm}
{\noindent \bf Posets, Lattices, and Upper-Closed Sets.\bf } We begin with the definition of a partially ordered set.

\begin{definition}[Partially Ordered Set (Poset)]
    A \emph{Partially Ordered Set} (\emph{Poset}) is given by a pair $(S, \geq)$, where $S$ is a set and $\geq$ is a partial order over the elements of $S$. When the partial order is clear from the context we call $S$ itself a poset.
\end{definition}

Lattices are special cases of posets defined as follows.

\begin{definition}[Meet, Join, Lattice]
    Let $(\cL, \geq)$ be poset and $x,y,z \in \cL$. If $z \geq x$, $z \geq y$, and for every $w \in \cL$ with $w \geq x,y$ it holds that $w \geq z$, we say that $z$ is the \emph{join} of $x,y$. If $z \leq x$, $z \leq y$, and for every $w \in \cL$ with $w\leq x,y$ we have $w \leq z$, we say that $w$ is the \emph{meet} of $x,y$. If every pair of elements of $\cL$ has a meet and a join, we say that $(\cL,\geq)$ is a \emph{lattice}.
\end{definition}

For a lattice $(\cL, \geq)$, we can define the \emph{join} (resp., \emph{meet}) \emph{operator} as the map that takes in input two elements of $\cL$ and outputs their join (resp., meet). An example of a lattice is the lattice of upper-closed subsets of a poset, defined as follows.

\begin{definition}[Lattice of Upper-Closed Subsets of a Poset]
    Let $(S, \geq)$ be a poset. We say that $X\subseteq S$ is \emph{upper-closed} if $x \in X, y \geq x$ implies $y \in X$. Let $\cU$ denote the set of all upper-closed subsets of $(S, \geq)$. The lattice of upper-closed subsets of $(S, \geq)$ is the poset $(\cU, \geq')$ where the partial order $\geq'$ is such that $A \leq' B$ if and only if $B \subset A$. The meet and join of this lattice are given by the set union and set intersection respectively.
\end{definition}

We next define the notion of lattice isomorphism.

\begin{definition}[Lattice Isomorphism]
    Given two lattices $(\cL^1, \geq^1)$ and $(\cL^2, \geq^2)$ whose meets  are denoted by $\wedge^1$ and $\wedge^2$ respectively and whose joins  are denoted by $\vee^1$ and $\vee^2$ respectively, we say that $\Phi: \cL^1 \rightarrow \cL^2$ is a lattice isomorphism if $\Phi$ is bijective and for every $x, y \in \cL^1$, it holds that $\Phi(x \wedge^1 y) = \Phi(x) \wedge^2 \Phi(y)$ and $\Phi(x \vee^1 y) = \Phi(x) \vee^2 \Phi(y)$.
\end{definition}

Finally, a sublattice is a subset of a lattice defined as follows,

\begin{definition}[Sublattice]
    Given a lattice $(\cL, \geq)$ with meet $\wedge$ and join $\vee$, we say that $\cF \subset \cL$ is a sublattice of $(\cL, \geq)$ if  for every $x, y \in \cF$, it holds that $x \vee y, x \wedge y \in \cF$.
\end{definition}

We now recall some classical notions and facts about stable matchings. 
In the following, we fix a school matching instance $I=(A,B,\succ, q)$ over students $A$ and schools $B$. Let $n = |A| + |B|$ denote the total number of agents, and let $\cS(I)$ denote the set of stable matchings of $I$.

\vspace{3mm}
{\noindent \bf The Lattice of Stable Matchings.} Many algorithmic results on stable matchings rely on the following lattice structure that can be associated to it. The proof of the following result can be found in~\cite{gusfieldbook}.

\begin{lemma}[Lattice of Stable Matchings]
\label{lem:lattice-stable-matchings}
    Consider the partial order $\geq$ over stable matchings where, for $M, M' \in {\cal S}(I)$, $M \geq M'$ if and only if $M(a) \succeq_a M'(a)$ for every $a \in A$. Then $(\cS(I), \geq)$ is a lattice, where the join (resp., the meet) of $M$ and $M'$, denoted by $M \vee M'$ (resp., $M \wedge M'$) is the stable matching where every $a \in A$ is paired with the partner they prefer (resp., dislike) the most between $M(a)$ and $M'(a)$, and every $b \in B$ with strictly less than $q_b$ partners in the new matching is also paired with the outside option.
\end{lemma}

Lemma~\ref{lem:lattice-stable-matchings} implies the existence of a largest and smallest stable matching with respect to $\geq$. These are classically known as the student optimal and student pessimal stable matchings and are defined as follows.

\begin{definition}[Student Optimal and Student Pessimal Stable Matchings] 
\label{def:student-optimal-pessimal}
The \emph{student-optimal} (resp., \emph{-pessimal}) stable matching denoted by $M_0$ (resp., $M_z$) are the elements of ${\cal S}(I)$ such that, for each $M \in \cS(I)$, it holds that $M_0 \geq M \geq M_z$.
\end{definition}

\vspace{3mm}

{\noindent \bf The Poset of Rotations.} For $M,M'\in{\cal S}(I)$, we say that $M'$ is an \emph{immediate predecessor} of $M$ if $M' < M$ and there is no $M''\in \cS(I)$ such that $M'< M''< M$. We now define the notion of a {\em rotation}.

\begin{definition}[Rotations]\label{def:rotation-onemany}
    Let $M,M' \in {\cal S}(I)$ with $M'$ an immediate predecessor of $M$. The pair $\rho(M, M') = (M\backslash M', M' \backslash M)$ is referred to as a {\em rotation}.
    We denote by $$\cR(I)=\{\rho(M,M'):\textrm{$M'$ is an immediate predecessor of $M$ in $({\cal S}(I), \geq$)}\}$$ the set of all rotations.
\end{definition}

Given a rotation $\rho = \rho(M, M')$, we denote $\rho_+ = M\backslash M'$ (resp. $\rho_- = M'\backslash M$) the set of pairs that are added (resp., removed) when going from $M'$ to $M$. We note that $\rho_+$ and $\rho_-$ are pairs over the same set of students and schools. In fact, a well known result in the theory of stable matchings is the so-called \emph{Rural Hospitals Theorem} (see e.g., \cite[Theorem 1.6.3]{gusfieldbook}), which states that, if $M$ is a stable matching such that, for some school $b \in B$, $|M(b)|<q_b$, then $M'(b)=M(b)$ for every stable matching $M'$, and moreover the set of students assigned to the outside option is the same in all stable matchings. Therefore, the transition between a stable matching $M$ and an immediate predecessor $M'$ consist of a subset of students exchanging their respective assigned school seats. This also justifies the naming {\em rotation} for the object $\rho=(\rho_+, \rho_-)$. 

Although every pair of consecutive matchings in the stable matching lattice defines a rotation, many rotations coincide. Formally, the following holds.

\begin{lemma}[Size of ${\cal R}(I)$]
\label{lem:sizeofrotationposet}
$|{\cal R}(I)|\leq n^2$.
\end{lemma} 

Lemma~\ref{lem:sizeofrotationposet} is a direct consequence of~\cite[Lemma 6]{bansal2007polynomial}. We now introduce the concept of elimination of a rotation.

\begin{definition}[Elimination of a rotation]\label{def:rotation-elimination-onemany}
    For $\rho\in\cR(I)$ and $M\in{\cal S}(I)$, we say $\rho$ is exposed in $M$ if $\rho=\rho(M, M')$ for some immediate predecessor $M'$ of $M$ that we denote by $M/\rho$. We say that $M/\rho$ is the matching obtained by \emph{eliminating} $\rho$ from $M$.
\end{definition}

The following lemma shows that one can traverse the stable matchings lattice by successive elimination of rotations. The proof can be found in~\cite{bansal2007polynomial}.

\begin{lemma}[Traversal of Stable Matchings Using Rotations]\label{lem:matching-transversal}

Recall that $M_0$ denotes the student-optimal stable matching. Then, for any $M\in{\cal S}(I)$, there exists a sequence $M_0,M_1,\dots,M_k$ such that $M_k=M$ and $M_{i}$ is an immediate predecessor of $M_{i-1}$, and hence, $M_{i}=M_{i-1}/\rho(M_{i-1},M_{i})$, for every $i\in[k]$. We write $M=M_0/R$ where $R=\{\rho(M_{i-1},M_{i}),i=1,\dots,k\}$. Moreover, for a given $M \in \cS(I)$, the set of rotations $R$ is unique and independent of the choice of the sequence $M_0,\dots,M_k$. 
\end{lemma}

The representation from the above lemma allows us to define a poset over the set of rotations. 

\begin{definition}[Poset of Rotations]\label{def:rotation-poset-onemany}
    The poset of rotations $(\cR(I), \trianglerighteq)$ is defined as follows: for all $\rho, \rho' \in \cR(I)$, we have $\rho \trianglerighteq \rho'$ if and only if either $\rho = \rho'$, or whenever $ \rho' $ is exposed in some $ M \in \cS(I) $ such that $ M=M_0/R $ then $\rho \in R$. In the latter case, we also write $\rho\triangleright\rho'$.
\end{definition}

Basically, the above definition states that a rotation $\rho$ is larger than a rotation $\rho'$ with respect to the order $\triangleright$ when in order to get to a stable matching where $\rho'$ is exposed starting from $M_0$ (and along any sequence of immediate predecessors) one needs first to eliminate $\rho$ along the way. Note that under this partial order, the sets $R$ given by Lemma~\ref{lem:matching-transversal} must be upper-closed in $(\cR(I), \trianglerighteq)$. In fact, a stronger result holds, namely, the existence of a lattice isomorphism between the lattice of upper-closed subsets of the poset $(\cR(I), \trianglerighteq)$ and the  lattice of stable matchings $(\cal S(I),\geq)$, as stated in the following theorem.
\begin{theorem}[Representation with Rotations]
\label{thm:representation_with_rotations} 
The lattice of upper-closed subsets of the poset $(\cR(I), \trianglerighteq)$ is isomorphic to the lattice of stable matchings $(\cS(I), \geq)$. This isomorphism is such that: (i) given a stable matching $M$ such that $M = M_0/R$, the set $R$ is upper-closed and is the image of $M$ by the isomorphism. We denote such set by $R_M$. (ii) given an upper-closed set of rotations $R$, there exists a stable matching, denoted by $M_R$, such that $M_R=M_0/R$, which is the image of $R$ by the isomorphism. %We denote such stable matching by $M_R$.
\end{theorem}

Theorem 1 follows from~\cite[Theorem 3]{bansal2007polynomial}. In particular,~\cite[Theorem 3]{bansal2007polynomial} shows a one-to-one correspondence, then it is immediate from the definition of the correspondence that the stable matching corresponding to the union (resp. intersection) of two upper-closed subsets of rotations is exactly the meet (resp. join) of the two stable matchings, implying that the one-to-one correspondence is also a lattice isomorphism.

An immediate consequence of Theorem~\ref{thm:representation_with_rotations} is that one needs to eliminate all rotations starting from $M_0$ to get to the student-pessimal matching $M_z$. In particular, the following corollary holds.

\begin{corollary}\label{cor:all-rotations}
    $M_z=M_0/\cR(I)$.
\end{corollary}

\vspace{3mm}
{\noindent \bf Sublattices of Stable Matchings and the Poset of Meta-Rotations.}
For a sublattice $\cF$ of the stable matchings lattice $(\cS, \geq)$, and similarly to the lattice isomorphism between the lattice of all stable matchings $\cS(I)$ and the lattice of upper-closed subsets of rotations, a lattice isomorphism can be defined between the stable matchings of $\cF$ and the lattice of upper-closed subsets of the so-called poset of {\em proper meta-rotations}. The poset of {\em proper meta-rotations} already appeared for example in \cite{scott2005study,metarot} for the one-to-one setting (i.e., $q_b = 1$ for every school $b \in B$) where {\em meta-rotations} were defined as sets of rotations encountered between consecutive stable matchings of $\cF$. We give here a slightly different definition that is more convenient for our purposes (and which reduces to the definition of \cite{scott2005study,metarot} in the one-to-one setting) which is based on grouping rotations into equivalent classes of a properly defined equivalence relation. The proofs of the results concerning meta-rotations are given in the Appendix for completeness.

\begin{definition}[Meta-Rotations] Consider a sublattice $\cF$ of the lattice of stable matchings $(\cS(I), \geq)$. Consider the binary relation $\sim$ on the set of rotations such that $\rho \sim \rho'$ if and only if $\rho$ and $\rho'$ are always eliminated together in any stable matching of $\cF$, i.e., if and only if for every stable matching $M \in \cF$ either $\rho, \rho' \in R_M$ or $\rho, \rho' \notin R_M$. This is clearly an equivalence relation. The equivalence classes of $\sim$ partition the set of rotations $\cR(I)$ into $k+2$ groups of rotations $\theta^{\cF}_0, \theta^{\cF}_1, \dots, \theta^{\cF}_k, \theta^{\cF}_{z}$ such that: $\theta^{\cF}_0 = \{\rho \;|\; \forall M \in \cF, \rho \in R_M\}$, $\theta^{\cF}_{z} = \{\rho \;|\; \forall M \in \cF, \rho \notin R_M\}$.
We refer to $\theta^{\cF}_0, \dots, \theta^{\cF}_{z}$ as \emph{meta-rotations}. In particular, we refer to $\theta^{\cF}_0$ and $\theta^{\cF}_{z}$ as the \emph{trivial} meta-rotations and to $\theta^{\cF}_1, \dots, \theta^{\cF}_{k}$ as the \emph{proper} meta-rotations.
\end{definition}
Let $\cR^{\cF}(I)$ (resp. $\cR^{\cF}_{\sf prop}(I)$) denote the set of all (resp. proper) meta-rotations corresponding to a sublattice $\cF$ of ${\cal S}(I)$. We define the following partial order over the proper meta-rotations.

\begin{definition}[Poset of Proper Meta-Rotations]
\label{def:poset_metarotations}
    The poset of proper meta-rotations denoted by $(\cR_{\sf prop}^{\cF}(I), \trianglerighteq^{\cF})$ is such that $\forall \theta, \theta' \in \cR_{\sf prop}^{\cF}(I), \;\; \theta \trianglerighteq^{\cF} \theta'$ if and only if for each $M \in \cF$ such that $\theta'\subseteq R_M$, we have $\theta\subseteq R_M$. That is, whenever the rotations of $\theta'$ are eliminated in a stable matching $M \in \cF$, then the rotations $\theta$ are also eliminated in $M$.
\end{definition}

It can be shown that the poset from Definition~\ref{def:poset_metarotations} is well-defined.
\begin{lemma}
\label{lem:meta_rotations_is_poset}
     $(\cR_{\sf prop}^{\cF}(I), \trianglerighteq^{\cF})$ from Definition~\ref{def:poset_metarotations} is a poset.
\end{lemma}
Similar to the isomorphism between the lattice of upper-closed subsets of rotations and stable matchings, there exists an isomorphism between the lattice of upper-closed subsets of proper meta-rotations and the stable matchings of $\cF$. 
In particular, the following theorem holds.
\begin{theorem}[Representation with Proper Meta-Rotations]
\label{thm:representation_meta_rotations}
The lattice of upper-closed subsets of proper meta-rotations is isomorphic to the sublattice $\cF$. This isomorphism is such that: (i) given a stable matching $M \in \cF$, the set $R_M$ can be partitioned into meta-rotations $\theta_0^{\cF}, \theta_1, \dots, \theta_r$ such that $\{\theta_1, \dots, \theta_r\}$ is an upper-closed subset of proper meta-rotations $(\cR_{\sf prop}^{\cF}(I), \trianglerighteq^{\cF})$. The image of $M$ by the isomorphism is the set $\{\theta_1, \dots, \theta_r\}$ that we denote by $\Theta_M$. (ii) given an upper-closed subset of proper meta-rotations $\Theta$, the subset of rotations $\theta_0^{\cF} \cup \bigcup_{\theta \in \Theta} \theta$ is upper-closed and the corresponding stable matching $M_{\{\theta_0^{\cF} \cup \bigcup_{\theta \in \Theta} \theta\}}$ is in $\cF$. The image of $\Theta$ via the isomorphism is the stable matching $M_{\{\theta_0^{\cF} \cup \bigcup_{\theta \in \Theta} \theta\}}$ that we denote by $M_\Theta$.
\end{theorem}

\section{Minimum Cut Representability}

\label{sec:minimum cut-representability}
In this section, we formally define minimum cut representability and give necessary and sufficient conditions for an optimization problem to be minimum cut representable. 
In the following, we fix a school matching instance $I=(A,B,\succ,q)$ over students $A$ and schools $B$. Unless otherwise stated, proofs can be found in the appendix.

\subsection{Definitions}
Let us begin by formally defining minimum cut representability of an optimization problem over stable matchings.

\begin{definition}[Minimum Cut Representability]\label{def:mcs}
Consider a feasibility set $\cF \subset \cS(I)$ and an objective function $f: \cF \rightarrow \mathbb{R}$. We say that the optimization problem \eqref{problem} is \emph{minimum cut representable} if and only if there exists a capacitated digraph $\cD(\cR(I) \cup \{s,t\}, A)$ with non-negative (possibly infinite) capacities and a real number $C$ such that,
\begin{enumerate}[(i)]
    \item For every $R \subset \cR(I)$, the $s$-$t$ cut $ \{s\} \cup R$ has finite capacity if and only if $R$ is an upper-closed set of the poset of rotations $(\cR(I), \trianglerighteq)$ and $M_R \in \cF$.
    \item The value of a finite capacity $s$-$t$ cut $ \{s\} \cup R$ is given by $f(M_R) + C$.
\end{enumerate}
We refer to $\cD(\cR(I) \cup \{s,t\}, A)$ as a \emph{witness digraph}. 
\end{definition}

One classical example of a minimum cut representable problem is the Minimum Weight Stable Matching problem (MWSM) (e.g. \cite{gusfieldbook}). We are interested in characterizing under which conditions  \eqref{problem} is minimum cut representable, and if so, construct a witness digraph.

Consider a real valued function $f:\cF \rightarrow \mathbb{R}$ over a sublattice $\cF$ of the lattice of stable matchings, and let $\theta^{\cF}_0, \theta^{\cF}_1, \dots, \theta^{\cF}_k, \theta^{\cF}_{z}$ be the corresponding meta-rotations. For every proper meta-rotation $\theta \in \cR_{\sf prop}^{\cF}(I)$, let $\Theta^\theta = \{\theta' \in \cR_{\sf prop}^{\cF}(I) \;|\; \theta' \trianglerighteq^{\cF} \theta\}$ be the set of all proper meta-rotations that are larger or equal than $\theta$, and $\Theta_\theta = \{\theta' \in \cR_{\sf prop}^{\cF}(I)\;|\; \theta' \triangleright^{\cF} \theta\}$ the set of all proper meta-rotations that are strictly larger than $\theta$. Note that both of these sets are upper-closed. Let $M^\theta = M_{\Theta^{\theta}}$ and $M_\theta = M_{\Theta_{\theta}}$, and let $R^\theta = R_{M^{\theta}}$ and $R_\theta = R_{M_{\theta}}$. Note that, by Theorem~\ref{thm:representation_meta_rotations}, $R_\theta = \theta_0^{\cF} \cup \bigcup_{\theta' \triangleright^{\cF} \theta} \theta'$ and $ R^\theta = \theta_0^{\cF} \cup \bigcup_{\theta' \trianglerighteq^{\cF} \theta} \theta' = R_\theta \cup \theta.$
In words, $R^\theta$ (resp. $R_{\theta}$) consists of all rotations contained in all proper meta-rotations larger (resp. strictly larger) than $\theta$, plus rotations in the trivial meta-rotation $\theta_0^{\cF}$ that are eliminated in every stable matching of $\cF$. We use the two special stable matchings $M_\theta$ and $M^\theta$ to define the notion of {\em first and second order differentials} of $f$. In the following, given a set $X$, we write $a \neq b \in X$ to denote that $a$ and $b$ are two distinct elements of $X$.

\begin{definition}[First Order Differential]\label{def:fod}
The first order differential of $f$ at a proper meta-rotation $\theta \in \cR_{\sf prop}^{\cF}(I)$ is defined as
$$
\partial f_{\theta} = f(M^{\theta})- f(M_{\theta}).
$$
\end{definition}

\begin{definition}[Second Order Differential]\label{def:sod}
The second order differential of $f$ at a pair of distinct proper meta-rotations $\theta \neq \theta' \in \cR_{\sf prop}^{\cF}(I)$, is defined as
$$\partial^2 f_{\theta, \theta'} = f(M^{\theta} \wedge M_{\theta'}) + f(M_{\theta} \wedge M^{\theta'}) - f(M_{\theta} \wedge M_{\theta'}) - f(M^{\theta} \wedge M^{\theta'}).
$$
\end{definition}

Before giving our characterization of minimum cut representability, we define a digraph associated to the pair $(f, \cF)$, that we refer to as the {\em Cut Digraph} of $(f, \cF)$. For every meta-rotation $\theta \in \cR^{\cF}(I)$, let $\widehat{\theta} \in \cR(I)$ denote some arbitrary (representative) rotation in $\theta$. The {\em Cut Digraph} of $(f, \cF)$ is defined as follows.

\begin{definition}[Cut Digraph]\label{def:Cut-Graph}
    Given a sublattice $\cF$ of $(\cS(I), \geq)$ and a function $f: \cF \rightarrow \mathbb{R}$, the cut digraph associated to the pair $(f,\cF)$ is the digraph $\cD^{f,{\cal F}}$ with vertex set $\cR(I) \cup \{s,t\}$ and arc set $A$ defined as follows:

    \vspace{2mm}
    {\em Arcs Within Meta-Rotations (WMR):} For every meta-rotation $\theta \in \cR^{\cF}(I)$, add an arc of infinite capacity between every pair of rotations $\rho, \rho' \in \theta$. 

    \vspace{2mm}
    {\em Arcs Between Proper Meta-Rotations (BPMR):} For every pair of distinct proper meta-rotations $\theta \neq \theta' \in \cR_{\sf prop}^{\cF}(I)$ such that $\theta \triangleleft^{\cF} \theta'$, add an arc of infinite capacity from $\widehat{\theta}$ to $\widehat{\theta'}$. Additionally, for every pair of distinct proper meta-rotations $\theta \neq \theta' \in \cR_{\sf prop}^{\cF}(I)$, add an arc of capacity $\frac{1}{2} \partial^2 f_{\theta, \theta'}$ from $\widehat{\theta}$ to $\widehat{\theta'}$.

    \vspace{2mm}
    {\em Arcs Between Source/Sink and Meta-Rotations (BSSMR):} First, add an arc of infinite capacity from $s$ to $\widehat{\theta^{\cF}_0}$ and an arc of infinite capacity from $\widehat{\theta^{\cF}_{z}}$ to $t$. Second, for every proper meta-rotation $\theta \in \cR_{\sf prop}^{\cF}(I)$, add an arc
    from $\widehat{\theta}$ to $t$
    of capacity $\partial f_\theta - \frac{1}{2}\sum_{\theta' \in \cR^{\cF}(I) \setminus \{\theta\}} \partial^2 f_{\theta, \theta'} + \gamma$ and an arc
    from $s$ to $\widehat{\theta}$
    of capacity $\gamma$. Here, $\gamma$ is any non-negative constant such that $\partial f_\theta - \frac{1}{2}\sum_{\theta' \in \cR^{\cF}(I) \setminus \{\theta\}} \partial^2 f_{\theta, \theta'} + \gamma \geq 0$ for every $\theta \in \cR_{\sf prop}^{\cF}(I)$. 
\end{definition}

\subsection{Necessary and Sufficient Conditions for Minimum Cut Representability}
We are now ready to state our main theorem, which gives necessary and sufficient conditions for minimum cut representability, and, in case of minimum cut representability, gives a witness digraph.

\begin{theorem}[Characterization of minimum cut representability]
\label{thm:main1}
The optimization problem \eqref{problem} is minimum cut representable if and only if
\begin{enumerate}[(i)]
    \item $\cF$ is a sublattice of $(\cS(I), \geq)$.
    \item For every $\theta \neq \theta' \in \cR_{\sf prop}^{\cF}(I)$, $$\partial^2 f_{\theta, \theta'} \geq 0.$$
    \item[(iii)] For every $M \in \cF$,
    $$f(M) = f(M_0^{\cF}) + \sum_{\theta \in \Theta_M} \partial f_\theta - \frac{1}{2} \sum_{\substack{\theta \neq \theta' \in \Theta_M}} \partial^2 f_{\theta, \theta'}.$$
\end{enumerate}
Moreover, when \eqref{problem} is minimum cut representable, the cut digraph $\cD^{f, \cF}$ associated to the pair $(f, \cF)$ is a witness digraph.
\end{theorem}

The function $f^{\sf apx}: M \rightarrow f(M_0^{\cF}) + \sum_{\theta \in \Theta_M} \partial f_\theta - \frac{1}{2} \sum_{\substack{\theta \neq \theta' \in \Theta_M}} \partial^2 f_{\theta, \theta'}$ in the right-hand-side of (iii) can be understood as a ``second order expansion'' of $f$. This motivates the naming ``first and second order differentials'' of the quantities $\partial f_\theta$ and $\partial^2 f_{\theta,\theta'}$. A further motivation is the fact that problems with a zero second order differential are exactly linear objective problems (MWSM) as we discuss next.

\smallskip 

{\noindent \bf Relation to the (MWSM) problem.} As mentioned before, every instance of the (MWSM) problem is minimum cut representable. The following theorem gives a characterization of when a minimum cut representable problem is linearizable, i.e., is an instance of (MWSM).

\begin{theorem}[Characterization of linearizable problems]
\label{thm:linear-to-second-order}
An optimization problem~\eqref{problem} is an instance of the minimum weight stable matching problem (MWSM) if and only if it is minimum cut representable and for every pair of distinct proper-meta rotations $\theta \neq \theta' \in \cR_{\sf prop}^{\cF}(I)$, it holds that $
\partial^2f_{\theta, \theta'} = 0
$.
\end{theorem}

\smallskip 

{\noindent \bf Additivity of minimum cut representable/linear problems.} The following lemma shows that any non-negative linear (i.e., conic) combination of minimum cut representable/linear problems is a minimum cut representable/linear problem. It will be used on several occasions in our applications. Its proof is immediate, thus omitted.

\begin{lemma}[Conic combinations of minimum cut representable/linear problems]\label{lem:additivity} 
Consider a sequence of minimum cut representable (resp., (MWSM)) problems  $(\Pi(f^1, \cF)),\dots, (\Pi(f^r,\cF))$ with witness digraphs $\cD^1(\cR(I)\cup \{s,t\} , A^1), \dots, \cD^K(\cR(I)\cup \{s,t\} , A^{I_k})$ and arc capacities $u_1, \dots, u_K$ respectively. Let $\lambda_1, \dots, \lambda_n \geq 0$. Then the problem $(\Pi(\sum_{k=1}^K \lambda_i f^i,\cF))$ is minimum cut representable (resp., a (MWSM)) problem with a witness digraph $\cD(\cR(I)\cup \{s,t\} , \cup_{k=1}^K A^{I_k})$ and arc capacities $\sum_{i=1}^n \lambda_i u_i$.
\end{lemma}

\smallskip

{\noindent \bf A Necessary Condition for Minimum Cut Representability.} 
The following necessary condition is implied by properties of rotations associated to meets and joins of matchings and of the cut function of a digraph. It will be employed to show that the certain problems are not minimum cut representable. 
\begin{lemma}\label{lem:if-mincut-then-lattice}
If the optimization problem \eqref{problem} is minimum cut representable, then the minima of $f$ over ${\cal F}$ form a sublattice of the lattice of stable matchings $(\cS(I), \geq)$.    
\end{lemma}

\smallskip

{\noindent \bf Algorithmic Implications of Our Framework.} Consider an optimization problem \eqref{problem} over a school matching instance $I$. %over stable matchings of $I$. 
To apply our framework, we first need to check that $\cF$ is a sublattice. 
Then, we compute the quantities  $f(M_0^{\cF})$, $\{\partial f_{\theta}\}_{\theta \in \cR_{\sf prop}^{\cF}(I)}$, $\{\partial^2 f_{\theta, \theta'}\}_{\theta \neq \theta' \in \cR_{\sf prop}^{\cF}(I)}$ and check the non-negativity of the latter. Next, we check that  $f(M) = f(M_0^{\cF}) + \sum_{\theta \in \Theta_M} \partial f_\theta - \frac{1}{2} \sum_{\theta \neq \theta' \in \Theta_M} \partial^2 f_{\theta, \theta'}$ for every $M \in \cF$.
The final step is to construct the witness digraph $\cD^{f, \cF}$ and solve a minimum $s$-$t$ cut problem over this digraph. 

While checking that $\cF$ is a sublattice and that $f(M) = f(M_0^{\cF}) + \sum_{\theta \in \Theta_M} \partial f_\theta - \frac{1}{2} \sum_{\theta \neq \theta' \in \Theta_M} \partial^2 f_{\theta, \theta'}$ for every $M \in \cF$ depend on the set of all stable matchings, which can be exponentially large (in the number of agents), it is often possible to check the properties implicitly (see for example our applications in Section~\ref{sec:applications}).  
One useful analogy here is that of submodular minimization. In particular, given a set function $g$, it is usually possible to check if $g$ is submodular efficiently without explicitly checking the submodularity equation for all subsets. 
It is also usually possible to compute the quantities $f(M_0^{\cF})$, $\{\partial f_{\theta}\}_{\theta \in \cR_{\sf prop}^{\cF}(I)}$, $\{\partial^2 f_{\theta, \theta'}\}_{\theta \neq \theta' \in \cR_{\sf prop}^{\cF}(I)}$ and the digraph $\cD^{f, \cF}$ as explicit polynomially sized functions of the input (for examples, see again our applications in Section~\ref{sec:applications}).

\section{Applications}\label{sec:applications}

To demonstrate the practical relevance of our framework, we introduce and study several real-world problems that motivated our work and whose solutions involve minimum cut representable problems. We start in Section~\ref{sec:siblings} with three problems motivated by the presence, among the students, of siblings who want to be matched to the same, or related, schools. We then discuss in Section~\ref{sec:2-stage} a two-stage stochastic stable matching problem where agents randomly enter and leave the market between the two stages.

\subsection{Matching Siblings}\label{sec:siblings}

Our first set of applications are problems that arise when siblings wish to enroll in the same, or in related, programs. In particular, we investigate the following three problems:
\begin{itemize}
\item In Section \ref{sec:matching-siblings-1}, given a set of pairs of siblings, we consider the problem of deciding on the existence of a stable matching assigning each pair of siblings to the same school. We refer to this problem as (MSSS). Through Theorem \ref{thm:linear-to-second-order}, we show that (MSSS) can be reduced to solving an instance of the minimum weight stable matching problem (MWSM) whose witness digraph can be constructed efficiently, and in particular, can be solved in polynomial time. Missing proofs can be found in the Appendix~\ref{ec:proofs-siblings-1}.
\item In Section~\ref{sec:matching-siblings-2}, we consider a setting where each pair of siblings have the same strict preference list over a series of after-school activities, but various constraints (e.g. age, skill-level, scheduling restrictions) may require each sibling to register in different sections for the same activity. We show that the problem of deciding whether there exists a stable matching where each pair of siblings can be assigned to (possibly different sections of) the same activity, referred to as (MSSP), can be reduced to solving a (non-linearizable) minimum cut representable problem whose witness digraph can be constructed efficiently, and in particular, can be solved in polynomial time. Missing proofs can be found in the Appendix~\ref{ec:MSSP-proofs}.
\item In Section~\ref{sec:matching-siblings-3}, we consider the generalization of the problem from Section~\ref{sec:matching-siblings-2}, where siblings may have different preference lists over activities. We refer to this problem as (MSDP). We discuss why our framework does not apply to (MSDP) and show that (MSDP) is in fact NP-Hard. Missing proofs can be found in the Appendix~\ref{app:MSDP-omitted}.
\end{itemize}

Throughout the section, we define a \emph{sibling instance} to be a school matching instance $I=(A,B,\succ,q)$, along with a set of sibling pairs $\cC = \{(a_1,\overline a_1), (a_2,\overline a_2),\dots, (a_m,\overline a_m)\}
\subseteq A \times A$ where for every $i \in [m]$, $(a_i, \overline{a}_i)$ is a pair of distinct students. We denote a sibling instance by $(I, \cC)$.
\subsubsection{Matching Siblings to the Same School (MSSS).}
\label{sec:matching-siblings-1}
The input of (MSSS) is a sibling instance $(I, \cC)$.
The goal is to compute a stable matching of $I$ where each pair of siblings are matched to the same school or show that no such matching exists. In this section, we show the following theorem.

\begin{theorem}\label{thm:MSSS}
    (MSSS) can be solved in polynomial time.
\end{theorem}

To solve (MSSS) one can equivalently solve \eqref{problem}, where for every stable matching $M \in \cS(I)$, $f(M)$ is the number of pairs of siblings that are not assigned to the same school in $M$, and $\cF$ is the set of all stable matchings. The optimal value of \eqref{problem} is zero if and only if the optimal stable matching sends each pair of siblings to the same school.

We now use the minimum cut representability framework to show the somehow surprising result that \eqref{problem} is an instance of the minimum weight stable matching problem (MWSM) whose witness digraph can be constructed efficiently, implying that (MSSS) can be solved in polynomial time. In the reminder, we assume without loss of generality the following.

\begin{assumption}\label{ass:5cases}
    Let $M_0,M_z$ be the student-optimal and student-pessimal matching, respectively. For every pair of siblings $(a,\bar{a}) \in \cC$, it holds that $M_0(a)\neq M_0(\bar{a})$ and $M_z(a)\neq M_z(\bar{a})$.
\end{assumption}

Assumption~\ref{ass:5cases} is without loss of generality in the sense that given an instance $\cI$ of (MSSS), we can efficiently construct a new instance $\cI'$ that verifies Assumption~\ref{ass:5cases} and such that  the maximum number of pairs of siblings matched to the same school in some stable matching of $\cI$ is the same over stable matchings of $\cI'$. The details are given in the Appendix~\ref{ec:reduction-5cases}.

Under Assumption~\ref{ass:5cases}, when a pair of siblings $(a,\bar{a})$ are matched to a school $b$ in some stable matching, there exists a couple of rotations $\rho_{in}(a,\bar{a},b)$ and $\rho_{out}(a,\bar{a},b)$ whose elimination assigns $a, \bar a$ together into $b$ and away from $b$ respectively. In particular, the following lemma holds.

\begin{lemma}\label{lem:2cases}
    Under Assumption~\ref{ass:5cases}, for every pair of siblings $(a,\bar{a}) \in \cC$ and school $b\in B$, exactly one of the following holds:
    \begin{enumerate}
        \item There is no stable matching $M\in\cS(I)$ such that $M(a)=M(\bar{a})=b$;
        \item There exists a unique pair of rotations $\rho_{in}(a,\bar{a},b),\rho_{out}(a,\bar{a},b)$ such that for every $M\in\cS(I)$, $M(a)=M(\bar{a})=b$ if and only if $\rho_{in}(a,\bar{a},b)\in R_M$ and $\rho_{out}(a,\bar{a},b)\notin R_M$. Moreover, $\rho_{in}(a,\bar{a},b) \triangleright \rho_{out}(a,\bar{a},b)$.
    \end{enumerate}
\end{lemma}

Based on Lemma~\ref{lem:2cases}, the set of rotations can be written as a union of three subsets of rotations, depending on whether the elimination of a rotation assigns $a,\overline a$ together into the same school, away from the same school, or none of the above. Formally, the rotations set can be written as $\cR(I)=\cR_{in}(a,\bar{a})\cup \cR_{out}(a,\bar{a})\cup\cR_{other}(a,\bar{a})$, where $\cR_{in}(a,\bar{a})=\{\rho_{in}(a,\bar{a},b):b\in B\textrm{ $s.t.$ $\rho_{in}(a,\bar{a},b)$ exists}\}$, $\cR_{out}(a,\bar{a})=\{\rho_{out}(a,\bar{a},b):b\in B\textrm{ $s.t.$ $\rho_{out}(a,\bar{a},b)$ exists}\}$ and $\cR_{other}(a,\bar{a})=\cR(I)\setminus(\cR_{in}(a,\bar{a})\cup \cR_{out}(a,\bar{a}))$. We remark that $\cR_{in}(a,\bar{a})$ and $\cR_{out}(a,\bar{a})$ may intersect.

We are now ready to show that our optimization problem~\eqref{problem} is indeed an instance of (MWSM) whose witness digraph can be constructed efficiently. We begin by decomposing the objective function $f$ as sum of simpler indicator functions. In particular, for each pair of siblings $(a_i, \overline{a}_i) \in \cC$, define the following indicator function,
\begin{align*}
    f^i(M)=\left\{\begin{array}{cc}
       0  & \quad\textrm{if $M(a_i)=M(\bar{a}_i)$} \\
       1  & \textrm{otherwise}
    \end{array}.\right.
\end{align*}
Observe that $f = \sum_{i\in[m]} f^i$.  We henceforth fix $i \in [m]$ and focus on the problem $(\Pi(f^i, \cF))$. Note that since $\cF = \cS(I)$, proper meta-rotations are simply rotations in our case. The following lemma gives the first and second order differentials of $f^i$ w.r.t. rotations and pairs of rotations, respectively. 

\begin{lemma}\label{lem:differentials-sibling}
The first and second order differentials for the indicator function $f^i$ are given by,
\begin{align*}
    \partial f^i_\rho=\left\{
    \begin{array}{cc}
      1,  & \textrm{if $\rho\in \cR_{out}(a_i,\bar{a}_i)$ and $\rho\notin \cR_{in}(a_i,\bar{a}_i)$} \\
      -1,   & \textrm{if $\rho\in \cR_{in}(a_i,\bar{a}_i)$ and $\rho\notin \cR_{out}(a_i,\bar{a}_i)$} \\
      0, & \textrm{otherwise}
    \end{array},\right.
\end{align*}
and
$$\partial^2 f^i_{\rho,\rho'}=0, \quad \forall \rho \neq \rho'\in\cR(I).$$
\end{lemma}

We have now all the ingredients to sketch a proof of Theorem~\ref{thm:MSSS}, with the full details given in the Appendix~\ref{ec:proofs-siblings-1-main}. In particular, we start by showing that the conditions of Theorem~\ref{thm:main1} hold for $(\Pi(f^i, \cF))$, implying the minimum cut representability of $(\Pi(f^i, \cF))$. Moreover, since by Lemma~\ref{lem:differentials-sibling} the second-order differentials of $f^i$ are all zero, by Theorem~\ref{thm:linear-to-second-order} the problem $(\Pi(f^i, \cF))$ is an instance of (MWSM). 
Then we show that the witness digraph $\cD^{f^i, \cF}$ can be constructed in time polynomial in the input size. Thus, by Lemma~\ref{lem:additivity}, \eqref{problem} is an instance of (MWSM) whose witness digraph can be constructed efficiently. Hence, (MSSS) can be solved in polynomial time. 

We remark that when there exists no stable matching where all pairs are matched to the same school, that is, when the instance of (MSSS) is infeasible, our approach allows us to efficiently find a stable matching that maximizes the number of pairs of siblings that are matched to the same school, since: (i) the maximum number of pairs of siblings assigned to the same school in some stable matching is the same in $\cI$ and in $\cI'$, the instance verifying Assumption~\ref{ass:5cases} that we construct in our proof of the assumption. (ii) for every $M \in {\cal S}(I)$, $f(M)$ constructed above counts exactly the number of pairs of siblings assigned to different schools by $M$.

\subsubsection{Matching Siblings with Same Preferences to the Same Activities (MSSP).}\label{sec:matching-siblings-2} 
In this section, we define and study a problem that arises, for instance, when siblings want to participate in the same activity (e.g., soccer, tennis, music, etc.), but there may be constraints (e.g., age, skill-level, scheduling restrictions) that require them to attend different classes where the same activity is performed. This model also pertains to a scenario where siblings want to register in different activities that are co-located.

Formally, an instance of (MSSP) is given by a sibling instance $(I,\cC)$ where $I=(A, B, \succ, q), \cC = \{(a_1,\overline a_1), (a_2,\overline a_2),\dots, (a_m,\overline a_m)\}\subseteq A \times A$ and such that additional conditions are satisfied, as detailed below. For consistency with the motivating application, in this section we will call elements of $B$ \emph{classes}. Let $Q=\{c_1,\dots, c_k\}$ be a set of $k$ activities then the following conditions holds:
\begin{itemize}
\item \emph{(Activities partition classes)} For every activity $j \in [k]$, there is a non-empty set $\{b_j^1,\dots, b_j^{\ell_j}\}\subseteq B$ of \emph{same-activity classes} (e.g., beginner class, intermediate class, advanced class, etc.), where moreover $\{b_j^1,\dots, b_j^{\ell_j}\}_{j \in [k]}$ is a partition of $B$. 
\item \emph{(Every student is eligible for exactly one same-activity class)} For each activity $c_j$ and each student $a$, there exists exactly one class $b_j^a \in \{b_j^1,\dots, b_j^{\ell_j}\}$ such that $b_j^a \succ_a \emptyset$. %This model the fact that $a$ is eligible for exactly one same-class  activity, for every activity. 
\item \emph{(Siblings rank activities in the same order)} Let $(a,\bar a) \in {\cal C}$ be a pair of siblings and $c_j, c_{\bar \jmath}$ be a pair of activities. Let $b^a_j$ be the unique class corresponding to activity $c_j$ for which $a$ is eligible, i.e., such that $b^a_j \succ_a \emptyset$. Define $b^{a}_{\bar \jmath}, b^{\bar a}_{j}, b^{\bar a}_{\bar \jmath}$ analogously. Then $b^a_j \succ_a b^a_{\bar \jmath}$ if and only if $b^{\bar a}_j \succ_a b^{\bar a}_{\bar \jmath}$.
\end{itemize}
$M \in \cS(I)$ is called \emph{activity-stable} if it matches each pair of siblings to same-activity classes. See Figure~\ref{fig:MSSP} for an example. The goal of the problem is to decide on whether an activity-stable matching exists. By building on Theorem~\ref{thm:main1}, we show the following.

\begin{theorem}\label{thm:after-school-siblings}
(MSSP) can be solved in polynomial time.
\end{theorem}

To show Theorem~\ref{thm:after-school-siblings}, we begin by defining, for each $i \in [m]$, the problem $\Pi_i$ of deciding whether there exists a matching $M \in \cS(I)$ where the sibling students $(a_i,\overline a_i)$ are matched to same-activity classes. The following lemma shows that the solution to problem $\Pi_i$ can be inferred from the solution of a minimum cut representable problem whose witness digraph can be constructed efficiently.
\begin{lemma}\label{lem:minimum cut-representable-after-school}
Fix $i \in [m]$. There exists an objective function $f^i: \cS(I) \rightarrow \mathbb{Z}_{\geq 0}$ such that the following holds:
\begin{enumerate}
    \item\label{it:lem:activity-siblings-function} For $M \in \cS(I)$, $f^i(M)=0$ if and only if, in $M$, $a_i$ and $\overline a_i$ are matched to same-activity classes.
    \item\label{it:lem:activity-siblings-polytime} $\Pi(f^i,\cS(I))$ is a minimum cut representable problem whose witness digraph can be constructed efficiently.
\end{enumerate}\end{lemma}

Now let $f:=\sum_{i \in [m]}f^i$. Using Lemma~\ref{lem:additivity} and Lemma~\ref{lem:minimum cut-representable-after-school}, part~\ref{it:lem:activity-siblings-polytime}, we know that $\Pi(f,\cS(I))$ can be solved in polynomial time. From Lemma~\ref{lem:minimum cut-representable-after-school}, part~\ref{it:lem:activity-siblings-function}, we know that the optimal solution of $\Pi(f,\cS(I))$ has objective function value $0$ if and only if there is a matching in $\cS(I)$ where each pair of siblings is matched to same-activity classes. Hence, (MSSP) can be solved in polynomial time, concluding the proof of Theorem~\ref{thm:after-school-siblings}.

\begin{figure}[t!]
\begin{center}
\begin{tabular}{c  l c | c c l }
$a_1$: & $b^{1}_{1} \AgentGreater{a_1}{}  b^{1}_{2}   \AgentGreater{a_1}{} \emptyset$ & \hspace{.15cm} &  & $b^{1}_{1}$: & $a_3 \AgentGreater{b^{1}_{1}}{} a_1 \AgentGreater{b^{1}_{1}}{} \emptyset$ \\ 
$\bar{a}_1$: &   $b^{2}_{1}  \AgentGreater{\bar{a}_1}{} b^{2}_{2}  \AgentGreater{\bar{a}_1}{} \emptyset$ & \hspace{.15cm} & & $b^{2}_{1}$: & $a_4 \AgentGreater{b^{2}_{1}}{} \bar{a}_1\AgentGreater{b^{2}_{1}}{} \emptyset$ \\
$a_2$: &  $ b^{3}_{1} \AgentGreater{a_2}{} \emptyset$  & \hspace{.15cm} & & $b^{3}_{1}$: & $a_2 \AgentGreater{b^{3}_{1}}{} \emptyset$ \\
$a_3$: &  $b^{1}_{2} \AgentGreater{a_3}{} b^{1}_{1}  \AgentGreater{a_3}{} \emptyset$  & \hspace{.15cm} & & $b^{1}_{2}$: & $a_1 \AgentGreater{b^{1}_{2}}{} a_3 \AgentGreater{b^{1}_{2}}{} \emptyset$ \\
$a_4$: & $b^{2}_{2} \AgentGreater{a_4}{} b^{2}_{1}  \AgentGreater{a_4}{} \emptyset$  & \hspace{.15cm} & & $b^{2}_{2}$: & $\bar{a}_1 \AgentGreater{b^{2}_{2}}{} a_4 \AgentGreater{b^{2}_{2}}{} \emptyset$
\vspace{3mm}
\end{tabular}
\end{center}
\caption{Example of an instance of (MSSP), with $\cC = \{(a_1,\bar{a}_1)\}$ and two activities: $b^{1}_{1}, b^{2}_{1},b^{3}_{1}$ are the classes of activity $c_1$; $b^{1}_{2}, b^{2}_{2}$ are the two classes of activity $c_2$. $a_1$ is eligible for $b^{1}_{1},b^{1}_{2}$; $\bar{a}_1$ is eligible for $b^{2}_{1},b^{2}_{2}$. All classes have a quota of $1$. Missing agents in the preference lists above can be added in any order after the outside option. The activity-stable matchings are the one-sided optimal stable matchings $M_0=\{a_1b^{1}_{1}, \bar{a}_1b^{2}_{1}, a_2b^{3}_{1}, a_3b^{1}_{2},a_4b^{2}_{2}\}$ (both $a_1$ and $\bar{a}_1$ attend a class of activity $c_1$) and $M_z=\{a_1b^{1}_{2}, \bar{a}_1b^{2}_{2}, a_2b^{3}_{1}, a_3b^{1}_{1},a_4b^{2}_{1}\}$ (both $a_1$ and $\bar{a}_1$ attend a class of activity $c_2$), while the other two stable matchings $M_1=\{a_1b^{1}_{2}, \bar{a}_1b^{2}_{1}, a_2b^{3}_{1}, a_3b^{1}_{1},a_4b^{2}_{2}\}$ and $M_2=\{a_1b^{1}_{1}, \bar{a}_1b^{2}_{2}, a_2b^{3}_{1}, a_3b^{1}_{2},a_4b^{2}_{1}\}$ are not activity-stable.}\label{fig:MSSP}
\end{figure}

Interestingly, unlike problem (MSSS) from Section~\ref{sec:matching-siblings-1}, we cannot linearize the function $f$ defined above. More generally, using Theorem~\ref{thm:main1} and Theorem~\ref{thm:linear-to-second-order}, we can show the following. %This follows from Theorem~\ref{thm:linear-to-second-order} and from the lemma below. 

\begin{lemma}\label{lem:counter-ex-not-MWSS}
There exists an instance of (MSSP) for which any minimum cut representable function whose minima coincide with the activity-stable matchings of the instance is not linearizable.
\end{lemma}

One might wonder whether the more general problem of maximizing the number of sibling pairs matched to the same-activity classes can be solved via minimum cut representability (recall that the analogous question for (MSSS) has a positive answer). The following lemma shows that this problem cannot be immediately addressed within our framework.

\begin{lemma}\label{lem:generalization-mssp-not-minc-representable}
There exists an instance of (MSSP) for which no minimum cut representable function $f$ over the rotation digraph satisfies that for every pair of stable matchings $M, M'$, $f(M)<f(M')$ if and only if $M$ matches strictly more pairs of sibling to the same activity than $M'$. 
\end{lemma}

\subsubsection{Matching Siblings with Different Preferences to Same Activities (MSDP).}\label{sec:matching-siblings-3}

The problem of Matching Siblings with Different Preferences to Same Activities (MSDP) is a modification of the (MSSP) problem where we drop the assumption that siblings rank activities in the same order. As before, the goal is to decide whether there exists a stable matching of students and classes so that each pair of sibling are enrolled in the same activity. We show the following. 

\begin{theorem}\label{thm:MSDP-Np-Complete}
(MSDP) is NP-Complete.
\end{theorem}

The proof of Theorem~\ref{thm:MSDP-Np-Complete} is based on a reduction from the problem of finding a matching that is stable in two different school matching instances~\citep{miyazaki2019jointly}. The next lemma employs Lemma~\ref{lem:if-mincut-then-lattice} to show that (MSDP) cannot be directly formulated as minimum cut representable optimization problem\footnote{Note that this fact is not implied by Theorem~\ref{thm:MSDP-Np-Complete} for two reasons: first, NP-Hardness implies that a problem cannot be solved in polynomial time only assuming P$\neq$NP; second, a problem may be minimum cut representable, but no algorithm may be known to find the witness graph in polynomial time.}. 

\begin{lemma}\label{lem:counter-ex-not-MSDP}
There exists an instance of (MSDP) for which any function whose minima coincide with the activity-stable matchings of the instance is not minimum cut representable.
\end{lemma}

\subsection{Two-stage Stochastic Stable Matching}\label{sec:2-stage}

We introduce and study a two-stage stochastic stable matching problem, where agents enter and leave the market between the two stages. The goal of the decision maker is, roughly speaking, to maximize the expected quality of the matchings across the two stages and minimize the expected students' discontent from being downgraded to a less preferred school when going from the first to the second-stage. Missing proofs can be found in the Appendix~\ref{ap:2-stage}. 

We formally define an instance of the two-stage stochastic stable matching problem as follows:

\noindent {\bf Aggregate market:} a set of students $A$, a set of schools $B$, and for every $a\in A$ (resp., $b \in B$), a strict order $\AgentGreater{a}{}$ (resp., $\AgentGreater{b}{}$) over $B^+ = B\cup\{\emptyset\}$ (resp., $A^+ = A\cup \{\emptyset\}$). $A$ and $B$ represent the set of all students and schools, those present in the first-stage and those that could potentially join in the second-stage.

\noindent {\bf First-stage instance:} a school matching instance $I=(A^I,B^I,\succ^I, q^I)$ of agents present in the first-stage where $A^I \subset A$, $B^I \subset B$ and $\succ^I$ is the restriction of the order $\succ$ to the agents $A^I \cup \{\emptyset\}$ and $B^I \cup \{\emptyset\}$\footnote{For every $a \in A^I$, $\succ^I_a$ (resp. $b \in B^I$, $\succ^I_b$) is a strict order over $B^I \cup \{\emptyset\}$ (resp. $A^I \cup \{\emptyset\}$) such that $b \succ^I_a b'$ if and only if $b \succ_a b'$ (resp. $a \succ^I_b a'$if and only if$a \succ_b a'$) for every $b, b' \in B^I \cup \{\emptyset\}$ (resp. $a, a' \in A^I \cup \{\emptyset\}$).}.

\noindent {\bf Second-stage instance:} a distribution $\cD$ over school matching instances. For every $J\sim\cD$, $J=(A^J, B^J, \succ^J, q^J)$ is a school matching instance where $A^J \subset A$, $B^J \subset B$ and $\succ^J$ is the restriction of $\succ$ to the agents $A^J \cup \{\emptyset\}$ and $B^J \cup \{\emptyset\}$. 

\noindent{\bf Cost:} cost functions $c_1, c_2 : A^+\times B^+ \rightarrow \mathbb{Q}$ and a penalty coefficient $\lambda \in \mathbb{Q}_+$. Let $J = (A^J, B^J, \succ^J, q^J)$ denote the sampled second-stage instance. Given a first-(resp., second-) stage matching $M^I$ (resp., $M^J$), we measure the quality of $M^I$ (resp., $M^J$) by $c_1(M^I) = \sum_{ab \in M^I} c_1(ab)$ (resp. $c_2(M^J) = \sum_{ab \in M^J} c_2(ab)$), 
and we measure the dissatisfaction of students for being moved to a less preferred school between $M^I$ and $M^J$ by 

\begin{align}
\label{eq:dissatisfation}
    d(M^I,M^J)=\lambda \sum_{a \in A^I \cap A^J} [R_a(M^J(a))-R_a(M^I(a))]^+.
\end{align}

Here, $R_a: B^+ \rightarrow \mathbb{N}$ is the \emph{rank function} of $a$, that is, $R_a(b)=i$ if and only if $b$ is the $i$-th most preferred choice (among the schools $B$ and the outside position) of $a$ and $[x]^+ = \max \{0, x\}$.

The coefficient $\lambda$ gives therefore the per unit of rank increase dissatisfaction of a student from switching to a school of higher rank between $M^I$ and $M^J$. This is a natural measure of dissatisfaction where students are unhappy to be downgraded to a less preferred school in the second-stage: the worse the school, the more unsatisfied they are. 

\noindent {\bf Objective function:} 
The goal is to solve the two-stage stochastic problem
{\small
\begin{align}
    \tag{2STO}
    \label{eq:2stagesto}
    \min_{M^I \in \cS(I)} c_1(M^I) + \mathbb{E}_{J \sim {\cal D}}\left[ \min_{M^J \in {\cS(J)}} c_2(M^J) + d(M^I, M^J)
    \right]
\end{align}
}

In~\eqref{eq:2stagesto}, we therefore wish to select a first-stage matching $M^I$ such that the cost of $M^I$ plus the expected cost that we have to pay in the second-stage is minimized.  
The second-stage cost is given by the cost of the second-stage matching plus the total dissatisfaction of the students for being downgraded to a less preferred school between the first- and second-stage. We note that once a first-stage matching $M^I$ is fixed and for every fixed second-stage instance $J$, the second-stage problem
$\min_{M^J \in {\cS(J)}} c_2(M^J) + d(M^I, M^J)
$
is an instance of (MWSM) and hence can be solved efficiently.

\begin{example}\label{ex:2-stage}
Consider the instance of 2-stage stochastic stable matching problem with all schools having a quota of $1$ and aggregate market coinciding with the first-stage instance $I$, with preferences given below:

\medskip 

\begin{center} 
\begin{tabular}{c  l c | c c l }
$a^1$: & $b^1 \AgentGreater{a^1}{}  b^2   \AgentGreater{a^1}{} b^4 \AgentGreater{a^1}{} b^3 \AgentGreater{a^1}{} \emptyset$ & \hspace{.15cm} &  & $b^1$: & $a^2 \AgentGreater{b^1}{} a^3 \AgentGreater{b^1}{} a^1 \AgentGreater{b^1}{} \emptyset$ \\ 
$a^2$: &   $b^2  \AgentGreater{a^2}{} b^3 \AgentGreater{a^2}{} b^5 \AgentGreater{a^2}{} b^1  \AgentGreater{a^2}{} \emptyset$ & \hspace{.15cm} & & $b^2$: & $a^3 \AgentGreater{b^2}{} a^1 \AgentGreater{b^2}{} a^4 \AgentGreater{b^2}{} a^2 \AgentGreater{b^2}{} \emptyset$ \\
$a^3$: &  $ b^3 \AgentGreater{a^3}{} b^5 \AgentGreater{a^3}{} b^1 \AgentGreater{a^3}{} b^2 \AgentGreater{a^3}{} \emptyset$  & \hspace{.15cm} & & $b^3$: & $a^1 \AgentGreater{b^3}{} a^5 \AgentGreater{b^3}{} a^4 \AgentGreater{b^3}{} a^2 \AgentGreater{b^3}{} a^3 \AgentGreater{b^3}{} \emptyset$. \\
$a^4$: &  $b^2\AgentGreater{a^4}{} b^3 \AgentGreater{a^4}{} b^4 \AgentGreater{a^4}{} b^5 \AgentGreater{a^4}{} \emptyset$  & \hspace{.15cm} & & $b^4$: & $a^1 \AgentGreater{b^4}{} a^5 \AgentGreater{b^4}{} a^4 \AgentGreater{b^4}{} \emptyset$ \\
$a^5$: & $b^5\AgentGreater{a^5}{} b^4 \AgentGreater{a^5}{} b^3  \AgentGreater{a^5}{} \emptyset$  & \hspace{.15cm} & & $b^5$: & $a^4 \AgentGreater{b^5}{} a^2 \AgentGreater{b^5}{} a^3 \AgentGreater{b^5}{} a^5 \AgentGreater{b^5}{} \emptyset$
\vspace{3mm}
\end{tabular}
\end{center} 

\smallskip 

${\cal D}$ assigns $1/2$ probability to each of the two instances $J_1$, given by 

\medskip 

\begin{center} 
\renewcommand{\arraystretch}{1.3}
\begin{tabular}{c  l c | c c l }
$a^1$: & $b^1 \AgentGreater{a^1}{J_1}  b^2   \AgentGreater{a^1}{J_1} b^3 \AgentGreater{a^1}{J_1} \emptyset$ & \hspace{.15cm} &  & $b^1$: & $a^2 \AgentGreater{b^1}{J_1} a^3 \AgentGreater{b^1}{J_1} a^1 \AgentGreater{b^1}{J_1} \emptyset$  \\
$a^2$: &   $b^2  \AgentGreater{a^2}{J_1} b^3 \AgentGreater{a^2}{J_1} b^1  \AgentGreater{a^2}{J_1} \emptyset$ & \hspace{.15cm} & & $b^2$: & $a^3 \AgentGreater{b^2}{J_1} a^1 \AgentGreater{b^2}{J_1} a^2 \AgentGreater{b^2}{J_1} \emptyset$ \\ 
$a^3$: &  $ b^3 \AgentGreater{a^3}{J_1} b^1 \AgentGreater{a^3}{J_1} b^2 \AgentGreater{a^3}{J_1} \emptyset$  & \hspace{.15cm} & & $b^3$: & $a^1  \AgentGreater{b^3}{J_1} a^2 \AgentGreater{b^3}{J_1} a^3 \AgentGreater{b^3}{J_1} \emptyset$
\vspace{3mm}
\end{tabular},
\end{center} 

\smallskip 

and $J_2$, given by

\medskip 

\begin{center} 
\renewcommand{\arraystretch}{1.3}
\begin{tabular}{c  l c | c c l }
$a^1$: & $b^1 \AgentGreater{a^1}{J_2}  b^4 \AgentGreater{a^1}{J_2} \emptyset$ & \hspace{.15cm} &  & $b^1$: & $ a^1 \AgentGreater{b^1}{J_2} \emptyset$ \\
$a^4$: &  $ b^4 \AgentGreater{a^4}{J_2} b^5 \AgentGreater{a^4}{J_2} \emptyset$  & \hspace{.15cm} & & $b^4$: & $a^1 \AgentGreater{b^4}{J_2} a^5 \AgentGreater{b^4}{J_2} a^4 \AgentGreater{b^4}{J_2} \emptyset$ \\
$a^5$: & $b^5\AgentGreater{a^5}{J_2} b^4  \AgentGreater{a^5}{J_2} \emptyset$  & \hspace{.15cm} & & $b^5$: & $a^4 \AgentGreater{b^5}{J_2} a^5 \AgentGreater{b^5}{J_2} \emptyset$
\vspace{3mm}
\end{tabular}.
\end{center} 

\smallskip 

Assume $\lambda=1$ and $c_2$ to be the all-zero vector. Consider the first-stage stable matching $M^I=\{a^1b^1, a^2b^3, a^3b^5, a^4b^2, a^5b^4\}$. With this special objective function, it is optimal to choose in the second-stage the student-optimal stable matching. Thus, if $J_1$ is realized, then we select the second-stage stable matching $M^{J_1}=\{a^1b^1, a^2b^2, a^3b^3\}$ with cost $0$ (none of the students get downgraded), while if $J_2$ is realized, we select $M^{J_2}=\{a^1b^1, a^4b^4, a^5b^5\}$ with cost $2$ (as $a^4$ is downgraded by $2$ positions). Thus, the objective function~\eqref{eq:2stagesto} evaluated in $M^I$ is equal to $c_1(M^I)+1$.

\end{example}

\begin{remark}\label{rem:third-cost}
We note that our algorithm and analysis extend seamlessly to capture the following generalization of the dissatisfaction function~\eqref{eq:dissatisfation}: given $\lambda\in \mathbb{Q}_+$, and a set of non-negative scores for the pairs $w:A^+\times B^+ \rightarrow \mathbb{Q}_+$, with the property that $w(ab)\geq w(ab')$ if $b \AgentGreater{a}{} b'$, define
$$d(M^I,M^J)=\lambda\sum_{a \in A^I \cap A^J} [w(aM^J(a)) - w(aM^I(a))]^+.$$ 
This more general version of the dissatisfaction function captures, for example, the setting when the per unit of rank change dissatisfaction from changing school is student-dependent (in which case $w(ab) = \alpha_a \cdot R_a(b)$ for some set of non-negative weights $\{\alpha_a\}_{a \in A}$) or where a student does not care if the change happens between, say, their top $5$ schools (in which case $w(ab_1)=w(ab_2)=\dots=w(ab_5)$, where $b_1,\dots,b_5$ are the top 5 schools in $a$'s list) but cares if it happens between one of their top $5$ schools and the rest. For simplicity of exposition, we restrict the presentation to the rank change dissatisfaction function~\eqref{eq:dissatisfation}.
\end{remark}

\subsubsection{Explicit Second-Stage Distribution.}

We first consider the case when the second-stage distribution $\cD$ is given explicitly by the list of all possible second-stage scenarios $\cal{J}=\{J_1,\dots,J_K\}$ and their probabilities of occurrence $\{p^{J_k}\}_{k \in [K]}$. In this setting, \eqref{eq:2stagesto} can be written as follows,

\begin{align*}
\tag{EXP-2STO}
\label{eq:exp2sto}
    \min_{\substack{
    M^I \in \cS(I)\\
    \{M^{J_k} \in \cS(J_k)\}_{k \in [K]}}} c_1(M^I) + \sum_{k\in [K]} p^{J_k}\left( c_2(M^{J_k}) + d(M^I, M^{J_k})\right),
\end{align*}

We now use our framework to show the following, \begin{theorem}\label{thm:two-stage-main1} 
Problem \eqref{eq:exp2sto} can be solved in time polynomial in the input size\footnote{For $x \in \mathbb{Q}$, let $\|x\|$ denote the encoding length of parameter $x$. Then ${\sf N}^{\sf exp} = O(\sum_{k \in K} \|p^{J_k}\| + \sum_{ab \in A \times B} (\|c_1(ab)\| + \|c_2(ab)\|) + \|\lambda\| + |A \cup B|^2)$.} ${\sf N}^{\sf exp}$.
\end{theorem}

To prove Theorem~\ref{thm:two-stage-main1}, we formulate~\eqref{eq:exp2sto} as an optimization problem of the form~\eqref{problem}, then show that the resulting problem is minimum cut representable. We start with discussing a special case. %\yfcomment{I'd omit the following. Everything we say below is correct and formal, so it's not just a sketch. We have already said where missing proofs are.} The full details are given in the Electronic Companion~\ref{ap:2-stage}. We give below a sketch of our approach.

\begin{figure}[ht]
    \centering
    \begin{tikzpicture}[node distance=2cm, >=stealth, thick]

    % Nodes positioned with more spacing between chains
    \node (s) at (0, -1.5) [circle, draw] {$s$};
    \node (rho1) at (-2, 0) [circle, draw] {$\rho_1$};
    \node (rho2) at (-2, 1.5) [circle, draw] {$\rho_2$};
    \node (rho3) at (-2, 3) [circle, draw] {$\rho_3$};
    \node (rho1p) at (2, 1.5) [circle, draw] {$\rho_{1}'$};
    \node (t) at (0, 4.5) [circle, draw] {$t$};

    % Directed edges with capacities
    \draw[->] (rho3) -- (t) node[midway, above] {$\frac{1}{2}$};
    \draw[->] (rho3) -- (rho2) node[midway, left] {$\infty$};
    \draw[->] (rho2) -- (rho1) node[midway, left] {$\infty$};
    
    \draw[->] (s) -- (rho1) node[midway, left] {$1$};
    \draw[->] (s) -- (rho2) node[midway, left] {$1$};
    \draw[->] (s) -- (rho3) node[midway, above] {$1$};
    \draw[->] (s) -- (rho1p) node[midway, below] {$1$};
    
    % Curved bidirectional arcs between rho3 and rho1'
    \draw[->] (rho3) to [bend left=20] node[midway, above] {$\frac{1}{2}$} (rho1p);
    \draw[->] (rho1p) to [bend left=20] node[midway, below] {$\frac{1}{2}$} (rho3);

    \draw[->] (rho1p) -- (t) node[midway, right] {$\frac{3}{2}$};

    \end{tikzpicture}
    \caption{The digraph $D(I,J_2)$ that our algorithm computes for the deterministic $2$-stage instance with first-stage instance  $I$ and second-stage instance $J_2$, defined as in Example~\ref{ex:2-stage}. We set $\lambda=1$ and $c_1,c_2$ be the $0$ vectors. Nodes on the left correspond to rotations $\rho_1=((a^1b^1,a^2b^3,a^3b^5,a^4b^2),(a^1b^2,a^2b^5,a^3b^1,a^4b^3))$, $\rho_2=((a^4b^3,a^5b^4),(a^4b^4,a^5b^3))$, and $\rho_3=((a^1b^2,a^2b^5,a^3b^1,a^4b^4),(a^1b^4,a^2b^1,a^3b^2,a^4b^5))$ from instance $I$. The node on the right corresponds to the unique rotation $\rho_1'=((a^4b^4,a^5b^5),(a^4b^5,a^5b^4))$ from $J_2$. Arcs of capacity $0$ are not represented.\\
    With this special objective function, the objective value of~\eqref{eq:exp2sto} computed in a pair of first- and second-stage stable matchings $M_1$ and $M_2$ is equal to the total rank drop of students between the first- and second-stage. The above digraph is such that the $s-t$ cut corresponding to any pair of first- and second-stage stable matchings $M_1,M_2$ has capacity equal to the objective value of $M_1, M_2$ in~\eqref{eq:exp2sto} plus $2$. For example, the $s$-$t$ cut $\{s,\rho_1,\rho_1'\}$ corresponds to the pair of stable matchings $M_0^{I}/\rho_1=M^{I}_1=\{a^1b^2,a^2b^5,a^3b^1,a^4b^3,a^5b^4\}$, $M^{J_2}_0/\rho_1'=M^{J_2}_1=\{a^1b^1,a^4b^5,a^5b^4\}$ whose objective value in~\eqref{eq:exp2sto} is $2$, and the cut has a value $2+2 = 4$.}
    \label{fig:digraph-2ss-example}
\end{figure}
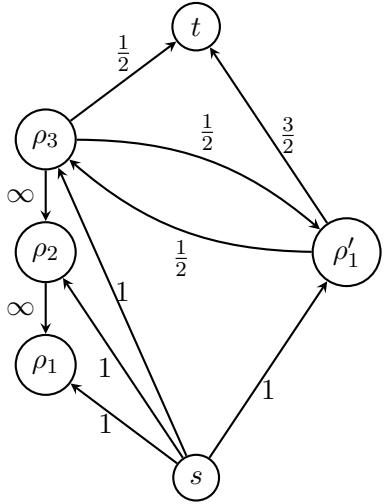

 \noindent{\bf The deterministic case.} When $K=1$, the second-stage instance consists of a single school matching instance $J$ that is known in hindsight. The high-level idea for solving this case is to reduce~\eqref{eq:exp2sto} to a minimum-cut representable problem over a school matching instance $I^\sqcup$ with witness digraph $\cD(I,J)$ such that (i) the vertices $\cD(I,J)$ are $s,t$, and a copy of each rotation from $\cR(I)$ and $\cR(J)$ (so a rotation present in both $\cR(I)$ and $\cR(J)$ appears twice in $\cD(I,J)$); moreover, (ii) every finite $s$-$t$ cut $S$ in $\cD(I,J)$ corresponds to a pair of upper-closed subsets $R^I,R^J$ of $(\cR(I),\trianglerighteq^I)$, $(\cR(J),\trianglerighteq^J)$, and thus to a pair of stable matchings $M^I,M^J$ of $I$ and $J$, respectively, and (iii) the value of $S$ equals (up to an additive constant) to the objective function \eqref{eq:exp2sto} computed in $M^I,M^J$. Therefore, from a minimum cut of $\cD(I,J)$ we can derive an optimal pair of stable matchings $M^I, M^J$. See Figure~\ref{fig:digraph-2ss-example} for an example. Thus,  $I^\sqcup$ can be thought of as the ``disjoint union'' of $I,J$, and the stable matchings of $I,J$ can be thought of as ``projections'' of the stable matchings in $I^\sqcup$. We next make these ideas formal and extend them to the general case.

{\noindent \bf Formulation of~\eqref{eq:exp2sto} as an optimization problem of the form~\eqref{problem}.} We begin with some introductory concepts.

\begin{definition}\label{def:disjoint-union}
We define the following elements.
    \begin{itemize}
        \item (Disjoint Union of Sets) Given a family $X_1, \dots, X_K$ of finite sets, we define the \emph{disjoint union} of $X_1, \dots, X_K$ as $\bigsqcup_{k=1}^K X_k = \bigcup_{k=1}^K \{(x_k, X_k) \;|\; x_k \in X_k\}.$
        \item (Disjoint Union of School Matching Instances) Given a family of school matching instances $I_1, \dots, I_K$, where $I_k = (A^{I_k}, B^{I_k}, \succ^{I_k}, q^{I_k})$ for every $k \in [K]$, denote $A^{\sqcup} = \bigsqcup_{k=1}^K A^{I_k}$ and $B^{\sqcup} = \bigsqcup_{k=1}^K B^{I_k}$. Define $q^{\sqcup}: B^{\sqcup} \rightarrow \mathbb{N}$ such that for every $(b, B^{I_k}) \in B^{\sqcup}$ where $k \in [k]$, we have $q^{\sqcup}_{(b, B^{I_k})} = q^{I_k}_b$. Consider a school $(b, B^{I_k}) \in B^{\sqcup}$ for some $k \in [K]$, we define the total order $\succ^{\sqcup}_{(b, B^{I_k})}$ over the students $A^{\sqcup}$ and outside option $\emptyset$ such that if 
        $a_1 \succ^{I_k}_b \dots \succ^{I_k}_b a_r \succ^{I_k}_b \emptyset \succ^{I_k}_b a_{r+1} \succ^{I_k}_b \dots \succ^{I_k}_b a_n$ denotes the preference list of $b$ in $I^{k}$ then $(b, B^{I_k})$ ranks the students from instance $I^{k}$ first in the same order given by $\succ^{I_k}_b$, i.e.,
        $$
        (a_1, A^{I_k}) \succ^{\sqcup}_{(b, B^{I_k})} \dots\succ^{\sqcup}_{(b, B^{I_k})} (a_r, A^{I_k})\succ^{\sqcup}_{(b, B^{I_k})} \emptyset\succ^{\sqcup}_{(b, B^{I_k})} (a_{r+1}, A^{I_k})\succ^{\sqcup}_{(b, B^{I_k})} \dots \succ^{\sqcup}_{(b, B^{I_k})} (a_n, A^{I_k}),
        $$ and then ranks the remaining students in any arbitrary order at the end of its list. A total order of students over the schools and the outside option can be defined analogously.
        We define the \emph{disjoint union of $I_1, \dots, I_K$} as the school matching instance
        $
        \bigsqcup_{k=1}^K I_k = (A^{\sqcup}, B^{\sqcup}, \succ^{\sqcup}, q^{\sqcup}).
        $
        \item (Restrictions of  Stable Matchings) Given a family of school matching instances $I_1, \dots, I_K$, where $I_k = (A^{I_k}, B^{I_k}, \succ^{I_k}, q^{I_k})$ for every $k \in [K]$, let $I^{\sqcup} = \bigsqcup_{k=1}^K I_k = (A^{\sqcup}, B^{\sqcup}, \succ^{\sqcup}, q^{\sqcup})$ denote their disjoint union, and let $M \in \cS(I^{\sqcup})$. Note that for every $(a, A^{I_k}) \in A^{\sqcup}$ and $(b, B^{I_{k'}}) \in B^{\sqcup}$, if $M((a, A^{I_k})) = (b, B^{I_{k'}})$ then $k=k'$ (as otherwise $(a, A^{I_k})$ and $(b, B^{I_{k'}})$ rank each other after the outside option). Hence, $M$ can be partitioned into $K$ subsets $M^{(I_1)}, \dots, M^{(I_K)}$ where for every $k \in [K]$, 
        $$
        M^{(I_k)} = \{(a, A^{I_k})(b, B^{I_k})\}_{\substack{ab \in A^{I_k} \times B^{I_k} \\ M((a, A^{I_k})) = (b, B^{I_k})}} \cup \{(a, A^{I_k})\emptyset\}_{\substack{ a \in A^{I_k} \\ M((a, A^{I_k})) = \emptyset}} \cup \{\emptyset (b, B^{I_k})\}_{\substack{b \in B^{I_k}\\ M((b, B^{I_k})) = \emptyset}}.
        $$
        We refer to $M^{(I_k)}$ as the \emph{restriction of $M$ to instance $I_k$}. By construction, every $a \in A^{I_k}$ appears in a unique pair $(a, A^{I_k})(b, B^{I_k})$ (or $(a, A^{I_k})\emptyset$) in $M^{(I_k)}$. We denote $M^{(I_k)}(a)=b$ (or $M^{(I_k)}(a)=\emptyset$).
        % \item ()
        \item (Extensions of Linear Functions) Given a family of school matching instances $I_1, \dots, I_K$, where $I_k = (A^{I_k}, B^{I_k}, \succ^{I_k}, q^{I_k})$ for every $k \in [K]$, let $I^{\sqcup} = \bigsqcup_{k=1}^K I_k = (A^{\sqcup}, B^{\sqcup}, \succ^{\sqcup}, q^{\sqcup})$ denote their disjoint union, let $A = \bigcup_{k=1}^K A^{I_k}$, $B = \bigcup_{k=1}^K B^{I_k}$, and let $M \in \cS(I^{\sqcup})$. Given a linear function $c: A^+ \times B^+ \rightarrow \mathbb{Q}$ and $k \in [K]$, we define $c^{(I_k)}: (A^{\sqcup})^+\times (B^{\sqcup})^+ \rightarrow \mathbb{Q}$ such that for every $k', k'' \in [K]$, $a \in A^{I_{k'}}$ and $b \in B^{I_{k''}}$,
        $$
        c^{(I_k)}((a, A^{I_{k'}}) (b, B^{I_{k''}})) 
            =
            \left\{
            \begin{matrix}
                c(ab) & \text{ if  $k'=k''=k$}\\
                0 & \text{ otherwise }
            \end{matrix},
            \right.
        $$
        and $c^{(I_k)}((a, A^{I_{k'}}) \emptyset)=c^{(I_k)}( \emptyset (b, B^{I_{k''}})) = c^{(I_k)}(\emptyset \emptyset) = 0$.
    \end{itemize}
\end{definition}

We are now ready to give our reformulation of~\eqref{eq:exp2sto}. Let $I^{\sqcup} = I \sqcup \bigsqcup_{k=1}^K J_k = (A^{\sqcup}, B^{\sqcup}, \succ^{\sqcup}, q^{\sqcup})$ denote the disjoint union of the first-stage instance $I$ and the second-stage instances $\{J_k\}_{k\in [K]}$. Define $\cF = \cS(I^{\sqcup})$, and define $f :{\cal F} \rightarrow \mathbb{R}$ such that,
$$f(M)=c_1^{(I)}(M) + \sum_{k \in [K]} p^{J_k}c^{(J_k)}_2(M) + \lambda \sum_{k \in [K]} p^{J_k}\sum_{a \in A^I \cap A^{J_k}} [R_a(M^{(J_k)}(a))-R_a(M^{(I)}(a))]^+.$$ 
The following lemma gives a reformulation of~\eqref{eq:exp2sto} as an optimization problem of the form~\eqref{problem}.
\begin{lemma}
\label{lem:reformulation}
    Let $f$ and $\cF$ be defined as above. Then 
    \eqref{eq:exp2sto} and~\eqref{problem} are equivalent in the sense that given a feasible solution of one, a feasible solution of the other with same objective value can be efficiently constructed.
\end{lemma}

\vspace{3mm}

\noindent{\bf Minimum cut representability of~\eqref{problem}.} Note that $$f = \;\; f^1 \;\; + \;\; \sum_{k \in [K]} p^{J_k} f_k^2 \;\; + \;
\;\lambda \sum_{k \in [K]} p^{J_k}  \sum_{a \in A^I \cap
A^{J_k}}f_{a, k}^3,$$ where
$f^1(M)=c^{(I)}_1(M)$, $f_k^2(M)= c^{(J_k)}_2(M) $, and $f_{a,k}^3(M)= [R_a(M^{(J_k)}(a))-R_a(M^{(I)}(a))]^+$ for every $M \in \cS(I^{\sqcup})$. Hence by Lemma~\ref{lem:additivity}, in order to show that~\eqref{problem} defined above is minimum cut representable, it is sufficient to show that $(\Pi(f^1, \cF))$, $(\Pi(f^2_k, \cF))$ for every $k \in [K]$ and $(\Pi(f^3_{a,k}, \cF))$ for every $a \in A^I \cap A^{J_k}$ and $k \in [K]$ are minimum cut representable. Note that $f^1$ and $f^2_k$ for every $k \in [K]$ are linear functions by definition and that $\cF$ is the lattice of all stable matchings of $I^{\sqcup}$ hence the corresponding problems are minimum cut representable. It remains to show that $(\Pi(f^3_{a,k}, \cF))$ for every $a \in A^I \cap A^{J_k}$ and $k \in [K]$ is minimum cut representable. For this purpose, we compute the first and second order differentials of $f^3$, %in Lemma~\ref{lem:fod-two-stage} and Lemma~\ref{lem:sod-two-stage}, respectively
and use these to verify the properties in Theorem~\ref{thm:main1}. Details are in the Appendix~\ref{ap:2-stage-exp}.

\subsubsection{Implicit Second-Stage Distribution.}

We consider in this section the more general model where the second-stage distribution $\cD$ is given by a sampling oracle. Let ${\sf N}^{\sf imp}$ denote the input size\footnote{${\sf N}^{\sf imp} = O(\sum_{ab \in A \times B} (\|c_1(ab)\|+\|c_2(ab)\|) + \|\lambda\| +|A\cup B|^2)$.} in this setting. 
We show the following hardness result for the two-stage problem \eqref{eq:2stagesto} under this more general model.

\begin{theorem}
\label{thm:imp-hard}
    Unless P$=$NP, there exists no algorithm that given any instance $\cI$ of \eqref{eq:2stagesto} where the second-stage distribution is specified implicitly by a sampling oracle, solves\footnote{In the strong sense that it gives at least one of the two, the optimal value or the optimal solution.} the problem in time and number of calls to the sampling oracle that is polynomial in the input size ${\sf N}^{\sf imp}$. This hardness result holds even if the cost parameters $\lambda$ and $\{c_1(ab)\}_{ab}, \{c_2(ab)\}_{ab}$ are in $\{0,1\}$.%\todo{Y: notation $\{c_1(ab)\}_{ab}$ is a little strange (where does $ab$ belong to?) just use $c_1,c_2$?}
\end{theorem}
Our proof of Theorem~\ref{thm:imp-hard} relies on a reduction from the problem of counting the number of vertex covers of an undirected graph, which is \#P-Hard \citep{provan1983complexity}.

On the positive side, we give an arbitrary good additive approximation to \eqref{eq:2stagesto} when the second-stage distribution is specified implicitly. Our algorithm runs in time pseudopolynomial in the input size. In particular, let $M^{I, *}$ be the optimal solution of the two stage problem \eqref{eq:2stagesto} and let $\val{}{M^I}$ denote the objective value of a first-stage stable matching $M^I$. We show the following result.

\begin{theorem}\label{thm:imp-fpras}
There exists an algorithm that, given a second-stage distribution that is specified implicitly by a sampling oracle, and two parameters $\epsilon > 0$ and $\alpha \in (0,1)$, gives a first-stage stable matching solution $M^I$ such that,
$$
\mathbb{P}\left(\val{}{M^I} \leq \val{}{M^{I,*}} + \epsilon \right) \geq 1-\alpha.
$$
The algorithm runs in time polynomial in ${\sf N}^{\sf imp}$, $\max_{ab\in A\times B} |c_2(ab)|$, $\lambda$, $1/\epsilon$ and $\log(1/\alpha)$.
\end{theorem}

The approximation algorithm from Theorem~\ref{thm:imp-fpras} runs in polynomial time when $\lambda$ and the values of $c_2$ are polynomially bounded (recall that, already in this case, by Theorem~\ref{thm:imp-hard}, no exact polynomial-time algorithm exists unless P$=$NP). Our algorithm employs the widely used \emph{Sample Average Approximation} method (e.g. \cite{charikar2005sampling,ravi2006hedging,swamy2005sampling}) 
to approximate an instance of \eqref{eq:2stagesto} with an implicitly specified second-stage distribution by an instance of explicitly specified second-stage distribution. Then we show that a relatively small number of samples (calls to the sampling oracle) is enough to get a good additive approximation with high probability.
Whether there exists and FPRAS for \eqref{eq:2stagesto} in this model is an interesting open question.

\vspace{3mm}
\noindent{\bf Algorithm.} 
Our starting point is to approximate the expected value in the two-stage stochastic problem \eqref{eq:2stagesto} by a sample average over $K$ i.i.d.~samples. 
In particular, let $J_1, \dots, J_K$ be $K$ i.i.d.~samples drawn from the distribution $\cD$. We replace the expected value in the two-stage stochastic problem \eqref{eq:2stagesto} with the sample average taken over $J_1, \dots, J_K$. We get the following sample average minimization problem,

\begin{align}
\tag{SAA}
\label{eq:saa}
    \min_{M^I \in \cS(I)}
     c_1(M^I) + \frac{1}{K} \sum_{k=1}^K \min_{
    M^{J_{k}} \in \cS(J_k)
    } c_2(M^{J_k}) + d(M^I, M^{J_k}).
\end{align}

Note that \eqref{eq:saa} is an instance of \eqref{eq:exp2sto} where $p^{J_k} = \frac{1}{K}$ for every $k \in [K]$, and can therefore be solved exactly using our algorithm for explicitly specified distributions, in time polynomial in $K$ and the other input parameters. It remains to show that, for any $\epsilon > 0$, a small number of samples $K$ is enough to get an $\epsilon$ additive approximation of~\eqref{eq:2stagesto} with high probability. Let $\hat{M}^{I, *}$ be an optimal solution of the sample average problem \eqref{eq:saa}. The following lemma bounds the value of the first-stage stable matching $\hat{M}^{I, *}$ in the original problem \eqref{eq:2stagesto} as a function of the number of samples $K$. 
\begin{lemma}\label{lem:samplecomplexity}
    Let $\alpha \in (0,1)$. Then with probability at least $1-\alpha$ it holds that,
    $$
    \val{}{\hat{M}^{I, *}} \leq \val{}{M^{I, *}} + 4|A|(\max_{ab} |c_2(ab)| + \lambda |B|)\sqrt{\frac{\max\{|A|, |B|\}\log(3.88) - \log(\alpha)}{K}}.
    $$
\end{lemma}
Theorem~\ref{thm:imp-fpras} follows immediately from Lemma~\ref{lem:samplecomplexity}. In particular, for any $\epsilon > 0$, a number of samples $$K = \left(4|A|(\max_{ab} |c_2(ab)| + \lambda |B|)\right)^2\frac{\max\{|A|, |B|\}\log(3.88)-\log(\alpha)}{\epsilon^2}$$ is enough to get an $\epsilon$ additive approximation with probability at least $1-\alpha$.

\subsubsection{Numerical Experiments.}\label{sec:numerical}
\paragraph{Experimental Setup.} We consider an instance of \eqref{eq:2stagesto} with $n=50$ students and schools, all of which have a quota of $1$. The preference lists are  generated uniformly at random. We suppose that each agent leaves the market in the second-stage with probability $p=0.25$ independently of the other agents. In terms of costs, we set the cost of a pair $ab$, $c_1(ab)$ and $c_2(ab)$, to be the average rank between the two sides\footnote{This is known as an \emph{egalitarian} stable matching (see, e.g.,~\cite{gusfieldbook}), since it balances between the utilities of the two sides of the market.}, and solve the problem for different values of the penalty coefficient $\lambda\in (0,2]$. We compare the performance of four first-stage stable matchings: the optimal matching computed by our algorithm denoted by $M^*$, the one-side optimal matchings $M_0$ and $M_z$ that are commonly used in practice, and the optimal matching %that maximizes the optimal stable matching 
in hindsight $M^{\sf off}$, that is, the first-stage stable matching that would have been optimal to take if one knew the realization of the second-stage distribution exactly.

\begin{figure}[h]
    \centering
    \includegraphics[scale=0.6]{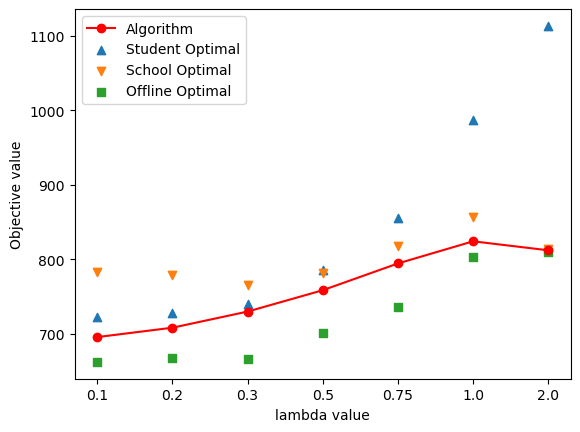}
    \caption{Comparison of the objective value of the optimal stable matching with other stable matchings commonly used in practice, on an instance with $50$ students and $50$ schools with randomly generated preferences.}
    \label{fig:Simulation}
\end{figure}

\paragraph{Results.} We observe that for a wide range of values of $\lambda$, the stable matching computed by our algorithm considerably outperforms the one-sided optimal matchings $M_0$ and $M_z$, which are commonly used in practice. Furthermore, we observe that the optimal value $\val{}{M^*}$ stays relatively close to the optimal value in hindsight, $\val{}{M^{\sf off}}$ for all values of $\lambda$. This suggests that if minimizing the rank change of students is a priority, our solution provides an effective approach to use in practice.

\paragraph{Choosing the coefficient $\lambda$.} The coefficient $\lambda$ balances the trade-off between minimizing changes across the two stages and maximizing the quality of the matchings. In practice, the appropriate value of $\lambda$ is application-dependent, and a good coefficient can be learned through qualitative feedback from participants (e.g., students and faculty).

\smallskip 

\paragraph{Acknowledgements.} Yuri Faenza acknowledges support from the NSF grant 2046146.

\printbibliography

\newpage
\appendix

{\noindent \Huge Appendix}

\section{Preliminaries: Omitted Proofs}

\subsection{Proof of Lemma~\ref{lem:meta_rotations_is_poset}}
%\afcomment{this is now checked}
We have,
\begin{itemize}
    \item Reflexivity: For every $\theta \in \cR_{\sf prop}^{\cF}(I)$,  $\theta \trianglerighteq^{\cF} \theta$ by definition.
    \item Antisymmetry: Suppose $\theta \trianglerighteq^{\cF} \theta'$ and $\theta' \trianglerighteq^{\cF} \theta$. Then, in any stable matching $M \in \cF$, if the rotations $\theta$ are eliminated, the rotations $\theta'$ are also eliminated, and vice versa. This implies that $\theta$ and $\theta'$ correspond to the same equivalence class of $\sim$ and hence $\theta=\theta'$.

    \item Transitivity: Suppose $\theta \trianglelefteq^{\cF} \theta'$ and $\theta' \trianglelefteq^{\cF} \theta''$. Then, for any stable matching $M \in \cF$, if the rotations $\theta$ are eliminated in $M$, the rotations $\theta'$ are also eliminated in $M$ by $\theta \trianglelefteq^{\cF} \theta'$, and then by $\theta' \trianglelefteq^{\cF} \theta''$ the rotations $\theta''$ are also eliminated in $M$. This implies that $\theta \trianglelefteq^{\cF} \theta''$.
\end{itemize}
\hfill $\square$
\subsection{Proof of Theorem~\ref{thm:representation_meta_rotations}}
Fix a stable matching $M \in \cF$, we show that $R_M$ can be partitioned into meta-rotations $\theta_0^{\cF}, \theta_1, \dots, \theta_r$ such that $\{\theta_1, \dots, \theta_r\}$ is an upper-closed subset of proper meta-rotations. Then given an upper-closed subset of proper meta-rotations $\Theta$, we show that the set of rotations $\theta_0^{\cF} \cup \bigcup_{\theta \in \Theta} \theta$ is upper-closed and the corresponding matching is in $\cF$. Then we show that the bijection $M \rightarrow \Theta_M$ (whose inverse is $\Theta \rightarrow M_\Theta$ by construction) is an isomorphism of lattices between $\cF$ and upper-closed subsets of proper meta-rotations.

For the first point, note that $R_M$ must contain the rotations $\theta_0^{\cF}$ that are eliminated in every stable matching of $\cF$ and not contain $\theta_{z}^{\cF}$ that are never eliminated in every stable matching of $\cF$, and for every other proper meta-rotation $\theta$, either all rotations in $\theta$ are in $R_M$ or outside $R_M$ (rotations of a meta-rotation are eliminated together). Hence, $R_M$ can be partitioned into a set of meta-rotations $\theta_0^{\cF}, \theta_1, \dots, \theta_r$. Now, fix $\theta \in \{\theta_1, \dots, \theta_r\}$ and let $\theta' \trianglerighteq^{\cF} \theta$. Since $\theta$ are eliminated in $M$, $\theta'$ must also be eliminated in $M$, hence $\theta' \in \{\theta_1, \dots, \theta_r\}$ and $\{\theta_1, \dots, \theta_r\}$ is indeed upper-closed.

For the second point, consider $\rho \in \theta_0^{\cF} \cup \bigcup_{\theta \in \Theta} \theta$, and let $\rho' \trianglerighteq \rho$. We show that $\rho' \in \theta_0^{\cF} \cup \bigcup_{\theta \in \Theta} \theta$. If $\rho' \in \theta_0^{\cF}$ we are done. Otherwise, if $\rho' \in \theta_{z}^{\cF}$, then because $\rho \notin \theta_{z}^{\cF}$, there exists $M \in \cF$ where $\rho$ is eliminated, then since $\rho' \trianglerighteq \rho$, $\rho'$ must also be eliminated in $M$, which is absurd since the rotations in $\theta_{z}^{\cF}$ are never eliminated in any matching of $\cF$. Finally, if $\rho'$ belongs to some proper meta-rotation (that we denote by $\overline{\rho'}$) then $\rho \notin \theta_0^{\cF}$ otherwise $\rho'$ must be also eliminated in every stable matching of $\cF$ as $\rho' \trianglerighteq \rho$. Let $\overline{\rho}$ denote the proper meta-rotation to which $\rho$ belongs. Then, consider a stable matching $M \in \cF$ where the rotations $\overline{\rho}$ are eliminated. Since $\rho' \trianglerighteq \rho$, $\rho'$ is also eliminated in $M$, implying that all the rotations in the proper meta-rotation $\overline{\rho'}$ are eliminated in $M$. This implies that $\overline{\rho'} \trianglerighteq^{\cF} \overline{\rho}$, and hence that $\overline{\rho'} \in \Theta$ (because $\Theta$ is upper-closed) and finally that $\rho' \in \theta_0^{\cF} \cup \bigcup_{\theta \in \Theta} \theta$. We now show that $M_{\theta_0^{\cF} \cup \bigcup_{\theta \in \Theta} \theta}$ is in $\cF$. First note that since for every $\theta \in \Theta$, there exists a stable matching $M(\theta) \in \cF$ in which $\theta$ are eliminated (which holds because $\theta \neq \theta_k^{\cF}$), the meet of all these matchings is a stable matching of $\cF$ (recall that $\cF$ is a sublattice) where all the rotations $\theta_0^{\cF} \cup \bigcup_{\theta \in \Theta} \theta$ are eliminated. Now let $M$ denote the join (i.e., largest) of all stable matching of $\cF$ in which all the rotations $\theta_0^{\cF} \cup \bigcup_{\theta \in \Theta} \theta$ are eliminated, i.e., such that $\theta_0^{\cF} \cup \bigcup_{\theta \in \Theta} \theta \subset R_M$. Note that $M \in \cF$ since $\cF$ is a sublattice. We claim that $R_M = \theta_0^{\cF} \cup \bigcup_{\theta \in \Theta} \theta$ and hence $M_{\theta_0^{\cF} \cup \bigcup_{\theta \in \Theta} \theta}=M \in \cF$. In fact, if $R_M \neq \theta_0^{\cF} \cup \bigcup_{\theta \in \Theta} \theta$, and since $R_M$ can be decomposed into meta-rotations $\theta_0^{\cF}$ plus a set of proper-meta rotations (as we show in the first point), there exists a proper meta-rotation $\hat{\theta} \subset R_M$ such that $\hat{\theta} \notin \Theta$. Now such meta-rotation cannot be larger than any $\theta \in \Theta$ otherwise it would belong to $\Theta$ as it is upper-closed. Hence, for every $\theta\in \Theta$, there exists a stable matching $M(\theta, \hat{\theta}) \in \cF$ in which $\theta$ are eliminated but not $\hat{\theta}$. The meet $\bigwedge_{\theta\in \Theta} M(\theta, \hat{\theta})\in \cF$ is a stable matching of $\cF$ in which all the rotations $\theta_0^{\cF} \cup \bigcup_{\theta \in \Theta} \theta$ are eliminated but not the rotations $\hat{\theta}$. Hence, the join $M \vee \left(\bigwedge_{\theta\in \Theta} M(\theta, \hat{\theta})\right) \in \cF$ is a stable matching of $\cF$ that is strictly larger than $M$ (as $\hat{\theta}$ are not eliminated in this matching) and where all the rotations $\theta_0^{\cF} \cup \bigcup_{\theta \in \Theta} \theta$ are eliminated contradicting the maximality of $M$.

Now for the third point. Fix $M, M' \in \cF$. Then $M \wedge M', M \vee M' \in \cF$. We show that $\Theta_{M \wedge M'} = \Theta_M \cup \Theta_{M'}$ and $\Theta_{M \vee M'} = \Theta_M \cap \Theta_{M'}$. In fact, $$R_{M \wedge M'} = R_M \cup R_{M'} = \left(\theta_0^{\cF} \cup \bigcup_{\theta \in \Theta_M}\theta\right) \cup \left(\theta_0^{\cF} \cup \bigcup_{\theta \in \Theta_{M'}}\theta\right)=\left(\theta_0^{\cF} \cup \bigcup_{\theta \in \Theta_M \cup \Theta_{M'}}\theta\right).$$ Hence, $\Theta_{M \wedge M'} = \Theta_M \cup \Theta_{M'}$. Similarly, $$R_{M \vee M'} = R_M \cap R_{M'} = \left(\theta_0^{\cF} \cup \bigcup_{\theta \in \Theta_M}\theta\right) \cap \left(\theta_0^{\cF} \cup \bigcup_{\theta \in \Theta_{M'}}\theta\right)=\left(\theta_0^{\cF} \cup \bigcup_{\theta \in \Theta_M \cap \Theta_{M'}}\theta\right).$$ Hence, $\Theta_{M \vee M'} = \Theta_M \cap \Theta_{M'}$.
\hfill $\square$

\section{Minimum-Cut Representability: Omitted Proofs}
\subsection{Proof of Theorem~\ref{thm:main1}}
\label{apx:main1}

Suppose \eqref{problem} is minimum cut representable. In particular, there exists a capacitated digraph $\cD$ over rotations and a source vertex $s$ and a sink vertex $t$ with non-negative (possibly infinite) capacities such that an $s$-$t$ cut $\{s\} \cup R$ has finite value if and only if $R$ is upper-closed and $M_R \in \cal F$. In addition, there exists a constant $C \in \mathbb{R}$ such that the value of a finite value $s$-$t$ cut $\{s\} \cup R$ in $\cD$ is given by $f(M_R)+C$. We show that (i)-(iii) hold.

Replace the capacity of the infinite capacity arcs of $\cD$ with a large positive number $Q$ ($Q > 2 \max_{M \in \cal F} f(M) + 2C$ is enough). For every subset of rotations $R \subset \cR(I)$, let $\cval(R)$ denote the value of the $s$-$t$ cut $\{s\} \cup R$ in the modified digraph.  $\cval(\cdot)$ is submodular over the sets of rotations (see, e.g. \cite{schrijver2003combinatorial}) and $\cval(R) = f(M_R)+C$ for every upper-closed subset of rotations $R$ such that $M_R \in {\cal F}$.

To show (i), consider $M, M' \in \cF$. By submodularity of $\cval(\cdot)$, we have:
$$
    \cval(R_M) + \cval(R_{M'}) \geq \cval(R_M \cup R_{M'}) + \cval( R_{M} \cap R_{M'}).
$$
Since $M, M' \in \cF$, the sets of rotations $R_M$ and $R_{M'}$ correspond to finite cuts in $\cD$. Thus, $\cval(R_M) < Q/2$ and $\cval(R_{M'}) < Q/2$, implying that $\cval(R_M \cup R_{M'}) < Q$ and $\cval(R_{M} \cap R_{M'}) < Q$, and hence that $\{s\} \cup R_M \cup R_{M'}$ and $\{s\} \cup R_{M} \cap R_{M'}$ are finite value $s$-$t$ cuts. Hence, $R_M \cup R_{M'}$ and $R_M \cap R_{M'}$ are upper-closed and their corresponding stable matchings $M_{R_M \cup R_{M'}}$ and $M_{R_M \cap R_{M'}}$ are in $\cF$. Since $M \rightarrow R_M$ is an isomorphism of lattices by Theorem~\ref{thm:representation_with_rotations}, we have that $M_{R_M \cup R_{M'}} = M_{R_{M \wedge M'}} = M \wedge M'$ and $M_{R_M \cap R_{M'}} = M_{R_{M \vee M'}} = M \vee M'$. We deduce that $M \vee M', M \wedge M' \in \mathcal{F}$, proving (i).

To show (ii), consider a pair of distinct proper meta-rotations $\theta \neq \theta' \in \cR_{\sf prop}^{\cF}(I)$. Since $\cval(R) = f(M_R) + C$ for every upper-closed subset of rotations $R$ and $M \rightarrow R_M$ is an isomorphism of lattices, we have $$
f(M^{\theta} \wedge M_{\theta'}) + f(M_{\theta} \wedge M^{\theta'}) =\cval(R_{M^{\theta} \wedge M_{\theta'}}) +  \cval(R_{M_{\theta} \wedge M^{\theta'}}) - 2C =\cval(R_{M^{\theta}} \cup R_{M_{\theta'}}) +  \cval(R_{M_{\theta}} \cup R_{M^{\theta'}}) - 2C.
$$
By definition of $R^\theta$ and $R_{\theta}$, we further have
$$
\cval(R_{M^{\theta}} \cup R_{M_{\theta'}}) +  \cval(R_{M_{\theta}} \cup R_{M^{\theta'}}) =\cval(R^{\theta} \cup R_{\theta'}) +  \cval(R_{\theta} \cup R^{\theta'}) = \cval(R_{\theta} \cup R_{\theta'} \cup \theta) +  \cval(R_{\theta} \cup R_{\theta'} \cup \theta').$$
By submodularity of $\cval(.)$ and since $\theta \neq \theta'$, we deduce
$$
\cval(R_{\theta} \cup R_{\theta'} \cup \theta) +  \cval(R_{\theta} \cup R_{\theta'} \cup \theta') \geq \cval(R_{\theta} \cup R_{\theta'} \cup \theta \cup \theta') +  \cval(R_{\theta} \cup R_{\theta'})=\cval(R^{\theta} \cup R^{\theta'}) + \cval(R_{\theta} \cup R_{\theta'}).
$$

Hence,
$$
f(M^{\theta} \wedge M_{\theta'}) + f(M_{\theta} \wedge M^{\theta'}) \geq  \cval(R^{\theta} \cup R^{\theta'}) + \cval(R_{\theta} \cup R_{\theta'}) - 2C =     f(M^{\theta} \wedge M^{\theta'}) + f(M_{\theta} \wedge M_{\theta'}).
$$

\iffalse 
\begin{align*}
    &f(M^{\theta} \wedge M_{\theta'}) + f(M_{\theta} \wedge M^{\theta'}) \\
    &=\cval(R_{M^{\theta} \wedge M_{\theta'}}) +  \cval(R_{M_{\theta} \wedge M^{\theta'}}) - 2C\\
    &=\cval(R_{M^{\theta}} \cup R_{M_{\theta'}}) +  \cval(R_{M_{\theta}} \cup R_{M^{\theta'}}) - 2C\\
    &= 
    \cval(R^{\theta} \cup R_{\theta'}) + \cval(R_{\theta} \cup R^{\theta'}) - 2C\\
    &=
    \cval(R_{\theta} \cup R_{\theta'} \cup \theta) + \cval(R_{\theta} \cup R_{\theta'} \cup \theta') - 2C\\
    &\geq
    \cval(R_{\theta} \cup R_{\theta'} \cup \theta \cup \theta') + \cval(R_{\theta} \cup R_{\theta'}) - 2C\\
    &=\cval(R_{M^{\theta}} \cup R_{M^{\theta'}}) +  \cval(R_{M_{\theta}} \cup R_{M_{\theta'}}) - 2C\\
    &=\cval(R_{M^{\theta} \wedge M^{\theta'}}) + \cval(R_{M_{\theta} \wedge M_{\theta'}})-2C\\
    &=
    f(M^{\theta} \wedge M^{\theta'}) + f(M_{\theta} \wedge M_{\theta'})
\end{align*}
The first and last equalities follows from the fact that $\cval(R) = f(M_R) + C$ for every upper-closed subset of rotations $R$. The second and second to last equalities hold since $M \rightarrow R_M$ is an isomorphism of lattices. The third and third to last equality follow from the definitions of $R^\theta$ and $R_{\theta}$. The inequality follows by submodularity of $\cval(.)$ and by the fact that $\theta \neq \theta'$. 
\fi 
Thus, $\partial^2 f_{\theta, \theta'} = f(M^{\theta} \wedge M_{\theta'}) + f(M_{\theta} \wedge M^{\theta'}) - f(M_{\theta} \wedge M_{\theta'}) - f(M^{\theta} \wedge M^{\theta'}) \geq 0
$, showing (ii).

Let us show (iii). For $A, B \subset \cR(I) \cup \{s,t\}$, let $u(A,B)$ denote the sum of the capacity of all arcs with tail in $A$ and head in $B$ in ${\cal D}$. Three claims are in order.
\begin{claim}
\label{clm:cut0}
    Let $R \subset \cR(I)$. Then, 
    $
    \cval(R) = u(R \cup \{s\}, \cR(I) \cup \{s,t\}) - u(R \cup \{s\}, R \cup \{s\})
    $.
\end{claim}
{\em Proof.} Immediate from the definition of an $s$-$t$ cut.
\hfill$\square$

\begin{claim}
\label{clm:cut1}
    Let $\theta \in \cR_{\sf prop}^{\cF}(I)$. Then, $\partial f_\theta = u(\theta, \cR(I) \cup \{s,t\}) - u(R_\theta \cup \{s\}, \theta) - u(\theta, R_\theta \cup \{s\}) - u(\theta, \theta).$
\end{claim}
{\em Proof. }
$\partial f_\theta 
    =f(M^{\theta}) - f(M_{\theta})
    = \cval(R_{M^\theta}) - C - \cval(R_{M_\theta}) + C 
    = \cval(R^{\theta}) - \cval(R_{\theta})$, where the second equality follows from the fact that $\cval(R) = f(R)+C$ for every upper-closed subset of rotations $R$ such that $M_R \in {\cal F}$ and the third follows from the definition of $R^\theta$ and $R_\theta$. By Claim~\ref{clm:cut0}, $\cval(R^{\theta}) - \cval(R_{\theta})$ coincides with 
    $$u(R^{\theta} \cup \{s\}, \cR(I) \cup \{s, t\})  - u(R^{\theta} \cup \{s\}, R^{\theta} \cup \{s\}) - u(R_{\theta} \cup \{s\}, \cR(I) \cup \{s, t\}) + u(R_{\theta} \cup \{s\}, R_{\theta} \cup \{s\}),
    $$
    which, since $R^\theta = R_\theta \cup \theta$, is equal to
$u(\theta, \cR(I) \cup \{s,t\}) - u(R_\theta \cup \{s\}, \theta) - u(\theta, R_\theta \cup \{s\}) - u(\theta, \theta)$
as required, concluding the proof. \hfill$\square$ 
\iffalse \begin{align*}
    \partial f_\theta 
    &=f(M^{\theta}) - f(M_{\theta})\\
    &= \cval(R_{M^\theta}) - C - \cval(R_{M_\theta}) + C\\
    &= \cval(R^{\theta}) - \cval(R_{\theta})\\
    &= u(R^{\theta} \cup \{s\}, \cR(I) \cup \{s, t\}) \\
    &\quad - u(R^{\theta} \cup \{s\}, R^{\theta} \cup \{s\})\\
    &\quad - u(R_{\theta} \cup \{s\}, \cR(I) \cup \{s, t\})\\ 
    &\quad - u(R_{\theta} \cup \{s\}, R^{\theta} \cup \{s\})\\
    &=u(\theta, \cR(I) \cup \{s,t\}) - u(R_\theta \cup \{s\}, \theta) \\
    &\quad - u(\theta, R_\theta \cup \{s\}) - u(\theta, \theta).
\end{align*}
The second equality follows from the fact that $\cval(R) = f(R)+C$ for every upper-closed subset of rotations $R$. The third equality follows from the definition of $R^\theta$ and $R_\theta$. The fourth equality follows from Claim~\ref{clm:cut0}. Finally the last equality follows from the fact that $R^\theta = R_\theta \cup \theta$. \fi

\begin{claim}
\label{clm:cut2}
    Let $\theta \neq \theta' \in \cR_{\sf prop}^{\cF}(I)$ be a pair of distinct meta-rotations. Then,
    $$
    \partial^2 f_{\theta, \theta'} = \left\{\begin{matrix}
        u(\theta,\theta') + u(\theta',\theta) & \quad \text{ if $\theta, \theta'$ are incomparable}\\
        0 & \text{ otherwise.}
    \end{matrix}\right.
    $$
\end{claim}
{\em Proof.} 
\underline{If $\theta$ and $\theta'$ are comparable, say without loss of generality $\theta \triangleright^{\cF} \theta'$:} We begin by showing that $R^{\theta} \subset R_{\theta'}$. Suppose not; then, there exists a rotation $\rho \in R^{\theta} \setminus R_{\theta'}$. Let $\overline{\rho}$ denote the meta-rotation to which $\rho$ belongs. Then, $\overline{\rho}$ is a proper meta-rotation (since $\rho$ is not eliminated in $M_{\theta'}$ but it is eliminated in $M^{\theta}$) and $\overline{\rho} \trianglerighteq^{\cF} \theta$ by definition of $R^{\theta}$. Thus, $\overline{\rho} \triangleright^{\cF} \theta'$ and  $\overline{\rho} \subset R^{\theta'}$. We conclude  that $\overline{\rho} \subset R_{\theta'} = R^{\theta'} \setminus \{\theta'\}$ since $\overline{\rho} \neq \theta'$ (recall that $\overline{\rho} \subset R^{\theta'}$). This implies that $\rho \in R_{\theta'}$, a contradiction. 

Now given $R^{\theta} \subset R_{\theta'}$,
$M^\theta \wedge M_{\theta'}  = M_{R^{\theta}} \wedge M_{R_{\theta'}}  = M_{R^{\theta} \cup R_{\theta'}} = M_{R_{\theta'}},$
where the second equality holds by the fact that $R \rightarrow M_R$ is an isomorphism of lattices. Since, $R_\theta \subset R^{\theta} \subset R_{\theta'} \subset R^{\theta'}$, we can apply similar arguments and deduce $M_\theta \wedge M^{\theta'} = M_{R^{\theta'}}$, 
$M_\theta \wedge M_{\theta'} = M_{R_{\theta'}}$, 
and $M^\theta \wedge M^{\theta'}= M_{R^{\theta'}}$. 
Hence, $\partial^2 f_{\theta, \theta'} = 0$.

\underline{If $\theta$ and $\theta'$ are incomparable:} Let $R=R_\theta \cup R_{\theta'}$.  
First, we show that the sets $R$, $\theta$, and $\theta'$ are disjoint: We already know that $\theta$ and $\theta'$ are disjoint, $\theta$ and $R_\theta$ are disjoint, and $\theta'$ and $R_{\theta'}$ are disjoint. Hence, we only need to check that $\theta$ and $R_{\theta'}$ are disjoint and $\theta$ and $R_{\theta'}$ are disjoint. We only show the first without loss of generality. Suppose $\theta$ and $R_{\theta'}$ are not disjoint. Then there exists a rotation $\rho \in R_{\theta} \cap \theta'$. By definition of $R_{\theta}$, the meta-rotation to which $\rho$ belongs, which is $\theta'$ here, is a subset of $R_{\theta}$. Then, because $R_{\theta}$ is the union of $\theta_{0}^{\cF}$ and all and only the proper-meta rotations that are larger than $\theta$ with respect to $\trianglerighteq^{\cF}$, it must be the case that $\theta' \trianglerighteq^{\cF} \theta$ contradicting the fact that $\theta'$ and $\theta$ are incomparable. 

From the definition of $M^\theta, M_{\theta}, R^{\theta}$ and $R_{\theta}$ and since $R \rightarrow M_R$ is an isomorphism of lattices, we have:
 $$\partial^2 f_{\theta, \theta'} 
= \cval(R^{\theta} \cup R_{\theta'})  + \cval(R_{\theta} \cup R^{\theta'}) - \cval(R_{\theta} \cup R_{\theta'}) - \cval(R^{\theta} \cup R^{\theta'}) = \cval(R \cup {\theta}) + \cval(R \cup \theta') - \cval(R) - \cval(R \cup \theta \cup \theta').$$
By Claim~\ref{clm:cut0}, we have 
$\cval(R \cup {\theta})= u(R \cup \theta \cup \{s\}, \cR(I) \cup \{s,t\}) - u(R \cup \theta \cup \{s\})$ and we can similarly rewrite $\cval(R \cup \theta')$, %= u(R \cup \theta' \cup \{s\}, \cR(I) \cup \{s,t\})- u(R \cup \theta' \cup \{s\}, R \cup \theta' \cup \{s\})$,
$\cval(R)$% = \cval(R \cup \theta \cup \theta')=u(R \cup \{s\}, \cR(I) \cup \{s,t\}) - u(R \cup \{s\}, R \cup \{s\})$ and
, and $\cval(R \cup \theta \cup \theta')$. %= u(R \cup \theta \cup \theta' \cup \{s\}, \cR(I) \cup \{s,t\}) - u(R \cup \theta \cup \theta' \cup \{s\}, R \cup \theta \cup \theta' \cup \{s\})$. 
Using the fact that $R$, $\theta$, and $\theta'$ are disjoint and by noticing that for every subsets of rotations $A, B, C$ such that $A \cap B = \emptyset$ it holds that $u(A \cup B,  C)=u(A, C) + u(B, C)$, simple calculations lead to  $\partial^2 f_{\theta, \theta'}=u(\theta, \theta') + u(\theta', \theta)$, as required. \hfill$\square$

We now prove (iii). In particular, for every $M \in \cF$, from Claim~\ref{clm:cut0} we can rewrite $f(M) - f(M_0^{\cF})$ as 
{\small{$$u(R_M \cup \{s\}, \cR(I) \cup \{s,t\}) 
    - u(R_M \cup \{s\}, R_M \cup \{s\}) 
    - u(\theta_0^{\cF} \cup \{s\}, \cR(I) \cup \{s,t\}) 
    + u(\theta_0^{\cF} \cup \{s\}, \theta_0^{\cF} \cup \{s\}).$$}}
    Since $R_M = \theta_0^{\cF} \cup \bigcup_{\theta \in \Theta_M} \theta$ (Theorem~\ref{thm:representation_meta_rotations}), we can further write  
$$f(M) - f(M_0^{\cF})=\sum_{\theta \in \Theta_M} u(\theta, \cR(I) \cup \{s,t\}) - \sum_{\theta, \theta' \in \Theta_M} u(\theta , \theta')-\sum_{\theta \in \Theta_M} u(\theta_0^{\cF} \cup \{s\}, \theta)- \sum_{\theta \in \Theta_M}u(\theta, \theta_0^{\cF} \cup \{s\}).$$
By separating pairs $\theta, \theta'$ that are comparable from those that are not, $- \sum_{\theta, \theta' \in \Theta_M} u(\theta , \theta')$ is equal to
 {\small{$$ -\sum_{\substack{\theta, \theta' \in \Theta_M\\\theta, \theta' \text{ incomparable}}} u(\theta, \theta') - \sum_{\substack{\theta, \theta' \in \Theta_M\\\theta' > \theta}} u(\theta, \theta') - \sum_{\substack{\theta, \theta' \in \Theta_M\\\theta' > \theta}} u(\theta', \theta) - \sum_{\theta \in \Theta_M} u(\theta, \theta).$$}}
Moreover, $R_\theta = \theta_0^{\cF} \cup \bigcup_{\substack{\theta' \in \cR^{\cF}_{\sf prop}(I)\\\theta'> \theta}} \theta$ implies $
R_\theta = \theta_0^{\cF} \cup \bigcup_{\substack{\theta' \in \Theta_M\\\theta'> \theta}} \theta
$ since $\Theta_M$ is upper-closed. Thus, 
{\small{$$-\sum_{\substack{\theta, \theta' \in \Theta_M\\\theta' > \theta}} u(\theta, \theta') -\sum_{\substack{\theta, \theta' \in \Theta_M\\\theta' > \theta}} u(\theta', \theta) 
    - \sum_{\theta \in \Theta_M} u(\theta_0^{\cF} \cup \{s\}, \theta) - \sum_{\theta \in \Theta_M}u(\theta, \theta_0^{\cF} \cup \{s\})= -\sum_{\theta \in \Theta_M} u(\theta, R_\theta \cup \{s\}) - \sum_{\theta \in \Theta_M} u(R_\theta \cup \{s\}, \theta).$$}}
Putting all together, we deduce  
$
f(M) - f(M_0^{\cF})=\sum_{\theta \in \Theta_M} \partial f_\theta - \frac{1}{2} \sum_{\theta \neq \theta' \in \Theta_M} \partial^2 f_{\theta, \theta'}.$ This concludes the proof that minimum cut representability implies that (i)-(iii) hold.

We now prove the converse. Suppose (i)-(iii) hold. That~\eqref{problem} is minimum cut representable immediately follows from the following lemma,
\begin{lemma}
\label{lem:onlyif}
    Suppose (i)-(iii) hold, then,
    \begin{enumerate}
        \item The cut digraph $\cD^{f, \cF}$ has non-negative (possibly infinite) capacities.
        \item A cut $\{s\} \cup R$ has finite capacity in $\cD^{f, \cF}$ if and only if $R$ is upper-closed and $M_R \in \cF$.
        \item The value of a finite capacity $s$-$t$ cut $\{s\} \cup R$ in $\cD^{f, \cF}$ is $f(M_R) - f(M_0^{\cF}) + \gamma \cdot |\cR_{\sf prop}^{\cF}(I)|$. Note that $- f(M_0^{\cF}) + \gamma \cdot |\cR_{\sf prop}^{\cF}(I)|$ is a constant independent of $R$.
    \end{enumerate}
\end{lemma}
{\em Proof.}
Suppose (i)-(iii) holds. We have,\\
1. Follows immediately from (ii).\\
2. Recall that the infinite capacity arcs of $\cD^{f, \cF}$ are: $(*)$ Arcs between every pair of rotations in a meta-rotation. $(**)$  An arc from $s$ to $\widehat{\theta_0^{\cF}}$ and an arc from $\widehat{\theta_{z}^{\cF}}$ to $t$, and finally, $(***)$ an arc from $\widehat{\theta}$ to $\widehat{\theta'}$ for every $\theta \triangleleft^{\cF} \theta'$. These arcs insure that a cut $\{s\} \cup R$ has finite capacity if and only if $(a)$ when $R$ contains a rotation then it contains all the meta-rotation to which it belongs and hence that $R$ can be partitioned into meta-rotations, $(b)$ $R$ includes every rotation in $\theta_0^{\cF}$ and excludes every rotation in $\theta_{z}^{\cF}$, and $(c)$ the proper meta-rotations whose rotations are included in $R$ form an upper-closed subset of meta-rotations. By Theorem~\ref{thm:representation_meta_rotations}, $(a)-(c)$ are equivalent to $R$ being upper-closed and corresponding to a stable matching $M_R \in \cF$.\\
3. For every $A, B \subset \cR(I)\cup \{s,t\}$, let $u(A, B)$ denote the sum of the capacities of arcs going from $A$ to $B$ in $\cD^{f, \cF}$. Consider a finite capacity $s$-$t$ cut $\{s\} \cup R$. The value of this $s$-$t$ cut is given by,
$$
        u(R \cup \{s\}, \cR(I) \cup \{s,t\}) - u(R \cup \{s\}, R \cup \{s\})
        = u(R, \cR(I) \cup \{t\}) + u(\{s\}, \cR(I)) - u(\{s\}, R) - u(R, R),
$$
where the equality holds since no arc is entering $s$ and no arc from $s$ to $t$ belong to $\cD^{f, \cF}$. Using moreover that $R = \theta_0^{\cF} \cup \bigcup_{\theta \in \Theta_{M_R}} \theta$ (Theorem~\ref{thm:representation_meta_rotations}) and the definition of $\cD^{f, \cF}$, the value of the cut is equal to $$\sum_{\theta \in \Theta_{M_R}} u(\theta, \cR(I)) + \sum_{\theta \in \Theta_{M_R}} u(\theta, \{t\}) + \gamma |{\textstyle \cR^{\cF}_{\sf prop}}(I)| - \gamma |\Theta_{M_R}| - \sum_{\theta, \theta' \in \Theta_{M_R}} u(\theta, \theta'),$$
which is further equal to,
$$\sum_{\theta \in \Theta_{M_R}} u(\theta, \cR(I) \setminus \theta) + \sum_{\theta \in \Theta_{M_R}} u(\theta, \{t\}) + \gamma |{\textstyle \cR^{\cF}_{\sf prop}}(I)| - \gamma |\Theta_{M_R}| - \sum_{\theta \neq \theta' \in \Theta_{M_R}} u(\theta, \theta').$$
Using again the definition of $\cD^{f, \cF}$, we have that, for $\theta \in \Theta_{M_R}$, 
$$u(\theta,\cR(I) \setminus \theta)+u(\theta,\{t\})=\sum_{\theta' \in \cR^{\cF}_{\sf prop}(I) \setminus \{\theta\}} \frac{1}{2} \partial^2 f_{\theta, \theta'}+ \partial f_{\theta} - \frac{1}{2} \sum_{\theta' \in \cR^{\cF}_{\sf prop}(I) \setminus\{\theta\}} \partial^2 f_{\theta, \theta'}+ \gamma,$$
% and,
% $$
% u(\theta, \theta') = \frac{1}{2} \partial^2 f_{\theta, \theta'}.
% $$
Putting everything together, we have that the value of the cut is $f(M) - f(M_0^{\cF}) + \gamma |{\textstyle \cR^{\cF}_{\sf prop}}(I)|$. 
\hfill $\square$

Lemma~\ref{lem:onlyif} implies that \eqref{problem} is minimum cut representable with witness digraph $\cD^{f, \cF}$. This also implies the second part of the theorem, namely, that whenever \eqref{problem} is minimum cut representable (in which case (i)-(iii) hold) as shown before, then $\cD^{f, \cF}$ is a witness digraph.

\subsection{Proof of Theorem~\ref{thm:linear-to-second-order}}

Consider an optimization problem $\eqref{problem}$. We show that $\eqref{problem}$ is an instance of (MWSM) if and only if it is minimum cut representable with vanishing second order differentials. Let us begin by showing some preliminary results.

Let $\cS(\cF) := \{ab 
\in A \times B| \exists M \in \cF, ab \in M\}$ denote the set of stable pairs, that is, pairs belonging to some stable matching of $\cF$. Note that for every proper meta-rotation $\theta \in \cR^{\cF}_{\sf prop}(I)$, there exists a subset of pairs $p_{\sf in}(\theta) \subset \cS(\cF)$ such that after all rotations from $\theta$ have been eliminated, $p_{\sf in}(\theta)$ enter the stable matching. Similarly, there exists a subset of pairs $p_{\sf out}(\theta) \subset \cS(\cF)$ such that after all rotations from $\theta$ have been eliminated, $p_{\sf out}(\theta)$ leave the stable matching. 

Consider a pair $ab \in \cS(\cF)$. Four cases are possible: (i) $ab \notin M_0^{\cF}$ and $ab \notin M_z^{\cF}$, where $M_z^{\cF}$ is the meet of all stable matchings of $\cF$. In this case, there exists a unique proper meta-rotation $\theta$ for which $ab \in p_{\sf in}(\theta)$ (resp. $ab \in p_{\sf out}(\theta)$), we denote such proper meta-rotation by $\theta_{\sf in}(ab)$ (resp. $\theta_{\sf out}(ab)$). (ii) $ab \in M_0^{\cF}$ and $ab \notin M_z^{\cF}$. In this case, there exists a unique proper meta-rotation $\theta$ for which $ab \in p_{\sf out}(\theta)$ that we denote by $\theta_{\sf out}(ab)$ and no proper meta-rotation for which $ab \in p_{\sf in}(\theta)$ and we write $\theta_{\sf in}(ab) = \emptyset$. (iii) $ab \notin M_0^{\cF}$ and $ab \in M_z^{\cF}$. In this case, there exists a unique proper meta-rotation $\theta$ for which $ab \in p_{\sf in}(\theta)$ that we denote by $\theta_{\sf in}(ab)$ and no proper meta-rotation for which $ab \in p_{\sf out}(\theta)$ and we write $\theta_{\sf out}(ab) = \emptyset$. (iv) $ab \in M_0^{\cF}$ and $ab \in M_z^{\cF}$. In this case, there exists a no proper meta-rotation $\theta$ for which $ab \in p_{\sf in}(\theta)$ and no proper meta-rotation for which $ab \in p_{\sf out}(\theta)$ and we write $\theta_{\sf in}(ab) = \emptyset$ and $\theta_{\sf out}(ab) = \emptyset$.

Consider now the matrix $K$ whose rows are the stable pairs $\cS(\cF)$, whose columns are the proper meta rotations $\cR^{\cF}_{\sf prop}(I)$, and such that for every $\theta \in \cR^{\cF}_{\sf prop}(I)$, column $\theta$ has $1$ in the rows $p_{\sf in}(\theta)$, $-1$ in the rows $p_{\sf out}(\theta)$ and $0$ everywhere else. For a stable matching $M \in \cF$, let $\chi(M) \in \{0,1\}^{|\cS(\cF)|}$ denote the indicator vector such that $\chi(M)_{ab} = 1$ if and only if $ab \in M$. Similarly, let $\xi(M) \in \{0,1\}^{|\cR^{\cF}_{\sf prop}(I)|}$ denote the indicator vector such that $\xi(M)_{\theta} = 1$ if and only if $\theta \in \Theta_M$. Note that for every $ab \in \cS(\cF)$, we have: (i) $(K\xi(M))_{ab} = 1$ if $\theta_{\sf in}(ab) \neq \emptyset$ and $\theta_{\sf in}(ab)$ is eliminated in $M$ and either $\theta_{\sf out}(ab)$ is not eliminated in $M$ or $\theta_{\sf out}(ab) = \emptyset$. (ii) $(K\xi(M))_{ab} = -1$ if $\theta_{\sf in}(ab) = \emptyset$ and $\theta_{\sf out}(ab) \neq \emptyset$ and $\theta_{\sf out}(ab)$ is eliminated in $M$. (iii) $(K\xi(M))_{ab} = 0$ otherwise. This implies in particular that 
\begin{align}\label{eq:chi-K}
\chi(M) = \chi(M_0^{\cF}) + K \xi(M),
\end{align}
as (i) $\chi(M)_{ab} - \chi(M_0^{\cF})_{ab} = 1-0 = 1$ if $\theta_{\sf in}(ab) \neq \emptyset$ and $\theta_{\sf in}(ab)$ is eliminated in $M$ and $\theta_{\sf out}(ab)$ is not eliminated in $M$ or $\theta_{\sf out}(ab) = \emptyset$; (ii) $\chi(M)_{ab} - \chi(M_0^{\cF})_{ab} = 0 - 1 = -1$ if $\theta_{\sf in}(ab) = \emptyset$ and $\theta_{\sf out}(ab) \neq \emptyset$ and $\theta_{\sf out}(ab)$ is eliminated in $M$; (iii) clearly $\chi(M)_{ab} - \chi(M_0^{\cF})_{ab} = 0$ in the remainder of the cases.

We now claim that $K$ is full column rank. Suppose not, then there exists a vector $\xi \in \mathbb{R}^{|\cR^{\cF}_{\sf prop}(I)|} \setminus \{{\bf 0}\}$ such that $K \xi = {\bf 0}$. Let $\theta^* \in \cR^{\cF}_{\sf prop}(I)$ be a minimal proper meta-rotation w.r.t $\trianglerighteq^{\cF}$ among all proper meta-rotations $\theta$ such that $\xi_{\theta} \neq 0$, and let $ab \in p_{\sf in}(\theta^*)$. Then, by definition $(K\xi)_{ab} = \xi_{\theta_{\sf in}(ab)} - \xi_{\theta_{\sf out}(ab)}$ if $\theta_{\sf out}(ab) \neq \emptyset$ or $(K\xi)_{ab} = \xi_{\theta_{\sf in}(ab)}$ otherwise. If $\theta_{\sf out}(ab) \neq \emptyset$, it must be the case that $\zeta_{\theta_{\sf out}(ab)} = 0$ by minimality of $\theta^*$ and the fact that $\theta_{\sf out}(ab) \triangleleft^F \theta_{\sf in}(ab) = \theta^*$. Hence, in both cases $(K\xi)_{ab} = \zeta_{\theta_{\sf in}(ab)} = \zeta_{\theta^*} > 0$; a contradiction.

Let $\widehat{K}$ denote the matrix $K$ augmented to a basis of $\mathbb{R}^{|\cS(\cF)| \times |\cS(\cF)|}$ by adding $|\cS(\cF)| - |\cR^{\cF}_{\sf prop}(I)|$ columns at the end. Similarly, for every vector $\xi \in \mathbb{R}^{|\cR^{\cF}_{\sf prop}(I)|}$, let $\widehat{\xi} \in \mathbb{R}^{|\cS(\cF)|}$ denote the vector $\xi$ augmented to a vector of $\mathbb{R}^{|\cS(\cF)|}$ by adding $|\cS(\cF)| - |\cR^{\cF}_{\sf prop}(I)|$ zeros at the end. Then, $\widehat{K}$ is non singular and for every $M \in \cF$, we have
$
\chi(M) = \chi(M_0^{\cF}) + \widehat{K} \widehat{\xi(M)}.
$

% for every $\xi \in \cR^{\cF}_{\sf prop}(I)$, there exists a stable matching $M_\x$ such that $\chi(M_{\xi}) = K \xi$.

We now prove the theorem. First, suppose that problem $\eqref{problem}$ is minimum cut representable and $\partial^2_{\theta, \theta'}f = 0$ for every $\theta 
 \neq \theta' \in \cR^{\cF}_{\sf prop}(I)$. We show that \eqref{problem} is an instance of (MWSM). In particular, we show the existence of weights $w_{ab}$ for every $a \in A$ and $b \in B$ such that $f(M) = \sum_{ab \in M} w_{ab}$ for every $M \in \cF$. Recall that the number of matched pairs is the same in every stable matching of $\cF$ (e.g. \cite{gusfieldbook}), let $m$ denote this number. Let $ {\bf \partial f} = [\partial_\theta f]_{\theta \in \cR^{\cF}_{\sf prop}(I)}$ denote the vector of all first order differentials. Now, since $\partial^2_{\theta, \theta'}f = 0$ for every $\theta \neq \theta' \in \cR^{\cF}_{\sf prop}(I)$, it holds that, 
{\small{$$
f(M) 
    = f(M_0^{\cF}) + \sum_{\theta \in \Theta_M} \partial f_\theta
    = f(M_0^{\cF}) +  {\bf \partial f}^T \xi(M)
    = f(M_0^{\cF}) +  \widehat{{\bf \partial f}}^T \widehat{\xi(M)}=f(M_0^{\cF}) +  \widehat{{\bf \partial f}}^T \widehat{K}^{-1}(\chi(M)-\chi(M_0^{\cF})),
 $$}}
and by letting ${\bf 1}$ denotes the column vector of all ones we conclude that $f$ is linear by writing, 
$$
f(M)=\left(\frac{f(M_0^{\cF}) - \widehat{{\bf \partial f}}^T \widehat{K}^{-1}\chi(M_0^{\cF})}{m} {\bf 1} +  \widehat{{\bf \partial f}}^T \widehat{K}^{-1}\right) \chi(M).
$$
\iffalse \begin{align*}
    &f(M) \\
    =& f(M_0^{\cF}) + \sum_{\theta \in \Theta_M} \partial f_\theta\\
    =& f(M_0^{\cF}) +  {\bf \partial f}^T \xi(M)\\
    =& f(M_0^{\cF}) +  \widehat{{\bf \partial f}}^T \widehat{\xi(M)}\\
    =& f(M_0^{\cF}) +  \widehat{{\bf \partial f}}^T \widehat{K}^{-1}(\chi(M)-\chi(M_0^{\cF}))\\
    =& \left(\frac{f(M_0^{\cF}) - \widehat{{\bf \partial f}}^T \widehat{K}^{-1}\chi(M_0^{\cF})}{m} {\bf 1} +  \widehat{{\bf \partial f}}^T \widehat{K}^{-1}\right) \chi(M)
\end{align*}
for every $M \in \cF$, where ${\bf 1}$ denotes the column vector of all ones. Hence, the linearity of $f$.
\fi 

For the reverse direction, suppose that $\eqref{problem}$ is an instance of (MWSM), i.e., there exists weights $w_{ab}$ for every $a \in A$ and $b \in B$ such that $f(M) = \sum_{ab \in M} w_{ab}$ for every $M \in \cF$. Clearly $\eqref{problem}$ is minimum cut representable. To conclude, we show that $\partial^2_{\theta, \theta'}f = 0$ for every $\theta \neq \theta' \in \cR^{\cF}_{\sf prop}(I)$. Let ${\bf w} = [w_{ab}]_{ab \in \cS(\cF)}$ denote the vector of all weights. For every $M \in {\cal F}$, we have $f(M)={\bf w}^T\chi(M)$. Thus, for every $\theta \neq \theta' \in \cR^{\cF}_{\sf prop}(I)$, we have,
$\partial^2 f_{\theta, \theta'} 
    = {\bf w}^T\chi(M^{\theta} \wedge M_{\theta'}) + {\bf w}^T\chi(M_{\theta} \wedge M^{\theta'}) - {\bf w}^T\chi(M_{\theta} \wedge M_{\theta'}) - {\bf w}^T\chi(M^{\theta} \wedge M^{\theta'}).$
Using~\eqref{eq:chi-K} we write
$\partial^2 f_{\theta, \theta'}= {\bf w}^TK\left(\xi(M^{\theta} \wedge M_{\theta'}) + \xi(M_{\theta} \wedge M^{\theta'})\right. \left.- \xi(M_{\theta} \wedge M_{\theta'}) - \xi(M^{\theta} \wedge M^{\theta'})\right)
$
and by letting ${\bf 1}(\Theta)$ be the indicator vector for a subset $\Theta \subset \cR^{\cF}_{\sf prop}(I)$, we conclude
$$
\partial^2 f_{\theta, \theta'} 
    = {\bf w}^TK\left({\bf 1}(\Theta^{\theta} \cup \Theta_{\theta'}) + {\bf 1}(\Theta_{\theta} \cup \Theta^{\theta'}) \right. \left.- {\bf 1}(\Theta_{\theta} \cup \Theta_{\theta'}) - {\bf 1}(\Theta^{\theta} \cup \Theta^{\theta'})\right)= 0.
$$
\iffalse \begin{align*}
    \partial^2 f_{\theta, \theta'} 
    &= f(M^{\theta} \wedge M_{\theta'}) + f(M_{\theta} \wedge M^{\theta'}) \\
    &\quad - f(M_{\theta} \wedge M_{\theta'}) - f(M^{\theta} \wedge M^{\theta'})\\
    &= {\bf w}^T\chi(M^{\theta} \wedge M_{\theta'}) + {\bf w}^T\chi(M_{\theta} \wedge M^{\theta'}) \\
    &\quad - {\bf w}^T\chi(M_{\theta} \wedge M_{\theta'}) - {\bf w}^T\chi(M^{\theta} \wedge M^{\theta'})\\
    &= {\bf w}^TK\left(\xi(M^{\theta} \wedge M_{\theta'}) + \xi(M_{\theta} \wedge M^{\theta'})\right.\\
    &\quad \left.- \xi(M_{\theta} \wedge M_{\theta'}) - \xi(M^{\theta} \wedge M^{\theta'})\right)\\
    &= {\bf w}^TK\left({\bf 1}(\Theta^{\theta} \cup \Theta_{\theta'}) + {\bf 1}(\Theta_{\theta} \cup \Theta^{\theta'}) \right.\\
    &\quad \left.- {\bf 1}(\Theta_{\theta} \cup \Theta_{\theta'}) - {\bf 1}(\Theta^{\theta} \cup \Theta^{\theta'})\right)\\
    &= 0
\end{align*} 
where ${\bf 1}(\Theta)$ is the indicator vector of the subset $\Theta \subset \cR^{\cF}_{\sf prop}(I)$.
\fi 

\subsection{Proof of Lemma~\ref{lem:if-mincut-then-lattice}}

%\afcomment{this is checked.. lets use colors for any new modifications..}

Let~\eqref{problem} be minimum-cut representable with witness digraph $\cD=\cD(\cR(I) \cup \{s,t\}, A)$ and $M',M'' \in {\cal F}$ achieve the minimum of $f$ in ${\cal F}$. We show that $M'\vee M''$, $M' \wedge M''$ also achieve the minimum of $f$ on ${\cal F}$. 

By definition of minimum cut representability, $\{s\} \cup R_{M'}, \{s\} \cup R_{M''}$ are minimum cuts of $\cD$. Then also their intersection $\{s\} \cup (R_{M'}\cap R_{M''})$ is (see, e.g.,~\cite{schrijver2003combinatorial}). Let $\overline R=R_{M'} \cap R_{M''}$. Since ${\cal F}$ is a lattice and $M \rightarrow R_M$ is an isomorphism of lattices, $\overline M=M' \vee M''$ and $\overline M \in {\cal F}$, with $R_{\overline M}=\overline R$. Again by minimum cut representability of~\eqref{problem}, $\overline M$ achieves the minimum of $f$ on ${\cal F}$. A similar argument shows that $M' \wedge M''$ also achieves the minimum of $f$ on ${\cal F}$.

\section{Matching Siblings: Omitted Proofs}\label{ap:siblings}

\subsection{Some Additional Definitions and Facts on School Matching Instances}

We next recall known properties and show new ones, that we will use in our matching siblings applications. Throughout this subsection, fix a school matching instance $I=(A,B,\succ,q)$. 
The first result is known as \emph{opposition of interest property}; its proof can found in \citep[Lemma 1.6.6]{gusfieldbook}. % \afcomment{Lemma 1.6.6 of \cite{gusfieldbook}}. 

For a rotation $\rho\in\cR(I)$, define $A(\rho)=\{a\in A:\exists b\in B,ab\in\rho_+\}$ to be the set of students that appear in some pair of $\rho_+$ (or equivalently $\rho_-$ by definition). Similarly, define $B(\rho)=\{b\in B:\exists a\in A,ab\in\rho_+\}$. 
 For $b \in B$ and sets $A',A''\subseteq A$, we say that \emph{$b$ prefers $A'$ to $A''$} if $b$ prefers all students in $A'$ to all students in $A'' \setminus A'$, i.e., $\hat{a}' \succ_b \hat{a}''$ for all $\hat{a}'\in A', \hat{a}'' \in A''\setminus A'$. The opposition of interest property states that given two stable matchings $M$ and $M'$ and a pair $ab \in M \backslash M'$ then $a$ prefers $M(a)$ to $M'(a)$ if and only if $b$ prefers $M'(b)$ to $M(b)$.

\begin{lemma}[Opposition of Interest]\label{lem:prefer-all-or-none}
Let $M,M' \in {\cal S}(I)$, and let $ab \in M \backslash M'$. Then,  \begin{enumerate}
    \item\label{it:opposition:1} either $M(a) \succ_a M'(a)$, and $b$ prefers $M'(b)$ to $M(b)$;% \sout{all students in $M'(b)$ to all students in $M(b)\setminus M'(b)$ (i.e., $\hat{a}' \succ_b \hat{a}$ for all $\hat{a}'\in M'(b), \hat{a} \in M(b)\setminus M'(b)$);} 
    \item\label{it:opposition:2} or $M'(a)\succ_a M(a)$ and $b$ prefers $M(b)$ to $M'(b)$.% \sout{all students in $M(b)$ to all students in $M'(b)\setminus M(b)$ (i.e., $\hat{a} \succ_b \hat{a}'$ for all $\hat{a}\in M(b), \hat{a}' \in M'(b)\setminus M(b)$).}
    \end{enumerate}
% \sout{When~\eqref{it:opposition:1} occurs, we say that $a$ prefers $M$ to $M'$ and $b$ prefers $M'$ to $M$.}
\end{lemma}

The next lemma shows that any two rotations involving the same student, or the same school, are comparable. While the statement for students is probably known and reproved here for completeness, the one about schools is, to the best of our knowledge, new, and less obvious.

\begin{lemma}[Comparability of Intersecting Rotations] \label{lem:dominace-onemany}
    Let $\rho\neq\rho' \in {\cal R}(I)$ be two rotations. Then,
    \begin{enumerate}
        \item If there is a student $a\in A$ such that $a\in A(\rho)\cap A(\rho')$, then $\rho$ and $\rho'$ are comparable. Moreover, if $\rho_+(a)$ denotes the partner of $a$ in $\rho_+$ and ${\rho'}_+(a)$ denotes the partner of $a$ in ${\rho'}_+$, then $\rho\triangleright\rho'$ if and only if $\rho_+(a) \succ_a {\rho'}_+(a)$.
        \item If there is a school $b\in B$ such that $b\in B(\rho)\cap B(\rho')$, then $\rho$ and $\rho'$ are comparable. Moreover, if $\rho_+(b)$ denotes the set of partners of $b$ in $\rho_+$ and $\rho'_+(b)$ denotes the set of partners of $b$ in ${\rho'}_+$, then $\rho\triangleright\rho'$ if and only if $a' \succ_b a$ for every $a \in \rho_+(b)$ and $a' \in \rho'_+(b)$.
    \end{enumerate}
\end{lemma}
\emph{Proof.} 
1.~Let $a\in A(\rho)\cap A(\rho')$ for some $\rho\neq \rho' \in \cR(I)$ and suppose without loss of generality that $ab\in\rho_+$ and $ab'\in\rho'_+$. Assume $b\succ_a b'$. We show that $\rho \triangleright \rho'$. Let $M\in{\cal S}(I)$ such that $\rho'$ is exposed in $M$. In particular, we have $M(a)=b'$. Let $(M_0,M_1,\dots,M_k)$ be a sequence of stable matchings obtained by eliminating the rotations from $R_M$ when they are exposed until $M_k=M$. Clearly $M_0(a)\succeq_a b$ and by assumption $b'=M_k(a)\prec_a b$. Consider the first stable matching $M_i$ in the sequence such that $M_i(a)\prec_a b$, and let $M_i=M_{i-1}/\rho''$. We claim that $\rho=\rho''$. In fact, by Corollary~\ref{cor:all-rotations}, starting from $M_0$, there exists a sequence $M_0, M_1', \dots, M_{|\cR(I)|}'=M_z$ 
from the student optimal matching $M_0$ to the student-pessimal matching $M_z$ where all rotations are eliminated sequentially. 
In such sequence, if $\rho$ is eliminated before $\rho''$, then after eliminating $\rho$, $a$ is assigned to $\rho_-(a)\prec_a b$ and remains assigned to a school less preferred than $b$ in the remainder of the sequence. By definition $\rho''_+(a)\succeq_a b$. Hence, it is not possible to eliminate $\rho''$ in the remaining sequence. Similarly if $\rho''$ appears earlier, we run into the same issue. Thus $\rho=\rho''$. Since $M$ was an arbitrary matching where $\rho'$ is exposed, we conclude that $\rho\triangleright\rho'$. If $b'\succ_a b$, we deduce $\rho'\triangleright\rho$ by a similar argument.

2.~Now suppose that $b\in B(\rho)\cap B(\rho')$. Recall that $M_\rho=M_0/ R_\rho$, $M^\rho=M_0/R^\rho$ where $R_\rho=\{\tilde{\rho}\in\cR(I):\tilde{\rho}\triangleright\rho\}$ and $R^\rho=\{\tilde{\rho}\in\cR(I):\tilde{\rho}\trianglerighteq\rho\}$ (similarly for $M_{\rho'}$ and $M^{\rho'}$). Consider the matching $M=M_{\rho}\wedge M_{\rho'}=M_0/(R_\rho\cup R_{\rho'})$. If $\rho,\rho'$ are not comparable then we have $\rho'\notin R_\rho$ and $\rho\notin R_{\rho'}$, and thus neither $\rho,\rho'$ is eliminated in $M$. Let $M_1=M^\rho\wedge M_{\rho'}=M_0/(R_\rho\cup R_{\rho'}\cup\{\rho\})$ and $M_1'=M_\rho\wedge M^{\rho'}=M_0/(R_\rho\cup R_{\rho'}\cup\{\rho'\})$. 
By Lemma~\ref{lem:prefer-all-or-none}, $b$ weakly prefers one of $M_1(b)$, $M'_1(b)$ to the other. 
Without loss of generality suppose $b$ weakly prefers $M_1(b)$ to $M'_1(b)$. By definition $M_1(b)=\bigl(M(b)\setminus\rho_+(b)\bigr)\cup\rho_-(b)$ and $M'_1(b)=\bigl(M(b)\setminus\rho'_+(b)\bigr)\cup\rho'_-(b)$. By part 1, since $\rho$ and $\rho'$ are not comparable, we have $A(\rho)\cap A(\rho')=\emptyset$ and thus $\rho_+(b),\rho_-(b),\rho'_+(b),\rho'_-(b)$ are pairwise disjoint subsets of students. Therefore, $b$ strictly prefers $M_1(b)$ to $M'_1(b)$, and thus $b$ prefers any student in $M_1(b)\setminus M'_1(b)=\rho_-(b)\cup\rho'_+(b)$ to any student in $M'_1(b)\setminus M_1(b)=\rho'_-(b)\cup\rho_+(b)$. In particular, for $a\in \rho'_+(b)$ and $a'\in\rho'_-(b)$ it holds $a\succ_b a'$. However, $b$ prefers $\rho'_-(b)$ to $\rho'_+(b)$ by definition of rotation, a contradiction. Thus $\rho$ and $\rho'$ are comparable. 

Now, having established that $\rho$ and $\rho'$ are comparable, we show the remaining results. If $\rho\triangleright\rho'$, then we first claim that $\rho_+(b)\cap \rho'_+(b)=\emptyset$. In fact, if there exists $a\in \rho_+(b)\cap \rho'_+(b)$, then $ab\in\rho_+$ and $ab\in \rho'_+$. By part 1, we have $b=\rho_+(a)\succ_a\rho'_+(a)=b$, a contradiction. Now, let $M\in\cS(I)$ be a stable matching such that $\rho'$ is exposed in $M$. Then, $\rho'_+(b)\subseteq M(b)$. We claim that $b$ prefers any student in $M(b)$ to any student in $\rho_+(b)$, and thus $b$ prefers $\rho'_+(b)$ to $\rho_+(b)$. Indeed, since $\rho$ has been eliminated in a sequence of rotation eliminations that, starting from $M_0$, lead to $M$, by Lemma~\ref{lem:prefer-all-or-none}, for any stable matching $M'$ that $\rho$ has been just eliminated in, $b$ prefers $M'(b)$ to $\rho_+(b)$,~i.e.~for any $a'\in M'(b)$ and $\hat{a}\in \rho_+(b)$, it satisfies $a'\succ_b \hat{a}$. Notice that for any $a\in M(b)$ and $a'\in M'(b)$, since $M'$ dominates $M$, we have $a\succeq_b a'$. Therefore, for any $a\in M(b)$ and $\hat{a}\in \rho_+(b)$, it holds $a\succ_b \hat{a}$, which completes the claim. On the other hand, if $\rho\triangleright\rho'$, then by symmetry we can deduce that $b$ prefers $\rho_+(b)$ to $\rho'_+(b)$, which shows the converse direction.
\hfill\Halmos

\begin{lemma}[All-or-Nothing Property]\label{lem:rotation-all-or-nothing}
    For two rotations $\rho,\rho'\in\cR(I)$, if $\rho\neq\rho'$, then $\rho_+\cap\rho'_+=\emptyset$ and $\rho_-\cap\rho'_-=\emptyset$.
\end{lemma}
\emph{Proof.} 
%\afcomment{this is now checked}
Let $\rho\neq \rho'\in\cR(I)$. Consider the student-pessimal matching $M_z=M_0/\cR(I)$. By Lemma~\ref{lem:matching-transversal} and Corollary~\ref{cor:all-rotations}, we can find a sequence of matchings $M_0,M_1,\dots,M_k$ in $\cS(I)$ such that $M_k=M_z$ and $\cR(I)=\{\rho(M_{i-1},M_i)|i\in [k]\}$. Therefore, there exists $i,i'\in[k]$ such that $\rho=\rho(M_{i-1},M_i)$ and $\rho'=\rho(M_{i'-1},M_{i'})$. Without loss of generality we assume $i<i'$. If there exists $ab\in\rho_+\cap\rho'_+$, then by definition of $ab\in\rho_+$ we have $ab\in M_{i-1}$ and $ab\notin M_{i}$. Notice that $M_0(a)\succeq_a M_1(a)\succeq_a\dots\succeq_a M_k(a)$, as a result $ab\notin M_{\bar{\imath}}$ for every $\bar{\imath}\ge i$, which contradicts with the fact that $ab\in M_{i'-1}$ as $ab\in\rho'_+$. Thus, $\rho_+\cap\rho'_+=\emptyset$. Similarly, one can show that $\rho_-\cap\rho'_-=\emptyset$.\hfill\Halmos

Next, we show the following lemma that we frequently use in our proofs.

\begin{lemma}\label{lem:rotation-monotone}
    Let $M\in\cS(I)$ and $ab\in M$. Let $\rho \in \cR(I)$. If $ab\in\rho_-$ then $\rho\in R_M$. If $ab\in\rho_+$, then $\rho\notin R_M$.
\end{lemma}
\emph{Proof.} %\afcomment{this is now checked} 
By Lemma~\ref{lem:matching-transversal}, we can find a sequence of stable matchings $M_0,M_1,\dots,M_k=M$ such that $R_M=\{\rho(M_{i-1},M_{i}):i\in[k]\}$. Along this sequence, we have $M_0(a)\succeq_{a}M_1(a)\succeq_{a}\dots\succeq_{a}M_k(a)=b$. If $ab\in\rho_-$, then $\rho_+(a)\succ_{a}b$ and thus $M_0(a)\succ_{a} b$. Let $j\in[k]$ be the smallest index such that $M_j(a)=b$. Then, $ab\in(\rho(M_{j-1},M_{j}))_-$. We can then conclude that $\rho=\rho(M_{j-1},M_{j})$, otherwise $\rho_-\cap(\rho(M_{j-1},M_{j}))_-\neq\emptyset$, 
a contradiction to Lemma~\ref{lem:rotation-all-or-nothing}. Therefore, $\rho\in R_M$. Now, if $ab\in\rho_+$, then there is no $i\in[k]$ such that $\rho=\rho(M_{i-1},M_i)$, otherwise $M_i(a)\prec_{a} b$, a contradiction. Therefore, $\rho\notin R_M$.\hfill\Halmos

\subsection{Matching Siblings to the Same School (MSSS): Omitted Proofs}\label{ec:proofs-siblings-1}

\subsubsection{Proof that Assumption~\ref{ass:5cases} is without loss of generality.}\label{ec:reduction-5cases}

Fix $(a,\bar{a})\in \cC$. We construct an instance $\widetilde \cI$ of (MSSS) over a new sibling instance $(\widetilde I,\cC)$ such that: \begin{itemize}
\item[1)] $\widetilde M_0(a)\neq \widetilde M_0(\bar{a})$ where $\widetilde M_0$ denotes the student-optimal matching of $\widetilde I$ and that for every other couple of siblings $(a', \bar a') \in \cC \setminus \{(a, \bar a)\}$, if $M_0(a') \neq M_0(\bar{a'})$ then $\widetilde M_0(a')\neq \widetilde M_0(\bar{a}')$.

\item[1')] $\widetilde M_z(a)\neq \widetilde M_z(\bar{a})$ where $\widetilde M_z$ denotes the student-pessimal matching of $\widetilde I$ and that for every other couple of siblings $(a', \bar a') \in \cC \setminus \{(a, \bar a)\}$, if $M_z(a') \neq M_z(\bar a')$ then $\widetilde M_z(a')\neq \widetilde M_z(\bar a')$.

\item[2)] If $M\in\cS(I)$ matches a set $\bar {\cal C}\subseteq {\cal C}$ of siblings to  same school, then we can efficiently find $\widetilde  M \in \cS(\widetilde  I)$ that matches $\widetilde{\cal C}$ to the same school with $\bar {\cal C} \subseteq \widetilde{\cal C} \subseteq {\cal C}$, and vice versa. 
\end{itemize}
First, we construct an instance $\hat \cI$ over a new siblings instance $(\hat I, \cC)$ as follows: if $M_0(a) \neq M_0(\bar{a})$, let $\hat \cI = \cI$. Otherwise, let $M_0(a)=M_0(\bar{a})=b$. Let $a_{worst}$ denote the least preferred student of $b$ in $M_0(b)$. To construct $\hat \cI$ we add a dummy student $d$ and a dummy school $s$ with quota $1$. That is, $\hat  I=(A\cup\{d\},B\cup\{s\},\hat 
\succ,\hat q)$ with $\hat  q_b=q_b$ for $b\in B$, and $\hat 
 q_s=1$. Starting from the old preference lists $\succ$, the new preference lists $\hat \succ$ are defined as follows: 
First, we modify the preference lists of $a$ and $b$ such that
$$\cdots \hat \succ_a \;s\; \hat \succ_a \;b \; \hat \succ_a\cdots \qquad  \hbox{and} \qquad \cdots \hat \succ_b \; a\; \hat 
\succ_b \cdots \hat \succ_b \; a_{worst} \; \hat \succ_b \; d \; \hat \succ_b \cdots$$
respectively, where we insert $s$ right before $b$ into $a$'s preference list and $d$ right after $a_{worst}$ (it may happen that $a=a_{worst}$) into $b$'s preference list. 
Second, for every student $a'\in A\setminus\{a\}$ we add $s$ at the end of its list (hence after $\emptyset$), and for every school $b'\in B\setminus\{b\}$ we add $d$ at the end of its list (hence after $\emptyset$). 
Finally, we define the preference lists of $d$ and $s$ as
$$b \hat \succ_d \;s \; \hat  \succ_d \; \emptyset \; \hat \succ_d \cdots ,\qquad \hbox{and} \qquad d \; \hat 
 \succ_s 
 \; a\; \hat \succ_s \; 
\emptyset \; \hat 
\succ_s \cdots,$$
where the dots can be filled with any arbitrary ordering of the remaining agents.

We begin by showing the following claim.

\begin{claim}\label{cl:no-new-blocking}
In the new instance $\hat \cI$,
\begin{itemize}
    \item[(i)] If $M \in \cS(I)$ then $M \cup \{ds\} \in \cS(\hat I)$.
    \item[(ii)] If $\hat M \in \cS(\hat I)$ such that $ds \in \hat M$ then $\hat M \setminus \{ds\} \in \cS(I)$.
    \item[(iii)] If $\hat M \in \cS(\hat I)$ such that $ds \notin \hat M$ then $\{as,db\} \subset \hat M$ and $\hat M \setminus \{as,db\} \cup \{ab\} \in \cS(I)$.
    \item[(iv)] $\hat M_0=M_0\setminus\{ab\}\cup\{as,db\}$ is the student-optimal stable matching in $\hat I$.
    \item[(v)] $\hat M_z=M_z \cup \{ds\}$ is the student-pessimal stable matching in $\hat I$.
\end{itemize}
\end{claim}
{\noindent \em Proof.}
(i) Observe that any pair blocking $M \cup \{ds\}$ in $\hat I$ either blocks $M$ in $I$ or contains at least one of $d,s$. However, $d$ is $s$ favorite partner, so $s$ cannot be part of a pair blocking $M\cup \{ds\}$ in $\hat I$. Moreover, the only agent that $d$ prefers to $s$ is $b$, for which $a'\; \hat \succeq_b \; a_{worst}\; \hat \succ_b \; d$ for every $a' \in M(b)$. Thus, $d$ cannot be part of a pair blocking $M\cup \{ds\}$ in $\hat I$ either, concluding the proof. 

(ii) Any pair blocking $\hat M \setminus \{ds\}$ in $I$ would also block $\hat M$ in $I'$, thus such a pair cannot exist by hypothesis.

(iii) 
By (i), $M_0\cup\{ds\}$ is a stable matching in $ \hat I$. Thus, by rural hospital theorem (see e.g., \cite[Theorem 1.6.3]{gusfieldbook}), $d$ and $s$ are always matched to a partner (different from $\emptyset$) in any stable matching $\hat M \in \cS(\hat I)$. If $ds \notin \hat M$, then $d$ and $s$ can only be matched to $b$ and $a$ respectively, as they are the only remaining agents in their preference lists preferred to $\emptyset$. Hence $\{as,db\} \subset \hat M$. 

Consider now $M=\hat M\setminus\{as,db\}\cup\{ab\}$ and suppose that a pair $a'b'$ blocks $M$, we show that $a'b'$ blocks $\hat M$ under $\hat \succ$. Indeed, if $a'\neq a$ and $b'\neq b$, then since $\hat \succ$ agrees with $\succ$ on $a'$ and $b'$, we have that $a'b'$ blocks $\hat M$. If $a'=a$ and $b'\neq b$, then $ab'$ blocks $M$ implies that $b'\succ_a b$ and $M(b')$ contains a student who is worse than $a$ under $\succ_{b'}$. Notice that $b' \succ_a b$ implies $b'\; \hat \succ_a \; s= \hat M(a)$, and $\hat M(b')=M(b')$, we have that $ab'$ blocks $\hat M$. If $a'\neq a$ and $b'=b$, then we observe that from $M$ to $\hat M$ $a'$ is assigned to the same school and $b$ is getting a worse partner. Therefore, the pair $a'b$ blocking $M$ also blocks $\hat M$ (as $a'$ and $b$ would still prefer each other to their partners). Finally, if $a'=a$ and $b'=b$ is not possible as $ab\in M$.

(iv) First we show that $\hat M_0$ is stable in $\hat I$. Suppose there is a pair $a'b'$ that blocks $\hat M_0$, then it has to be $a'\in\{a,d\}$ or $b'\in\{b,s\}$ otherwise $a'b'$ also blocks $M_0$. $a'=d$ is impossible as $d$ is matched with the best possible school, which also implies $b'=s$ is impossible. The only remaining possibility is $a'b' =ab$ which is also impossible as $a$ is matched to $s$ which they prefer to $b$.

Next, we show that $\hat M_0$ is the student-optimal stable matching in $\hat I$. Suppose to the contrary that there exists $\hat M \in {\cal S}(\hat I)$ such that $\hat M \; \hat >\; \hat M_0$, where $\hat >$ denotes the partial order of stable matchings in $\hat I$. %Thus, there exist at least two students $\tilde a$ such that $M'(\tilde a)\succ_{\tilde a} \hat M_0(\tilde a)$.
Assume first $ds \in \hat M$. Then, by (ii),  $M=\hat M\setminus \{ds\}\in {\cal S}(I)$. At least two students prefer $\hat M$ to $\hat M_0$ in $\hat I$. Let $a' \neq d$ be one such student. If $a' \neq a$, then $M(a') \succ_{a'} \hat M_0(a') = M_0(a')$, contradicting the student-optimality of $M_0$ in $I$. Thus, $a' = a$. But then $M(a)\succ_a \hat M_0(a)= s \succ_a b= M_0(a)$, again a contradiction. Thus, $ds \notin \hat M$. Then by (iii), $M=\hat M \setminus \{as,db\} \cup \{ab\} \in \cS(I)$. Since $\hat M \;\hat >\; \hat M_0$, there exist at least two students that prefer $\hat M$ to $\hat M_0$. Because $\hat M(d)=\hat M_0(d)$, these students are different from $d$. Hence, there exists $a' \neq a, d$ such that $M(a')=\hat M(a')\succ_{a'} \hat M_0(a')=M_0(a')$, again contradicting the student-optimality of $M_0$ in $I$. 

(v) Let $\hat M$ denote the student-pessimal stable matching of $\hat I$. Then, $\hat M(s) = d$ (since $ds$ is stable by (i) and $d$ is $s$'s most preferred partner). Hence, by (ii), $\hat M \setminus \{ds\} \in \cS(I)$, and hence $\hat M \setminus \{ds\} \geq M_z$. But then every student prefers $\hat M$ to $\hat M_z = M_z \cup \{ds\}$ in $\hat I$. By (i) $\hat M_z \in \cS(\hat I)$, hence $\hat M \; \hat \geq \;\hat M_z$ implying that $\hat M = \hat M_z$.
\hfill$\square$

\smallskip

Next, we construct an instance $\widetilde \cI$ over a new siblings instance $(\widetilde I, \cC)$ as follows: If $\hat M_z(a) \neq \hat M_z(\bar a)$, let $\widetilde \cI = \hat \cI$. Otherwise, let $\hat M_z(a) = \hat M_z(\bar a) = \bar{b}$. 
To construct $\widetilde \cI$, we add a dummy student $\bar d$ and a dummy school $\bar s$ with quota $1$. Starting from the old preference lists $\hat \succ$, the new preference lists $\widetilde \succ$ are defined as follows: 
$$\cdots\widetilde{\succ}_{\bar{a}} \; \bar{b} \; \widetilde{\succ}_{\bar{a}} \;\bar{s}\; \widetilde{\succ}_{\bar{a}}\cdots \quad  \hbox{and,} \quad \cdots\widetilde{\succ}_{\bar{b}} \;\bar{d}\; \widetilde{\succ}_{\bar{b}} \;\bar{a}\; \widetilde{\succ}_{\bar{b}} \cdots \quad \hbox{and,} \quad \;\bar{s}\;\widetilde{\succ}_{\bar{d}} \;\bar{b}\; \widetilde{\succ}_{\bar{d}} \;\emptyset\;\widetilde{\succ}_{\bar{d}} \cdots \quad \hbox{and,} \quad \;\bar{a}\;\widetilde{\succ}_{\bar{s}} \;\bar{d}\;\widetilde{\succ}_{\bar{s}} \;\emptyset\;\widetilde{\succ}_{\bar{s}} \cdots,$$
and finally for every other student $a' \neq \bar a, \bar d$ we add $\bar s$ at the end of its list (hence after $\emptyset$), and for every other school $b' \neq \bar b, \bar s$ we add $\bar d$ at the end of its list (hence after $\emptyset$). As before, we claim the following.

\begin{claim}\label{cl:no-new-blocking2}
In the new instance $\widetilde \cI$,
\begin{itemize}
    \item[(i)] If $\hat M \in \cS(\hat I)$ then $\hat M \cup \{\bar{d}\bar{s}\} \in \cS(\widetilde I)$.
    \item[(ii)] If $\widetilde M \in \cS(\widetilde I)$ such that $\bar{d}\bar{s} \in \widetilde M$ then $\widetilde M \setminus \{\bar d \bar s\} \in \cS(\hat I)$.
    \item[(iii)] If $\widetilde M \in \cS(\widetilde I)$ such that $\bar{d}\bar{s} \notin \widetilde M$ then $\{\bar a \bar s, \bar d \bar b\} \subset \widetilde M$ and $\widetilde M \setminus \{\bar a \bar s, \bar d \bar b\} \cup \{\bar a \bar b\} \in \cS(\hat I)$.
    \item[(iv)] $\widetilde M_z= \hat M_z\setminus\{\bar a \bar b\}\cup\{\bar a \bar s, \bar d \bar b\}$ is the student-pessimal stable matching in $\widetilde I$.
    \item[(v)] $\widetilde M_0=\hat M_0 \cup \{\bar d \bar s\}$ is the student-optimal stable matching in $\widetilde I$.
\end{itemize}
\end{claim}
{\noindent \em Proof.}
(i) Observe that any pair blocking $\hat M \cup \{\bar{d}\bar{s}\}$ in $\widetilde I$ either blocks $\hat M$ in $\hat I$ or contains at least one of $\bar d, \bar s$. However, $\bar d$'s favorite partner is $\bar s$, so $\bar d$ cannot be part of such blocking pair. Moreover, the only agent that $\bar s$ prefers to $\bar d$ is $\bar a$, for which $\hat M(\bar a) \; \hat \succeq_{\bar a} \; \hat M_z(\bar a) = \bar b\; \hat \succ_{\bar a} \; \bar s$. Thus, $\bar s$ cannot be part of such a blocking pair either, concluding the proof. 

(ii) Any pair blocking $\widetilde M \setminus \{\bar d \bar s\}$ in $\hat I$ would also block $\widetilde M$ in $\widetilde I$, thus such a pair cannot exist.

(iii) 
By (i), $\hat M_z\cup\{\bar d \bar s\}$ is a stable matching in $ \hat I$. Thus, by rural hospital theorem (see e.g., \cite[Theorem 1.6.3]{gusfieldbook}), $\bar d$ and $\bar s$ are always matched to a partner (different from $\emptyset$) in any stable matching $\widetilde M \in \cS(\widetilde I)$. If $\bar d \bar s \notin \widetilde M$, then $\bar d$ and $\bar s$ can only be matched to $\bar b$ and $\bar a$ respectively, as they are the only remaining agents in their preference lists preferred to $\emptyset$. Hence $\{\bar a \bar s, \bar d \bar b\} \subset \widetilde M$. 

Consider now $\hat M=\widetilde M \setminus \{\bar a \bar s, \bar d \bar b\} \cup \{\bar a \bar b\}$ and suppose that a pair $a'b'$ blocks $\hat M$, we show that $a'b'$ blocks $\widetilde M$ under $\widetilde \succ$. Indeed, if $a'\neq \bar a$ and $b'\neq \bar b$, then since $\widetilde \succ$ agrees with $\hat \succ$ on $a'$ and $b'$, we have that $a'b'$ blocks $\widetilde M$. If $a'=\bar a$ and $b'\neq \bar b$, then $\bar a b'$ blocks $\hat M$ implies that $b'\; \hat \succ_{\bar a} \;\bar b$ and $M(b')$ contains a student who is worse than $\bar a$ under $\hat \succ_{b'}$. Notice that $b'\; \hat \succ_{\bar a} \;\bar b$ implies $b'\; \hat \succ_{\bar a} \; \bar s= \widetilde M(\bar a)$, and $\widetilde M(b')= \hat M(b')$, hence $\bar a b'$ blocks $\widetilde M$. If $a'\neq \bar a$ and $b'= \bar b$, then observe that $a'$ is assigned to the same school in $\hat M$ and $\widetilde M$ hence $\bar b = b' \;\widetilde \succ_{a'} \; \hat M(a') = \widetilde M(a')$. Next, we have $a' \;\widetilde \succ_{\bar b} \; a''$ for some $a'' \in \hat M(\bar b)$. If $a'' \neq \bar a$ then $a'' \in \widetilde M(\bar b)$ and hence $a'\bar b$ is a blocking pair of $\widetilde M$. If $a'' = \bar a$ then because $a' \neq \bar d$ (as $d$ is not a student in $\hat I$) and $\bar d$ is the immediate student that $\bar b$ prefers to $\bar a$, it must be the case that $a' \;\widetilde \succ_{\bar b} \; \bar d$, which again implies that $a' \bar b$ is a blocking pair of $\widetilde M$. Finally, if $a'=\bar a$ and $b'=\bar b$ is not possible as $\bar a \bar b\in \hat M$.

(iv) First we show that $\widetilde M_z$ is stable in $\widetilde I$. Suppose there is a pair $a'b'$ that blocks $\widetilde M_z$, then it has to be $a'\in\{\bar a, \bar d\}$ or $b' \in\{\bar b, \bar s\}$ otherwise $a'b'$ also blocks $\hat M_z$. $b'=\bar s$ is impossible as $\bar s$ is matched with the best possible student, which also implies $a'=\bar d$ is impossible. The only remaining possibility is $a'b' =\bar a \bar b$, in this case, there exists $a'' \in \widetilde M_z(\bar b)$ such that $\bar a \; \widetilde \succ_{\bar b} \; a''$ but then since $a'' \neq \bar d$ (since $\bar a \; \widetilde \prec_{\bar b} \; \bar d$) this implies that $a'' \in \hat M_z(\bar b)$ such that $\bar a \; \hat \succ_{\bar b} \; a''$ and that $a'b'$ is blocking for $\hat M_z$ as well, a contradiction.

Next, we show that $\widetilde M_z$ is the student-pessimal stable matching in $\widetilde I$. Suppose to the contrary that there exists $\widetilde M \in {\cal S}(\widetilde I)$ such that $\widetilde M \; \widetilde <\; \widetilde M_z$, where $\widetilde >$ denotes the partial order of stable matchings in $\widetilde I$.
Assume first $\bar d \bar s \in \widetilde M$. Then, by (ii),  $\hat M=\widetilde M\setminus \{\bar d \bar s\}\in {\cal S}(\hat I)$. At least two students prefer $\widetilde M_z$ to $\widetilde M$ in $\widetilde I$. Let $a' \neq \bar d$ be one such student. If $a' \neq \bar a$, then $\hat M(a') = \widetilde M(a')\; \hat \prec_{a'} \; \widetilde M_z(a') = \hat M_z(a')$, contradicting the student-pessimality of $\hat M_z$ in $\hat I$. Thus, $a' = \bar a$. But then $\hat M(\bar a) = \widetilde M(\bar a) \; \hat \prec_{\bar a}\; \widetilde M_z(\bar a)= \bar{s} \; \hat \prec_{\bar a} \;\bar b= \hat M_z(\bar a)$, again a contradiction. Thus, $\bar d \bar s \notin \widetilde M$. Then by (iii), $\hat M=\widetilde M \setminus \{\bar a \bar s, \bar d \bar b\} \cup \{\bar a \bar b\} \in \cS(\hat I)$. Since $\widetilde M \;\widetilde <\; \widetilde M_z$, there exists at least two students that prefer $\widetilde M_z$ to $\widetilde M$. Because $\widetilde M(\bar d) \neq \bar s$ we have $\widetilde M(\bar d) = \bar b = \widetilde M_z(\bar d)$, these students are therefore different from $\bar d$. Hence, there exists $a' \neq \bar a, \bar d$ such that $\hat M(a')=\widetilde M(a') \; \hat \prec_{a'} \; \widetilde M_z(a')=\hat M_z(a')$, again contradicting the student-pessimality of $\hat M_z$ in $\hat I$.

(v) Let $\hat M$ denote the student-optimal stable matching of $\widetilde I$. Then, $\widetilde M(\bar d) = \bar s$ (since $\bar d \bar s$ is stable by (i) and $\bar s$ is $\bar d$'s most preferred partner). Hence, by (ii), $\widetilde M \setminus \{\bar d \bar s\} \in \cS(\hat I)$, and hence $\widetilde M \setminus \{\bar d \bar s\} \; \hat \leq \; \hat M_0$. But then every student prefers $\widetilde M_0 = \hat M_0 \cup \{\bar d \bar s\}$ to $\widetilde M$ in $\widetilde I$. By (i) $\widetilde M_0 \in \cS(\widetilde I)$, hence $\widetilde M \; \widetilde \leq \;\widetilde M_0$ implying that $\widetilde M = \widetilde M_0$.
\hfill$\square$

\smallskip

Note that both of the above constructions (of $\hat \cI$ and of $\widetilde \cI$) are efficient. The resulting instance $\widetilde \cI$ verifies 1), 1') and 2). Indeed, by Claim~\ref{cl:no-new-blocking}, (iv), the first construction of $\hat \cI$ guarantees that $(a, \bar a)$ are not matched to the same school in the student optimal matching of $\hat \cI$ and that given any other pair $(a', \bar a') \in \cC$, if the pair are not matched to the same school in the student optimal matching of $\cI$ then they are not in that of $\hat \cI$ either. Then given Claim~\ref{cl:no-new-blocking2}, (v), the construction of $\widetilde \cI$ does not change the pairs of siblings assigned to the same school in the student-optimal stable matching. Hence 1) holds. 1') holds by a symmetric argument. To see that 2) holds, consider $M \in \cS(I)$ which matches a set $\bar \cC$ of siblings to the same school. Then, by Claim~\ref{cl:no-new-blocking} (i) and Claim~\ref{cl:no-new-blocking2} (i), the matching $M \cup \{ds, \bar d \bar s\} \in \cS(\widetilde I)$ and (by construction) matches exactly the set of siblings $\bar \cC$ to the same school. For the converse, let $\widetilde M \in \cS(\widetilde I)$, Claim~\ref{cl:no-new-blocking2} (ii) and (iii), give a construction of a stable matching $\hat M \in \cS(\hat I)$ that matches at least the set of siblings $\bar \cC$ to the same school. Then, starting from such matching, Claim~\ref{cl:no-new-blocking} (ii) and (iii) give a construction of a stable matching $M \in \cS(I)$ that matches at least the set of siblings $\bar \cC$ to the same school. This implies 2).

We can now iteratively repeat the construction for every pair of siblings so as to construct, in polynomial time, and instance ${\cal I}'$. Notice that by repeatedly applying property 1) (resp., 1')), once a pair of siblings $(a,\bar{a})$ is not matched to the same school in the student (resp. school)-optimal stable matching, they are not matched to the same school in the student (resp., school)-optimal matching in later rounds. Thus, no two pair of siblings are matched to the same school in either the student- or the student-pessimal stable matchings of ${\cal I}'$. Given 2), the maximum number of pairs of siblings matched to the same school in some stable matching of $\cI'$ coincides with its maximum over the stable matchings of $\cI$. Indeed, if there exists $k$ pairs of siblings matched to the same school in $M \in {\cal S}(I)$, then by repeated application of 2), there exists $M' \in {\cal S}(I')$ that also matches at least $k$ pairs of siblings to the same school, and vice versa.

\subsubsection{Proof of Lemma~\ref{lem:2cases}.}

Fix a pair of siblings $(a, \bar{a}) \in \cC$ and a school $b \in C$, we show that if there exists a stable matching $M\in\cS(I)$ with $M(a)=M(\bar{a})=b$, then case 2 of the lemma happens.  Assume there exists $M\in\cS(I)$ such that $M(a)=M(\bar{a})=b$.

\underline{Defining $\rho_{in}$:} 
Notice that by Assumption~\ref{ass:5cases}, we have $M_0(a)\neq M_0(\bar{a})$. We distinguish the following cases: (i) $M_0(a)\neq b$, $M_0(\bar{a})= b$, in which case there exists a rotation $\rho_1\in R_M$ such that $ab\in\rho_{1-}$; (ii) $M_0(a)=b$, $M_0(\bar{a})\neq b$, in which case there exists a rotation $\rho_2\in R_M$ such that $\bar{a}b\in\rho_{2-}$; (iii) $M_0(a)\neq b$, $M_0(\bar{a})\neq b$ in which case there exists two rotations $\rho_3, \bar{\rho}_3\in R_M$ ($\rho_3$ and $\bar{\rho}_3$ may be identical) such that $ab\in\rho_{3-}$ and $\bar{a}b\in\bar{\rho}_{3-}$. %(iv) $M_0(a)=M_0(\bar{a})=b$. 
If (i) happens, we let $\rho_{in}=\rho_1$. If (ii) happens, we let $\rho_{in}=\rho_2$. If (iii) happens, then since $b\in B(\rho_3)\cap B(\bar{\rho}_3)$, by Lemma~\ref{lem:dominace-onemany}, $\rho_3$ and $\bar{\rho}_3$ must be comparable (including the case where $\rho_3=\bar{\rho}_3$). We let $\rho_{in}$ be the \emph{smaller} one (w.r.t. $\trianglerighteq$) between $\rho_3$ and $\bar{\rho}_3$. 

\underline{Defining $\rho_{out}$:} Similarly, to define $\rho_{out}$, we consider three cases: (i) $M_z(a)\neq b$, $M_z(\bar{a})=b$, and there is a rotation $\rho_1\in \cR(I) \setminus R_M$ such that $ab\in{\rho}_{1+}$; (ii) $M_z(a)=b$, $M_z(\bar{a})\neq b$, and there is a rotation $\rho_2\in \cR(I) \setminus R_M$ such that $\bar{a}b\in{\rho}_{2+}$; (iii) $M_z(a)\neq b$, $M_z(\bar{a})\neq b$, and there are $\rho_3, \bar{\rho}_3\in \cR(I) \setminus R_M$ ($\rho_3$ and $\bar{\rho}_3$ may be identical) such that $ab\in{\rho}_{3+}$ and $\bar{a}b\in\bar{\rho}_{3+}$. %(iv) $M_z(a)=M_z(\bar{a})=b$. 
When (i) or (ii) happens, $\rho_{out}$ is defined as the rotation $\rho_1$ or $\rho_2$, respectively. If (iii) happens, then since $b\in B(\rho_3)\cap B(\bar{\rho}_3)$, by Lemma~\ref{lem:dominace-onemany}, we have that $\rho_3$ and $\bar{\rho}_3$ are comparable (including the case $\rho_3=\bar{\rho}_3$). Then, we let $\rho_{out}$ be the \emph{larger} one (w.r.t. $\trianglerighteq$) between $\rho_3$ and $\bar{\rho}_3$. 

Let $M'\in\cS(I)\setminus \{M_0,M_z\}$ be a stable matching such that $M'(a)=M'(\bar{a})=b$. Then, since $M_0$ does not match both $a$ and $\bar{a}$ to $b$, by definition of $\rho_{in}$, at least one of $ab\in(\rho_{in})_-$ and $\bar{a}b\in(\rho_{in})_-$ holds. By Lemma~\ref{lem:rotation-monotone}, $\rho_{in}\in R_{M'}$. Similarly, $M_z(a)\neq b$ or $M_z(\bar{a})\neq b$ implies $ab\in(\rho_{out})_+$ or $\bar{a}b\in(\rho_{out})_+$. Again by Lemma~\ref{lem:rotation-monotone}, we have $\rho_{out}\notin R_{M'}$.

Now for the converse, consider a stable matching $M'\in\cS(I)$, such that $\rho_{in}\in R_{M'}$ and $\rho_{out}\notin R_{M'}$. If $M'(a)=M'(\bar{a})=b$ does not hold, let us first analyze the case $M'(a)\succ_a b$. Recall the definition of $\rho_{in}$. If we are in case (i) or (iii), then $ab\in(\rho_{in})_-$, which implies $M'(a)\preceq_a b$ since $\rho_{in}\in R_{M'}$, a contradiction. Case (ii) is impossible since $M_0(a)=b\prec_a M'(a)$ contradicts with the student-optimality of $M_0$. Therefore, $M'(a)\succ_a b$ never holds true. By symmetry we can rule out $M'(\bar{a})\succ_{\bar{a}} b$. By applying a similar argument to the definition of $\rho_{out}$, we can also eliminate the possibility of $M'(a)\prec_a b$ or $M'(\bar{a})\prec_{\bar{a}} b$. Thus, we must have $M'(a)=M'(\bar{a})=b$. 

Now for the unicity of the pair $\rho_{in},\rho_{out}$, if there exists another pair $\rho,\rho'$ such that for every $M' \in \cS(I)$, $M'(a)=M'(\bar{a})=b$ if and only if $\rho\in R_{M'}$ and $\rho'\notin R_{M'}$. Then, it cannot be the case that $\rho'\triangleright \rho$, otherwise, given $M'$ such that $M'(a)=M'(\bar{a})=b$, $\rho\in R_{M'}$ implies $\rho'\in R_{M'}$. Consider the stable matching $M^\rho$, then $\rho'\notin R_{M^\rho}$. Therefore $M^{\rho}(a)=M^{\rho}(\bar{a})=b$, which implies $\rho_{in}\in R_{M^\rho}$ and thus $\rho_{in}\triangleright \rho$. Similarly, by applying the argument for $M^{\rho_{in}}$, we can obtain $\rho\triangleright \rho_{in}$. Hence, $\rho=\rho_{in}$. A symmetric argument can be shown to prove $\rho'=\rho_{out}$. Thus the uniqueness is guaranteed. % Hence, we are in case 5.

Finally, for the fact that $\rho_{in}(a,\bar{a},b) \triangleright \rho_{out}(a,\bar{a},b)$, by definition of $\rho_{in}(a,\bar{a},b)$ and $\rho_{out}(a,\bar{a},b)$, $b \in B(\rho_{in}(a,\bar{a},b)) \cap B(\rho_{out}(a,\bar{a},b))$, hence, by Lemma~\ref{lem:dominace-onemany}, $\rho_{in}(a,\bar{a},b)$ and $\rho_{out}(a,\bar{a},b)$ are comparable. Let $M'$ such that $M'(a)=M'(a')=b$. If $\rho_{in}(a,\bar{a},b) \triangleleft \rho_{out}(a,\bar{a},b)$ then because $\rho_{in}(a,\bar{a},b) \in R_{M'}$ we would have that $\rho_{out}(a,\bar{a},b) \in R_{M'}$ which contradicts the fact that $\rho_{out}(a,\bar{a},b) \notin R_{M'}$. Hence, $\rho_{in}(a,\bar{a},b) \triangleright \rho_{out}(a,\bar{a},b)$.
\hfill\Halmos

\subsubsection{Proof of Lemma~\ref{lem:differentials-sibling}}

We have,

\noindent {\bf First order differentials:}
By definition, $\partial f^i_\rho=f^i(M^\rho)-f^i(M_\rho)$. Notice that $\partial f^i_\rho\in\{0,\pm 1\}$. In the following, we show that $\rho\in\cR_{out}(a_i,\bar{a}_i)$ and $\rho\notin\cR_{in}(a_i,\bar{a}_i)$ if and only if $\partial f_\rho^i = 1$ and $\rho\notin\cR_{out}(a_i,\bar{a}_i)$ and $\rho\in\cR_{in}(a_i,\bar{a}_i)$ if and only if $\partial f_\rho^i = -1$. The lemma follows.

\underline{$\rho\in\cR_{out}(a_i,\bar{a}_i)$ and $\rho\notin\cR_{in}(a_i,\bar{a}_i)$ if and only if $\partial f_\rho^i = 1$: } Let $b$ such that $\rho=\rho_{out}(a_i,\bar{a}_i,b)$. 

By Lemma~\ref{lem:2cases}, $\rho_{in}(a_i,\bar{a}_i,b)$ exists and is such that $\rho_{in}(a_i,\bar{a}_i,b)\triangleright \rho_{out}(a_i,\bar{a}_i,b)$. Consider the matching $M_\rho=M_0/R_\rho$. By definition, we have $\rho_{in}(a_i,\bar{a}_i,b)\in R_\rho$ and $\rho_{out}(a_i,\bar{a}_i,b)\notin R_\rho$, thus by Lemma~\ref{lem:2cases} again, we have $M_\rho(a_i)=M_\rho(\bar{a}_i)=b$, which implies $f^i(M_\rho)=0$. 

Now, since $R^\rho=R_\rho\cup\{\rho_{out}(a_i,\bar{a}_i,b)\}$, by Lemma~\ref{lem:2cases} at least one of $M^\rho(a_i)= b$ or $M^\rho(\bar{a}_i)=b$ does not hold. Without loss of generality we assume that $M^\rho(a_i)=b'\neq b$, we argue that $M^\rho(\bar{a}_i)\neq b'$. Indeed, if $M^\rho(\bar{a}_i)=b'$ also holds, then we apply Lemma~\ref{lem:2cases} to $M^\rho$, we deduce that $\rho_{in}(a_i,\bar{a}_i,b')\in R^\rho$ and $\rho_{out}(a_i,\bar{a}_i,b')\notin R^\rho$. Again, by applying Lemma~\ref{lem:2cases} to $M_\rho$, we have either $\rho_{in}(a_i,\bar{a}_i,b')\notin R_\rho$ or $\rho_{out}(a_i,\bar{a}_i,b')\in R_\rho$. The latter is impossible, since we have deduced that $\rho_{out}(a_i,\bar{a}_i,b')\notin R^\rho$ and $R_\rho\subset R^\rho$, thus $\rho_{in}(a_i,\bar{a}_i,b')\notin R_\rho$. Now, since $R^\rho=R_\rho\cup\{\rho\}$, we have $\rho=\rho_{in}(a_i,\bar{a}_i,b')$, which contradicts with $\rho\notin\cR_{in}(a_i,\bar{a}_i)$. Thus, $M^\rho(\bar{a}_i)=b'$ is false. We then have $M^\rho(a_i)\neq M^\rho(\bar{a}_i)$, which gives $f^i(M^\rho)=1$. Hence, $\partial f^i_\rho=f^i(M^\rho)-f^i(M_\rho)=1$. 

Conversely, suppose that $\partial f^i_\rho=1$. There exists $b\in B$ such that $M^\rho(a_i)=M^\rho(\bar{a}_i)=b$ and $M_\rho(a_i)\neq M_\rho(\bar{a}_i)$. By Lemma~\ref{lem:2cases}, we have $\rho_{in}(a_i,\bar{a}_i,b)\in R^\rho$ and $\rho_{out}(a_i,\bar{a}_i,b)\notin R^\rho$. By definition, we have $R^\rho=R_\rho\cup\{\rho\}$. If $\rho\neq \rho_{in}(a_i,\bar{a}_i,b)$, then $\rho_{in}(a_i,\bar{a}_i,b)\in R_\rho$ and $\rho_{out}(a_i,\bar{a}_i,b)\notin R_\rho$. 
Thus by Lemma~\ref{lem:2cases}, $M_\rho(a_i)=M_\rho(\bar{a}_i)=b$, a contradiction. Thus $\rho=\rho_{in}(a_i,\bar{a}_i,b)$, which implies $\rho\in \cR_{in}(a_i,\bar{a}_i)$. If in addition $\rho\in \cR_{out}(a_i,\bar{a}_i)$, suppose that $\rho=\rho_{out}(a_i,\bar{a}_i,b')$. Then we have $\rho_{in}(a_i,\bar{a}_i,b')\in R_\rho$ since $\rho_{in}(a_i,\bar{a}_i,b')\triangleright\rho$. Again by Lemma~\ref{lem:2cases}$, \rho_{in}(a_i,\bar{a}_i,b')\in R_\rho$ and $\rho_{out}(a_i,\bar{a}_i,b')\notin R_\rho$ gives $M_\rho(a_i)=M_\rho(\bar{a}_i)=b'$, a contradiction. Therefore, $\rho\notin \cR_{out}(a_i,\bar{a}_i)$, as required. 

\underline{$\rho\in\cR_{out}(a_i,\bar{a}_i)$ and $\rho\notin\cR_{in}(a_i,\bar{a}_i)$ if and only if $\partial f_\rho^i = 1$:}
This case follows by a similar argument.

{\noindent \bf Second order differentials: } For the second order differentials, our goal is to show that for every pair of distinct rotations $\rho \neq \rho'$, it holds that $f^i(M^\rho\wedge M_{\rho'})+f^i(M_\rho\wedge M^{\rho'})=f^i(M_\rho\wedge M_{\rho'})+f^i(M^\rho\wedge M^{\rho'})$.

Assume first that $\rho$ and $\rho'$ are comparable. Without loss of generality, assume $\rho\triangleright\rho'$. Then, for every $\rho''$ such that $\rho''\trianglerighteq\rho$, we have $\rho''\triangleright\rho'$. Therefore, by definition, $M^\rho\wedge M_{\rho'}=M_{\rho'}$, $M_\rho\wedge M^{\rho'}=M^{\rho'}$. Similarly, $M_\rho\wedge M_{\rho'}=M_{\rho'}$ and $M^\rho\wedge M^{\rho'}=M^{\rho'}$. Hence, we obtain $f^i(M^\rho\wedge M_{\rho'})+f^i(M_\rho\wedge M^{\rho'})=f^i(M_{\rho'})+f^i(M^{\rho'})=f^i(M_\rho\wedge M_{\rho'})+f^i(M^\rho\wedge M^{\rho'})$.

Next assume that $\rho,\rho'$ are not comparable. We distinguish the following cases.

\underline{(a). $f^i(M^\rho\wedge M_{\rho'})+f^i(M_\rho\wedge M^{\rho'})=2$:} Then neither $M^\rho\wedge M_{\rho'}$ nor $M_\rho\wedge M^{\rho'}$ match $a$ and $\bar{a}$ to the same school. Suppose $f^i(M_\rho\wedge M_{\rho'})=0$. That is, $M_\rho\wedge M_{\rho'}$ matches both $a$ and $\bar{a}$ to a school $b$. Since $\rho$ and $\rho'$ are not comparable, we have $\rho\notin R_{\rho'}$. 
Therefore, $R_{M^{\rho}\wedge M_{\rho'}} \setminus R_{M_{\rho} \wedge M_{\rho'}} = R^{\rho} \cup R_{\rho'} \setminus R_{\rho} \cup R_{\rho'} = \{\rho\}$ implying that $M^\rho\wedge M_{\rho'}$ is an immediate predecessor of $M_\rho\wedge M_{\rho'}$ such that $(M_\rho\wedge M_{\rho'})/\rho=M^\rho\wedge M_{\rho'}$. By Lemma~\ref{lem:2cases}, $\rho=\rho_{out}(a,\bar{a},b)$, since by eliminating $\rho$, the siblings $a$ and $\bar{a}$ move from the same school $b$ to different schools. In addition, since $\rho_{in}(a,\bar{a},b)\triangleright\rho_{out}(a,\bar{a},b)$, we have $\rho_{in}(a,\bar{a},b)\in R_{\rho}$. Now consider the matching $M_\rho\wedge M^{\rho'}$. Let $R$ be the set of rotations such that $M_\rho\wedge M^{\rho'}=M_0/ R$, then since $\rho\notin R$, we have $\rho_{out}(a,\bar{a},b)\notin R$. Since $R_{\rho}\subseteq R$, we have $\rho_{in}(a,\bar{a},b)\in R$. Therefore, by Lemma~\ref{lem:2cases}, the matching $M_\rho\wedge M^{\rho'}$ assigns both $a$ and $\bar{a}$ to the same school $b$, which implies $f^i(M_\rho\wedge M^{\rho'})=0$, a contradiction to the assumption that $f^i(M^\rho\wedge M_{\rho'})+f^i(M_\rho\wedge M^{\rho'})=2$. Therefore, we obtain $f^i(M_\rho\wedge M_{\rho'})=1$. If $f^i(M^\rho\wedge M^{\rho'})=0$, then $M^\rho\wedge M^{\rho'}$ assigns both $a$ and $\bar{a}$ to the same school, say $b$. By an argument similar to the one above, we can show that $\rho=\rho_{in}(a,\bar{a},b)$, which leads to $f^i(M^\rho\wedge M_{\rho'})=0$, a contradiction.

\underline{(b). $f^i(M^\rho\wedge M_{\rho'})+f^i(M_\rho\wedge M^{\rho'})=1$:} Without loss of generality, let us assume $f^i(M^\rho\wedge M_{\rho'})=0$ and $f^i(M_\rho\wedge M^{\rho'})=1$. The former equation means $M^\rho\wedge M_{\rho'}$ assigns both $a$ and $\bar{a}$ to the same school, say $b$, and the latter equation means the opposite. Let $M^\rho\wedge M_{\rho'}=M_0/ R$, then by Lemma~\ref{lem:2cases}, $\rho_{in}(a,\bar{a},b)\in R$ and $\rho_{out}(a,\bar{a},b)\notin R$. We claim that exactly one of the following is true: (i) $\rho_{in}(a,\bar{a},b)=\rho$; (ii) $\rho_{out}(a,\bar{a},b)=\rho'$. Suppose both (i) and (ii) hold, then by Lemma~\ref{lem:2cases}, we have $\rho\triangleright\rho'$, a contradiction. If neither (i) and (ii) is true, then the matching $M_\rho\wedge M^{\rho'}$ should also match $a$ and $\bar{a}$ to $b$, as eliminating $\rho'$ from and adding $\rho$ to $M^\rho\wedge M_{\rho'}$ do not change the fact that $a$ and $\bar{a}$ are matched to $b$. Thus $f^i(M_\rho\wedge M^{\rho'})=0$, a contradiction. 

Suppose therefore that (i) holds but (ii) does not. Then, eliminating $\rho'$ from the matching $M^\rho\wedge M_{\rho'}$ will still match both $a$ and $\bar{a}$ to $b$, thus $f^i(M^\rho\wedge M^{\rho'})=0$. To complete the argument, we need $f^i(M_\rho\wedge M_{\rho'})=1$. Suppose by contradiction that $f^i(M_\rho\wedge M_{\rho'})=0$. Then $M_\rho\wedge M_{\rho'}$ assigns both $a$ and $\bar{a}$ to the same school, say $b'$. Notice that $b'\neq b$, since (i) implies that $\rho_{in}(a,\bar{a},b)=\rho$ has not been eliminated in $M_\rho\wedge M_{\rho'}$. Therefore, by Lemma~\ref{lem:2cases}, $\rho=\rho_{out}(a,\bar{a},b')=\rho_{in}(a,\bar{a},b)$. Consider the matching $M_\rho\wedge M^{\rho'}$ and set $M_\rho\wedge M^{\rho'}=M_0/ R_1$. Since $\rho\notin R_1$, we have $\rho_{out}(a,\bar{a},b')\notin R_1$. Also, since $M_\rho\wedge M_{\rho'}$ is such that $a$ and $\bar{a}$ are both matched to $b'$, we obtain that $\rho_{in}(a,\bar{a},b')\in R_1\setminus\{\rho'\}\subseteq R_1$. Therefore, by Lemma~\ref{lem:2cases}, $M_\rho\wedge M^{\rho'}$ matches both $a$ and $\bar{a}$ to $b'$, which gives $f^i(M_\rho\wedge M^{\rho'})=0$, a contradiction. Therefore, we have $f^i(M_\rho\wedge M_{\rho'})=1$ and the equality holds. Similarly, the case when (i) fails and (ii) holds leads to $f^i(M_\rho\wedge M_{\rho'})=0$ and $f^i(M^\rho\wedge M^{\rho'})=1$, which also gives the desired result.

\underline{(c). $f^i(M^\rho\wedge M_{\rho'})+f^i(M_\rho\wedge M^{\rho'})=0$:} In this case, both $M^\rho\wedge M_{\rho'}$ and $M_\rho\wedge M^{\rho'}$ assign $a$ and $\bar{a}$ to the same school. Let $M^\rho\wedge M_{\rho'}(a)=M^\rho\wedge M_{\rho'}(\bar{a})=b$ and $M_\rho\wedge M^{\rho'}(a)=M_\rho\wedge M^{\rho'}(\bar{a})=b'$. First, let us assume $b\neq b'$. Notice that $M_\rho\wedge M_{\rho'}=(M^\rho\wedge M_{\rho'})\vee(M_\rho\wedge M^{\rho'})$ and $M^\rho\wedge M^{\rho'}=(M^\rho\wedge M_{\rho'})\wedge(M_\rho\wedge M^{\rho'})$. Thus, $M_\rho\wedge M_{\rho'}$ matches $a$ (resp.~$\bar{a}$) to the better school of $a$'s (resp.~$\bar{a}$'s) list between $b$ and $b'$, while $M^\rho\wedge M^{\rho'}$ matches $a$ (resp.~$\bar{a}$) to the worse school of $a$'s (resp.~$\bar{a}$'s) list between $b$ and $b'$. Assume wlog that $b\succ_a b'$. If we also have $b\succ_{\bar{a}} b'$, then we obtain that $M_\rho\wedge M_{\rho'}$ matches both $a$ and $\bar{a}$ to $b$, and $M^\rho\wedge M^{\rho'}$ matches both $a$ and $\bar{a}$ to $b'$, which leads to $f^i(M_\rho\wedge M_{\rho'})+f^i(M^\rho\wedge M^{\rho'})=0$. Conversely, if we have $b'\succ_{\bar{a}} b$, then by the opposition of interests (Lemma~\ref{lem:prefer-all-or-none}), we deduce that $b$ prefers $M^\rho\wedge M_{\rho'}$ to $M_\rho\wedge M^{\rho'}$. Use the same argument for $b\succ_a b'$, we deduce that $b$ prefers $M_\rho\wedge M^{\rho'}$ to $M^\rho\wedge M_{\rho'}$, which is a contradiction. Therefore, this case ($b\succ_a b'$ and $b'\succ_{\bar{a}} b$) does not happen. Finally, if $b=b'$, since $M_\rho\wedge M_{\rho'}=(M^\rho\wedge M_{\rho'})\vee(M_\rho\wedge M^{\rho'})$, we have $M_\rho\wedge M_{\rho'}(a)=b$. Similarly, we can deduce that $M^\rho\wedge M^{\rho'}(a)=b$, and we have the desired equality.
\hfill\Halmos

\subsubsection{Proof of Theorem~\ref{thm:MSSS}}\label{ec:proofs-siblings-1-main} 
First, fix $(a_i,\bar{a}_i)\in\cC$, we verify that $(\Pi(f^i, \cF))$ satisfies the conditions of Theorem~\ref{thm:main1}: (i) is trivial since $\cF=\cS(I)$. Next, given Lemma~\ref{lem:differentials-sibling}, (ii) holds since every second order differential is $0$. Let us verify (iii). Since by Lemma~\ref{lem:differentials-sibling}, the second order differentials are always $0$, the goal is to show $f^i(M)=f^i(M_0)+\sum_{\rho\in R_M}\partial f^i_\rho=1+\sum_{\rho\in R_M}\partial f^i_\rho$, where the second equality follows from Assumption~\ref{ass:5cases}.

% Using the notations in Assumption~\ref{ass:5cases}, 

Let $B^i \subset B$ denote the subset of schools for which $\rho_{in}(a_i, \bar{a}_i, b)$ and $\rho_{out}(a_i, \bar{a}_i, b)$ exist (schools that verify case 2 of Lemma~\ref{lem:2cases}). For $b\in B^i$ and $\rho\in\cR(I)$, we define the indicator function $g(b,\rho)$ such that $g(b,\rho)=1$ if $\rho\in \{\rho_{in}(a_i,\bar{a}_i,b),\rho_{out}(a_i,\bar{a}_i,b)\}$, otherwise $g(b,\rho)=0$.

Fix $\rho \in \cR(I)$. When $\partial f_\rho^i\neq 0$, we claim that there exists a unique $b\in B^i$ such that $g(b,\rho)=1$. In fact, if $\partial f_\rho^i= 1$, then $\rho\in\cR_{out}(a_i,\bar{a}_i)$ and $\rho\notin\cR_{in}(a_i,\bar{a}_i)$. Therefore, there exists $b\in B^i$ such that $\rho=\rho_{out}(a_i,\bar{a}_i,b)$. Notice that if $b'\neq b$ also satisfies $\rho=\rho_{out}(a_i,\bar{a}_i,b')$, then in the stable matching $M_\rho$ corresponding to the upper-closed set of rotations $R_\rho:=\{\rho'\in\cR(I):\rho' \triangleright\rho\}$, we have $\rho=\rho_{out}(a_i,\bar{a}_i,b)\notin R_\rho$. By Lemma~\ref{lem:2cases}, it holds that $\rho_{in}(a_i,\bar{a}_i,b)\triangleright\rho$, thus $\rho_{in}(a_i,\bar{a}_i,b)\in R_\rho$. Therefore, we conclude, again by Lemma~\ref{lem:2cases}, that $M_\rho(a_i)=M_\rho(\bar{a}_i)=b$. Using the same argument it must also hold that $M_\rho(a_i)=M_\rho(\bar{a}_i)=b'$, which contradicts the fact that $b\neq b'$. Hence, $\rho\neq\rho_{out}(a_i,\bar{a}_i,b')$ for every $b' \in B^i \setminus \{b\}$. It is also impossible to have $\rho=\rho_{in}(a_i,\bar{a}_i,b')$ for any $b' \in B^i$ as $\rho\notin\cR_{in}(a_i,\bar{a}_i)$. Thus, by definition $g(b,\rho)=1$ and $g(b',\rho)=0$ for every $b'\in B^i \setminus \{b\}$. 

The case when $\partial f_\rho^i=-1$ can be treated similarly. Therefore, $\partial f_\rho^i\neq 0$ implies that $\sum_{b\in B^i}g(b,\rho)=1$ and we have,
\begin{align*}
    \sum_{\rho\in R_M}\partial f_\rho^i &=\sum_{\rho\in R_M}\left(\partial f_\rho^i\cdot\sum_{b\in B^i}g(b,\rho)\right) =\sum_{\rho\in R_M}\sum_{b\in B^i} \partial f_\rho^i\cdot g(b,\rho) =\sum_{b\in B^i}\sum_{\rho\in R_M}\partial f_\rho^i\cdot g(b,\rho) \\
    &=\sum_{b\in B^i}\sum_{\rho\in R_M}(\mathbbm{1}(\rho=\rho_{out}(a_i,\bar{a}_i,b))-\mathbbm{1}(\rho=\rho_{in}(a_i,\bar{a}_i,b))).
\end{align*}

If $f^i(M)=1$, then by Lemma~\ref{lem:2cases}, there is no $b\in B^i$ such that $\rho_{in}(a_i,\bar{a}_i,b)\in R_M$ and $\rho_{out}(a_i,\bar{a}_i,b)\notin R_M$. Therefore, for every $b\in B^i$, either $\{\rho_{in}(a_i,\bar{a}_i,b),\rho_{out}(a_i,\bar{a}_i,b)\}\subseteq R_M$ or $\{\rho_{in}(a_i,\bar{a}_i,b),\rho_{out}(a_i,\bar{a}_i,b)\}\subseteq \cR(I)\setminus R_M$, and both lead to $\sum_{\rho\in R_M}(\mathbbm{1}(\rho=\rho_{out}(a_i,\bar{a}_i,b))-\mathbbm{1}(\rho=\rho_{in}(a_i,\bar{a}_i,b)))=0$. 
Therefore, we have $f^i(M)=1+\sum_{\rho\in R_M}\partial f^i_\rho$.

If $f^i(M)=0$, then $M(a_i)=M(\bar{a}_i)=b^*$ for some $b^*\in B^i$. By Lemma~\ref{lem:2cases}, 
for every $b \in B^i \setminus \{b^*\}$, either $\{\rho_{in}(a_i,\bar{a}_i,b),\rho_{out}(a_i,\bar{a}_i,b)\}\subseteq R_M$ or $\{\rho_{in}(a_i,\bar{a}_i,b),\rho_{out}(a_i,\bar{a}_i,b)\}\subseteq \cR(I)\setminus R_M$, and hence $\sum_{\rho\in R_M}(\mathbbm{1}(\rho=\rho_{out}(a_i,\bar{a}_i,b))-\mathbbm{1}(\rho=\rho_{in}(a_i,\bar{a}_i,b)))=0$. Again 
by Lemma~\ref{lem:2cases}, since $M(a_i)=M(\bar{a}_i)=b^*$, we have $\rho_{in}(a_i,\bar{a}_i,b^*)\in R_M$ and $\rho_{out}(a_i,\bar{a}_i,b^*)\notin R_M$, thus $\sum_{\rho\in R_M}(\mathbbm{1}(\rho=\rho_{out}(a_i,\bar{a}_i,b^*))-\mathbbm{1}(\rho=\rho_{in}(a_i,\bar{a}_i,b^*)))=-1$. Hence, it holds that $\sum_{\rho\in R_M}\partial f_\rho^i=-1$, which gives the desired identity for this case as well.

We have shown that condition (iii) in Theorem~\ref{thm:main1} holds. Therefore, ($\Pi(f^i,\cS(I))$) is minimum cut representable. Moreover, it given that the second order differentials are all zero it is also an instance of (MWSM) by Theorem~\ref{thm:linear-to-second-order}. A witness digraph of ($\Pi(f^i,\cS(I))$) can be constructed in polynomial time. In fact, to constrcut the witness digraph $\cD^{f^i, \cF}$ given by cut graph of $(f^i, \cF)$ (Definition~\ref{def:Cut-Graph}) one needs to compute the second order differentials (which are all zero) and the first order differentials which can also be computed in polynomial time given Lemma~\ref{lem:differentials-sibling}. In fact, from Corollary~\ref{cor:all-rotations}, starting from $M_0$ and eliminating exposed rotations in any order until getting to $M_z$ will eliminate all rotations, hence, in order to check whether a rotation $\rho$ belongs or does not belong to $\cR_{in}(a_i, \bar{a}_i)$ and $\cR_{out}(a_i, \bar{a}_i)$, one can simply look at whether $\rho$ moves $a_i, \bar{a}_i$ into or out of the same school in such chain of stable matchings. By Lemma~\ref{lem:additivity} and since $f = \sum_{i \in m} f^i$,~\ref{problem} is minimum cut representable (and more specifically an instance of (MWSM)) whose witness digraph can be constructed efficiently, which implies the Theorem.\hfill\Halmos

\subsection{Matching Siblings with Same Preferences to the Same Activities (MSSP): Omitted Proofs}\label{ec:MSSP-proofs}

\subsubsection{Proof of Lemma~\ref{lem:minimum cut-representable-after-school}}
Let us first consider the following \emph{Implied Rotation Problem} (IRP): we are given a school matching instance $I'=(A,B,\succ,q)$ and $\theta, \overline \theta \in {\cal R}(I') \cup \{\emptyset,\infty\}$ such that at least one of $\theta$, $\overline \theta$ is a rotation, i.e, $\{\theta, \overline \theta\} \not\subseteq \{\emptyset, \infty\}$.
 The goal is to find a stable matching $M \in {\cal S}(I')$ such that:
\begin{itemize}
    \item if both $\theta$ and $\overline \theta$ are rotations, i.e., $\theta,\overline \theta \neq \emptyset, \infty$, then either $\theta, \overline{\theta} \in R_M$ or $\theta, \overline{\theta} \notin R_M$;
    \item else if $\theta = \emptyset$ (resp. $\overline \theta = \emptyset$) then $\overline \theta \in R_M$ (resp. $\theta \in R_M$).
    \item else if $\theta = \infty$ (resp. $\overline \theta = \infty$) then $\overline \theta \notin R_M$ (resp. $\theta \notin R_M$).
\end{itemize}

\begin{lemma}\label{lem:implied-rotation}
Let $J$ be an IRP instance. There exists a function $f: {\cal S}(I')\rightarrow \mathbb{Z}_{\geq 0}$ such that, for each $M \in {\cal S}(I')$, $f(M)=0$ if and only if $M$ is a feasible solution for $J$. Moreover, $\Pi(f,{\cal S}(I'))$ is minimum cut representable and its witness graph can be constructed in polynomial time.
\end{lemma}
\proof{Proof.} 
In the first case, i.e., $\theta, \overline \theta \neq \emptyset, \infty$, we can consider the cut graph constructed in Definition~\ref{def:Cut-Graph} for any arbitrary objective function and the feasibility set $\cF = \cS(I')$, where we only keep infinite capacity arcs (note that these do not depend on the objective function). Then we add to this graph an arc from $\theta$ to $\overline{\theta}$ and an arc from $\overline{\theta}$ to $\theta$ both of capacity $1$. An $s$-$t$ cut in this graph has finite value if and only if it is an upper-closed subset of rotations and has value $0$ if $\theta$ and $\overline \theta$ are either both inside the cut or both outside. Let ${\sf val}(R)$ denote the value of the $s$-$t$ cut $\{s\} \cup R$ in this digraph. Then the function $f$ such that $f(M) = {\sf val}(R_M)$ for every $M \in \cS(I')$, which is by construction minimum cut representable, is such that $f(M)=0$ if and only if either $\theta, \overline{\theta} \in R_M$ or $\theta, \overline{\theta} \notin R_M$.

In the second case, i.e., one of $\theta$ or $\overline{\theta}$ is $\emptyset$, say without loss of generality $\theta=\emptyset$, we can again consider the cut graph constructed in Definition~\ref{def:Cut-Graph} for any arbitrary objective function and the feasibility set $\cF = \cS(I')$, where we only keep infinite capacity arcs. Then we add to this graph an arc from $s$ to $\overline \theta$ of capacity $1$. An $s$-$t$ cut in this graph has finite value if and only if it is an upper-closed subset of rotations and has value $0$ if and only if $\overline \theta$ is in the cut. Let ${\sf val}(R)$ denote the value of the $s$-$t$ cut $\{s\} \cup R$ in this digraph. Then the function $f$ such that $f(M) = {\sf val}(R_M)$ for every $M \in \cS(I')$, which is by construction minimum cut representable, is such that $f(M)=0$ if and only if $\overline{\theta} \in R_M$.

In the third case, i.e., one of $\theta$ or $\overline{\theta}$ is $\infty$, say without loss of generality $\theta=\infty$, we can again consider the cut graph constructed in Definition~\ref{def:Cut-Graph} for any arbitrary objective function and the feasibility set $\cF = \cS(I')$, where we only keep infinite capacity arcs. Then we add to this graph an arc from $\overline \theta$ to $t$ of capacity $1$. An $s$-$t$ cut in this graph has finite value if and only if it is an upper-closed subset of rotations and has value $0$ if and only if $\overline \theta$ is outside the cut. Let ${\sf val}(R)$ denote the value of the $s$-$t$ cut $\{s\} \cup R$ in this digraph. Then the function $f$ such that $f(M) = {\sf val}(R_M)$ for every $M \in \cS(I')$, which is by construction minimum cut representable, is such that $f(M)=0$ if and only if $\overline{\theta} \notin R_M$.

Clearly, all of the above witness digraphs can be constructed in polynomial time. \hfill\Halmos
\endproof

\iffalse Then
$$\partial f_{\rho} = 
\begin{cases}
-1 & \hbox{ if $\rho = \theta$,}\\
0 & \hbox{otherwise.}
\end{cases} \quad \hbox{and} \quad
\partial^2 f_{\rho,\rho'}=\begin{cases}
0 & \hbox{if $\rho \neq \theta, \theta'$} \\

\end{cases}.$$

By symmetry, the only case left is $\theta\in R_M$, $\theta' =\infty$.

\fi \iffalse
Let $f: {\cal S}(I')\rightarrow \{0,1\}$ be defined as follows: $f(M)=1$ if $R_M \cap \{\theta,\overline \theta\}=1$; $f(M)=0$ otherwise. Assume first $\theta \succeq \overline \theta$. Then
$$\partial f_{\rho} = 
\begin{cases}
1 & \hbox{ if $\rho = \theta$,}\\
-1 & \hbox{ if $\rho = \overline \theta$,}\\
0 & \hbox{otherwise.}
\end{cases} \quad \hbox{and} \quad
\partial^2 f_{\rho,\rho'}=\begin{cases}
e    
\end{cases}.$$
We are left to show that $f(M)$
\fi

To ease the notation, we omit the subscript/superscript $i$ in $(a_i, \bar{a}_i)$ and in particular denote the fixed pair of siblings by $(a,\overline a)$. We say that an activity is stable for a student if there is a stable matching where the student is matched to a class corresponding to that activity. We also say that a student $x$ is matched to an activity $c_j$ in matching $M$ if $x$ is matched to a class corresponding to that activity. Let $L_a \subset \{c_1, \dots, c_k\}$ be the subset of stable activities for $a$. Define $L_{\overline a}$ similarly. Note that we can assume without loss of generality that $|L_a \cap L_{\overline{a}}| \geq 2$, otherwise, either $|L_a \cap L_{\overline{a}}|=0$ in which case no stable matching matches $a$ and $\overline{a}$ to the same activity and Lemma~\ref{lem:minimum cut-representable-after-school} clearly holds by considering the linear function where weight $1$ is given to every pair of student and class, or $|L_a \cap L_{\overline{a}}|=1$, in which case there is only one activity $c_t$ to which $a$ and $\overline{a}$ can be both matched at the same time and Lemma~\ref{lem:minimum cut-representable-after-school} clearly holds again by considering the linear function where weight $1$ is given to every pair $ab$ (resp. $\overline a \overline b$) where class $b$ (resp. $\overline b$) corresponds to an activity other than $c_t$ and gives weight $0$ to the rest of the pairs. In both of these cases the resulting problems are minimum cut representable whose witness digraph can be constructed in polynomial time as they are instances of (MWSM).

The next lemma relates (MSSP) to a family of IRP problems. 

\begin{lemma}\label{lem:from-implied-rotation-to-MSSP}
There exists a family $\cK$ of at most $k+1$ IRP instances on $I$ such that, for each $M \in \cS(I)$, $M$ matches $a$ and $\overline a$ to same-activity classes if and only if $M$ is a feasible solution to all the IRP instances. Moreover, $\cK$ can be constructed in polynomial time.
\end{lemma}
%\afcomment{this is now checked}
\proof{Proof.} Without loss of generality, assume the siblings $a,\bar a$ rank activities in the order $c_1,\dots, c_k$. Let $Q(a)$ denote the set of indices $\ell$ for which there exists a rotation $\theta^{a, \ell}$ whose elimination matches $a$ to activity $c_{\ell}$. Define $Q(\overline{a})$ and $\theta^{\overline{a}, \ell}$ analogously. Note that because $|L_a| \geq 2$, there must exist a rotation whose elimination changes the activity to which $a$ is matched, and hence $Q(a) \neq \emptyset$. Similarly, $Q(\overline a) \neq \emptyset$.

For every $\ell \in [k]$, such that $\ell \in Q(a) \cup Q(\overline{a})$, define the pair of rotations $\theta_\ell=\theta^{a,\ell_a}$ and $\overline \theta_\ell=\theta^{\overline a,\ell_{\overline a}}$, where $\ell_a$ (resp., $\ell_{\overline a}$) is the smallest integer $r \geq \ell$ such that $r \in Q(a)$ (resp., $r \in Q(\overline a)$), and $\theta_\ell=\infty$ (resp., $\overline \theta_\ell=\infty$) if such $r$ does not exist. In words, when $\theta_\ell = \infty$, then there is no rotation whose elimination matches $a$ to an activity worse or equal to $c_{\ell}$, and hence, $a$ is always matched to an activity strictly better than $c_\ell$ since $Q(a) \neq \emptyset$. Otherwise, $\theta_\ell$ is the rotation whose elimination matches $a$ to the best activity for $a$ among all stable activities that are worse or equal to $c_{\ell}$.

Let $c_j$ (resp., $c_{\bar \jmath}$) be the activity $a$ (resp., $\overline a$) is matched to in $M_0$. If $j \neq \bar \jmath$, define an additional pair $\theta_0, \overline \theta_0$  as follows. If $j<\bar \jmath$, let $\theta_0$ be the rotation whose elimination takes $a$ from an activity strictly better than $c_{\bar \jmath}$ to an activity worse or equal to $c_{\bar \jmath}$. Note that such rotation exists because $a$ is matched to activity $c_j$ which is strictly better than $c_{\bar \jmath}$ in $M_0$ and matched to an activity $c_{{j}'}$ worse or equal to $c_{\bar \jmath}$ in $M_z$ since $|L_a \cap L_{\overline a}| \geq 1$. Moreover, let $\overline \theta_0=\emptyset$; else if $\bar \jmath < j$, we let $\overline \theta_0$ be the rotation whose elimination takes $\overline a$ from an activity strictly better than $c_j$ to one that is worse or equal to $c_j$, which similarly, exists, and let $\theta_0 = \emptyset$.

Let $\cK$ be the family of IRP instances corresponding to the pairs $\theta_l, \overline{\theta}_l$ defined above.

To show one direction of the lemma, assume that $M$ is not a feasible solution for the IRP instance $I_\ell$ corresponding to the pair $\theta_\ell,\overline \theta_\ell$ for some $\ell \in \{0, \dots, k\}$. We show that $M$ does not match $a$ and $\overline{a}$ to the same activity. Assume first that $\theta_\ell\notin \{\emptyset, \infty\}$ and $\overline\theta_\ell \notin \{\emptyset, \infty\}$ and, without loss of generality, that $\theta_\ell \in R_M$, $\overline \theta_\ell \notin R_M$. $\theta_\ell \in R_M$ implies that $a$ is matched in $M$ to a class corresponding to an activity $c_t$ worse or equal to $c_\ell$, i.e., such that $t \geq \ell$; $\overline\theta_\ell \notin R_M$ that $\overline a$ is matched in $M$ to a class corresponding to an activity $c_t$ strictly better than $c_\ell$, i.e., such that $t<\ell$. Since $a,\overline a$ share the same preferences over activities, they are not matched to same-activity classes in $M$, as required. Now suppose that one of $\theta_\ell, \overline \theta_\ell$ is equal to $\emptyset$ - say, without loss of generality, $\overline \theta_\ell = \emptyset$ and $\theta_\ell \notin \{\emptyset, \infty\}$. Then, $M$ being infeasible for $I_0$ implies $\theta_\ell \notin R_M$. Thus, $a$ is matched in $M$ to an activity $c_t$ strictly better than $c_{\bar \jmath}$, i.e., such that $t<\bar \jmath$. On the other hand, and by definition of $\overline{j}$, $\overline a$ is matched in $M$ to an activity $c_{\overline t}$ with $\overline t \geq \bar \jmath$. In this case again we conclude that $a,\overline a$ are not matched to same-activity classes. Last, suppose that one of $\theta_\ell, \overline \theta_\ell$ is equal to $\infty$ - say, without loss of generality, $\overline \theta_\ell = \infty$ and $\theta_\ell \notin \{\emptyset, \infty\}$. Then, since $\overline \theta_\ell = \infty$, $\overline a$ is matched to an activity strictly better than $c_\ell$. On the other hand, $M$ being infeasible for $I_\ell$ implies that $\theta_\ell \in R_M$ and hence $a$ is matched in $M$ to an activity $c_t$ that is worse or equal to $c_\ell$, i.e., such that $t\geq \ell$, and again we conclude that $a,\overline a$ are not matched to same-activity classes in $M$.

To show the opposite direction of the lemma, suppose that $a, \overline a$ are not matched to same-activity classes in $M$. We show that $M$ is not feasible for some IRP instance in $\cK$. Let $a$ be matched to activity $c_t$, $\overline a$ to activity $c_{\overline t}$, and assume without loss of generality that $\overline t < t$. Recall that $a$ (resp., $\overline a$) is matched in $M_0$ to activity $c_j$ (resp., $c_{\bar \jmath}$). There are two cases, $\overline t < j$, in which case $\bar \jmath \leq \overline t< j \leq t$ or $\overline t \geq j$, in which case $\bar \jmath, j \leq \overline t < t$. In the first case when $\bar \jmath \leq \overline t< j \leq t$, %Then, $a$ and $\overline a$ are matched to different activities in $M_0$, with $\overline a$ being matched to an activity they prefer. 
$I_0$ is well defined and $M$ is feasible for $I_0$ if and only if $R_M$ contains rotation $\overline{\theta}_0$. By definition, $\overline{\theta}_0 \in R_M$ implies 
 $\overline t\geq j$. Thus, $M$ is not feasible for $I_0$. We now consider the second case when $\bar \jmath, j \leq \overline t < t$. Assume first that there is no stable activity $c_{\overline t'} \in L_{\overline a}$ for $\overline a$ worse or equal to $c_t$, i.e., such that $\overline t'\geq t$. Then, $\overline{\theta}_t = \infty$ and hence any stable matching $M'$ such that $a$ is matched to an activity $c_{t'}$ with $t'\geq t$ (i.e., such that $\theta_t \in R_{M'}$) is not feasible for $I_t$. In particular, $M$ is not feasible for $I_\ell$.  We are left to consider the case when there exists a stable activity $c_{\overline{t}'} \in L_{\overline{a}}$ for $\overline a$ such that $\overline t' \geq t$. In this case, there exists a rotation whose elimination matches $\overline a$ to a worse or equal activity to $c_t$. Hence, $\overline{\theta}_t \neq \infty$. Moreover, since $a$ is matched to activity $c_t$ in $M$ and $t>j$, there exists a rotation whose elimination matches $a$ to a worse or equal activity to $c_t$ implying that $\theta_t \neq \infty$ as well.
 Since $\overline t <t$, $\overline \theta_t \notin R_M$. On the other hand, $\theta_t \in R_M$. Thus, $M$ is not feasible for $I_t$, concluding the proof. \hfill\Halmos
\endproof

We can now conclude the proof of Lemma~\ref{lem:minimum cut-representable-after-school}. Construct the family of instances $\cK$ as in Lemma~\ref{lem:from-implied-rotation-to-MSSP} and the corresponding functions $\{f_J\}_{J \in \cK}$, where $f_J$ is the function obtained by applying Lemma~\ref{lem:implied-rotation} to instance $J \in \cK$. Let $f^i=\sum_{J \in \cK} f_J$. By Lemma~\ref{lem:from-implied-rotation-to-MSSP}, we deduce that, for each $M \in {\cal S}(I)$, $f^i(M)=0$ if and only if $M$ matches $a,\overline a$ to the same activity. This shows 1. To conclude that $2$ holds, we apply Lemma~\ref{lem:implied-rotation}, Lemma~\ref{lem:additivity} and Lemma~\ref{lem:from-implied-rotation-to-MSSP}  to deduce that $\Pi(f^i,S)$ is minimum cut representable whose witness digraph can be constructed in polynomial time. \hfill \Halmos

\subsubsection{Proof of Lemma~\ref{lem:counter-ex-not-MWSS}.} 
Consider again the (MSSP) instance from Figure~\ref{fig:MSSP}. Notice that $M_0 \geq M_1, M_2 \geq M_z$, while $M_1, M_2$ are incomparable. Thus, the poset of rotations of the instance contains exactly two incomparable rotations $\rho_1,\rho_2$. 
Recall that $M_0$, $M_z=M/\{\rho_1,\rho_2\}$ are the only activity-stable matchings. Let $f: \{M_0,M_1,M_2,M_z\}\rightarrow \mathbb{Z}_{\geq 0}$ be a minimum cut representable function such that $f(M_0)=f(M_z)<f(M_1),f(M_2)$. If $f$ is linearizable, by Theorem~\ref{thm:linear-to-second-order} we know $\partial^2 f_{\rho_1,\rho_2}=0$. Using Theorem~\ref{thm:main1}, we deduce that 
$$f(M_1)=f(M_0)+w_1; \quad f(M_2)=f(M_0)+w_2; \quad f(M_z)=f(M_0)+w_1+w_2,$$
for some $w_1,w_2$. Hence, $f(M_0)+f(M_z)=f(M_1)+f(M_2)$,
contradicting the fact that $f(M_0)$, $f(M_z)$ are the only minima of $f$. \hfill \Halmos

\subsubsection{Proof of Lemma~\ref{lem:generalization-mssp-not-minc-representable}} %\afcomment{this is now checked} 
Consider the instance of (MSSP) in Figure~\ref{fig:MSSP-counting}. Let $M_1=\{a_1b_1,a_2b_2,a_3b_3\}$, $\overline{M}_1=\{\bar{a}_1\bar{b}_1,\bar{a}_2\bar{b}_2,\bar{a}_3\bar{b}_3\}$, $M_2=\{a_1b_2,a_2b_3,a_3b_1\}$, $\overline{M}_2=\{\bar{a}_1\bar{b}_2,\bar{a}_2\bar{b}_3,\bar{a}_3\bar{b}_1\}$, $M_3=\{a_1b_3,a_2b_1,a_3b_2\}$, $\overline{M}_3=\{\bar{a}_1\bar{b}_3,\bar{a}_2\bar{b}_1,\bar{a}_3\bar{b}_2\}$. One can verify that $M_i\cup \overline{M}_j$ for any $i,j\in\{1,2,3\}$ are all stable matchings and the only stable matchings. In addition, $M_i\cup\overline{M}_i$ ($i=1,2,3$) are the activity-stable matchings, while the other stable matchings match no pairs of siblings to the same activity. 

Consider a function $f$ such that for any stable matchings $M,M'$, $f(M)<f(M')$ if and only if $M$ matches strictly more pairs of siblings to the same activity than $M'$. Then we must have
\begin{align*}
    f(M)=\left\{\begin{array}{cc}
      x   & \quad  M=M_i\cup \overline{M}_i\;\textrm{ for some\; $i = 1, 2, 3$} \\
      y   &  \textrm{otherwise}
    \end{array}\right.
\end{align*}
for some $x<y$. On the other hand, the school matching instance corresponding to the above (MSSP) instance has $4$ rotations: $$\rho_1=(\{a_1b_1,a_2b_2,a_3b_3\},\{a_1b_2,a_2b_3,a_3b_1\}), \rho_2=(\{a_1b_2,a_2b_3,a_3b_1\},\{a_1b_3,a_2b_1,a_3b_2\}),$$ $$\bar{\rho}_1=(\{\bar{a}_1\bar{b}_1,\bar{a}_2\bar{b}_2,\bar{a}_3\bar{b}_3\},\{\bar{a}_1\bar{b}_2,\bar{a}_2\bar{b}_3,\bar{a}_3\bar{b}_1\}), \bar{\rho}_2=(\{\bar{a}_1\bar{b}_2,\bar{a}_2\bar{b}_3,\bar{a}_3\bar{b}_1\},\{\bar{a}_1\bar{b}_3,\bar{a}_2\bar{b}_1,\bar{a}_3\bar{b}_2\}),$$ with $\rho_1 \triangleright \rho_2$ and $\bar{\rho}_1 \triangleright \bar{\rho}_2$, and the other pairs of rotations being incomparable. 

One easily checks that $M_{\rho_1}=M_1 \cup \overline{M_1}$, $M^{\rho_1}= M_2 \cup \overline{M_1}$, $M_{\bar{\rho}_2}=M_1 \cup \overline{M}_2$, and $M^{\bar{\rho}_2}=M_1 \cup \overline{M}_3$.
Hence, by Definition~\ref{def:sod}, $\partial^2 f_{\rho_1,\bar{\rho}_2}=f(M_2\cup\overline{M}_2)+f(M_1\cup\overline{M}_3)-f(M_1\cup\overline{M}_2)-f(M_2\cup\overline{M}_3)=x-y<0,$ implying by Theorem~\ref{thm:main1} that $f$ is not minimum cut representable. \hfill \Halmos

\begin{figure}[t!]
\begin{center}
\begin{tabular}{c  l c | c c l }
$a_1$: & $b_1 \AgentGreater{a_1}{}  b_2   \AgentGreater{a_1}{} b_3 \AgentGreater{a_1}{} \emptyset$ & \hspace{.15cm} &  & $b_1$: & $a_2 \AgentGreater{b_1}{} a_3 \AgentGreater{b_1}{} a_1 \AgentGreater{b_1}{} \emptyset$ \\ 
$a_2$: & $b_2 \AgentGreater{a_2}{}  b_3   \AgentGreater{a_2}{} b_1\AgentGreater{a_2}{} \emptyset$ & \hspace{.15cm} &  & $b_2$: & $a_3 \AgentGreater{b_2}{} a_1 \AgentGreater{b_2}{} a_2 \AgentGreater{b_2}{} \emptyset$ \\ 
$a_3$: & $b_3 \AgentGreater{a_3}{}  b_1   \AgentGreater{a_3}{} b_2 \AgentGreater{a_3}{} \emptyset$ & \hspace{.15cm} &  & $b_3$: & $a_1 \AgentGreater{b_3}{} a_2 \AgentGreater{b_3}{} a_3 \AgentGreater{b_3}{} \emptyset$ \\ 
$\bar{a}_1$: & $\bar{b}_1 \AgentGreater{\bar{a}_1}{}  \bar{b}_2   \AgentGreater{\bar{a}_1}{} \bar{b}_3 \AgentGreater{\bar{a}_1}{} \emptyset$ & \hspace{.15cm} &  & $\bar{b}_1$: & $\bar{a}_2 \AgentGreater{\bar{b}_1}{} \bar{a}_3 \AgentGreater{\bar{b}_1}{} \bar{a}_1 \AgentGreater{\bar{b}_1}{} \emptyset$ \\ 
$\bar{a}_2$: & $\bar{b}_2 \AgentGreater{\bar{a}_2}{}  \bar{b}_3   \AgentGreater{\bar{a}_2}{} \bar{b}_1\AgentGreater{\bar{a}_2}{} \emptyset$ & \hspace{.15cm} &  & $\bar{b}_2$: & $\bar{a}_3 \AgentGreater{\bar{b}_2}{} \bar{a}_1 \AgentGreater{\bar{b}_2}{} \bar{a}_2 \AgentGreater{\bar{b}_2}{} \emptyset$ \\ 
$\bar{a}_3$: & $\bar{b}_3 \AgentGreater{\bar{a}_3}{}  \bar{b}_1   \AgentGreater{\bar{a}_3}{} \bar{b}_2 \AgentGreater{\bar{a}_3}{} \emptyset$ & \hspace{.15cm} &  & $\bar{b}_3$: & $\bar{a}_1 \AgentGreater{\bar{b}_3}{} \bar{a}_2 \AgentGreater{\bar{b}_3}{} \bar{a}_3 \AgentGreater{\bar{b}_3}{} \emptyset$ \\ 
\vspace{3mm}
\end{tabular}
\end{center}
\caption{An instance of (MSSP) with $\cC =\{(a_i,\bar{a}_i)\}_{i=1,2,3}$ and three activities $c_1, c_2, c_3$, with two classes $b_i,\bar{b}_i$ for each activity $c_i$. Each student $a_i$ (resp.,~$\bar{a}_i$) is eligible for class $b_1,b_2,b_3$ (resp., $\bar{b}_1,\bar{b}_2,\bar{b}_3$). All classes have quota of $1$.}\label{fig:MSSP-counting}
\end{figure}

\subsection{Matching Siblings with Different Preferences to Same Activities (MSDP): Omitted Proofs}\label{app:MSDP-omitted}

\subsubsection{Proof of Theorem~\ref{thm:MSDP-Np-Complete}} %\afcomment{this is now checked}

We reduce from the NP-complete \emph{Jointly Stable Matching} (JSM) problem (see~\cite{miyazaki2019jointly}). In (JSM), we are given two school matchings instances $(A,B,\succ^1,q)$, $(A,B,\succ^2,q)$, over the same sets of students and schools, but with possibly different preference lists, and all schools having a quota of $1$ . The goal is to decide if there exists a matching that is stable in both instances. Given an instance of (JSM), define the following instance (MSDP) over the school matching instance $I = (A', B', \succ', q')$ and siblings $\cC$. For each $a \in A$, create a pair of siblings $(a',\overline a')$. For every $b \in B$, create an activity $c$ and two classes $b',\overline b'$ that correspond to $c$. Every class has a quota of $1$. For each pair of siblings $(a',\overline a')$, the preference list of $a'$ (resp., $\overline a'$) is obtained from the preference list of $a$ in $\succ^1$ (resp., $\succ^2$), by replacing $b_i$ with $b_i'$ (resp., $\overline b_i'$). Thus, the class corresponding to activity $c$ for which $a'$ (resp., $\overline a'$) is eligible is $b'$ (resp., $\overline b'$). For each $b \in B$, the preference list of $b' \in B'$ is obtained by setting $a_1'\succ'_{b'} a_2'$ if and only if $a_1 \succ^1_b a_2$. Similarly, the preference list of $\overline b' \in B'$ is obtained by setting $\overline a_1'\succ'_{b'} \overline a_2'$ if and only if $a_1 \succ^2_b a_2$. 

Let $M$ be a matching that is stable in both $(A,B,\succ^1,q)$ and $(A,B,\succ^2,q)$. Construct a matching $M'$ for the (MDSP) instance defined above by adding, for each pair $ab 
\in M$, pairs $a'b'$ and $\overline a'\overline b'$ to $M'$. We claim that $M'$ is activity-stable. First observe that each pair of siblings is matched to a class corresponding to the same activity. Thus, we only need to show that $M'$ is stable. By contradiction, assume without loss of generality that $a'b'$ is a blocking pair for $M'$. Equivalently, $b' \succ'_{a'} M'(a')$ and $a' \succ'_{b'} M'(b')$. Thus, $b\succ^1_{a} M(a)$ and $a \succ^1_{b} M(a)$, contradicting the stability of $M$ in $(A,B,\succ^1,q)$. 

Conversely, assume that $M'$ is an activity-stable matching in the (MSDP) instance. In particular, for each pair of siblings $(a',\overline a')$ and same-activity classes $b',\overline b'$, we have $a'b' \in M'$ if and only if $\overline a'\overline b' \in M'$. Create a matching $M$ over $A\times B$ by adding $ab$ to $M$ for each $a'b' \in M'$. We claim that $M$ is stable in both $(A,B,\succ^1,q)$ and $(A,B,\succ^2,q)$. Suppose by contradiction that it is not stable in the former. There exists therefore $ab \in M$ such that $b \succ^1_{a} M(a)$ and $a\succ^1_{b} M(b)$. Then, $b'\succ'_{a'} M'(a')$ and $a'\succ'_{b'} M'(b')$, a contradiction. The case when $M$ is not stable in $(A,B,\succ^2,q)$ follows analogously. \hfill \Halmos

\subsubsection{Proof of Lemma~\ref{lem:counter-ex-not-MSDP}.} 

\begin{figure}[t!]
\begin{center}
\begin{tabular}{c  l c | c c l }
$a_1$: & $b_1 \AgentGreater{a_1}{}  b_2   \AgentGreater{a_1}{} b_3 \AgentGreater{a_1}{} \emptyset$ & \hspace{.15cm} &  & $b_1$: & $a_2 \AgentGreater{b_1}{} a_1 \AgentGreater{b_1}{} \emptyset \AgentGreater{b_1}{} a_3$ \\ 
$a_2$: & $b_2 \AgentGreater{a_2}{}  b_1   \AgentGreater{a_2}{} \emptyset \AgentGreater{a_2}{} b_3$ & \hspace{.15cm} &  & $b_2$: & $a_3 \AgentGreater{b_2}{} a_1 \AgentGreater{b_2}{} a_2 \AgentGreater{b_2}{} \emptyset$ \\ 
$a_3$: & $b_3 \AgentGreater{a_3}{}  b_2   \AgentGreater{a_3}{} \emptyset \AgentGreater{a_3}{} b_1$ & \hspace{.15cm} &  & $b_3$: & $a_1 \AgentGreater{b_3}{} a_3 \AgentGreater{b_3}{} \emptyset  \AgentGreater{b_3}{} a_2$ \\ 
$\bar{a}_1$: & $\bar{b}_3 \AgentGreater{\bar{a}_1}{}  \bar{b}_2   \AgentGreater{\bar{a}_1}{} \bar{b}_1 \AgentGreater{\bar{a}_1}{} \emptyset$ & \hspace{.15cm} &  & $\bar{b}_1$: & $\bar{a}_1 \AgentGreater{\bar{b}_1}{} \bar{a}_3 \AgentGreater{\bar{b}_1}{} \emptyset \AgentGreater{\bar{b}_1}{} \bar{a}_2$ \\ 
$\bar{a}_2$: & $\bar{b}_2 \AgentGreater{\bar{a}_2}{}  \bar{b}_3   \AgentGreater{\bar{a}_2}{} \emptyset \AgentGreater{\bar{a}_2}{} \bar{b}_1$ & \hspace{.15cm} &  & $\bar{b}_2$: & $\bar{a}_3 \AgentGreater{\bar{b}_2}{} \bar{a}_1 \AgentGreater{\bar{b}_2}{} \bar{a}_2 \AgentGreater{\bar{b}_2}{} \emptyset$ \\ 
$\bar{a}_3$: & $\bar{b}_1 \AgentGreater{\bar{a}_3}{}  \bar{b}_2   \AgentGreater{\bar{a}_3}{} \emptyset \AgentGreater{\bar{a}_3}{} \bar{b}_3$ & \hspace{.15cm} &  & $\bar{b}_3$: & $\bar{a}_2 \AgentGreater{\bar{b}_3}{} \bar{a}_1 \AgentGreater{\bar{b}_3}{} \emptyset \AgentGreater{\bar{b}_3}{} \bar{a}_3$ \\ 
\vspace{3mm}
\end{tabular}
\end{center}
\caption{An instance of (MSDP) with $\cC =\{(a_1,\bar{a}_1)\}$, and three activities $c_1, c_2, c_3$, with two classes $b_i,\bar{b}_i$ for each activity $c_i$. Student $a_1$ (resp.,~$\bar{a}_1$) is eligible for class $b_1,b_2,b_3$ (resp., $\bar{b}_1,\bar{b}_2,\bar{b}_3$). All classes have quota of $1$.}\label{fig:MSDP-counterexample}
\end{figure} 

Consider the (MSDP) instance given in Figure~\ref{fig:MSDP-counterexample}. Define 
$$M_1=\{a_1b_1,a_2b_2,a_3b_3\}, M_2=\{a_1b_2,a_2b_1,a_3b_3\}, %M_3=\{a_1b_3,a_2b_1,a_3b_2\},$$$$\overline{M}_1=\{\bar a_1\bar b_3,\bar a_2\bar b_2,\bar a_3\bar b_1\}, 
\overline{M}_2=\{\bar a_1\bar b_2,\bar a_2\bar b_3,\bar a_3\bar b_1\}, \overline{M}_3=\{\bar a_1\bar b_1,\bar a_2\bar b_3,\bar a_3\bar b_2\}.$$ One easily verifies that $M_2 \cup \overline M_2$, $M_1 \cup \overline{M}_3$ are activity-stable matchings, while their join $M_1 \cup \overline M_2$ is not activity-stable. The claim then follows from Lemma~\ref{lem:if-mincut-then-lattice}. \hfill \Halmos

\iffalse The poset of rotations is composed of for elements:
$$\rho_1=(\{a_1b_1,a_2b_2\},\{a_1b_2,a_2b_1\}), \rho_2=(\{a_1b_2,a_3b_3\},\{a_1b_3,a_3b_2\}),$$ $$\bar{\rho}_1=(\{\bar{a}_1\bar{b}_3,\bar{a}_2\bar{b}_2\},\{\bar{a}_1\bar{b}_2,\bar{a}_2\bar{b}_3\}), \bar{\rho}_2=(\{\bar{a}_1\bar{b}_2,\bar{a}_3\bar{b}_1\},\{\bar{a}_1\bar{b}_1,\bar{a}_3\bar{b}_2\}),$$ with $\rho_1 \trianglerighteq \rho_2$ and $\bar{\rho}_1 \trianglerighteq \bar{\rho}_2$, and the other pairs of rotations being incomparable. One can verify that $M_{\rho_1}=M_1 \cup \overline{M}_1$, $M^{\rho_1}=M_2 \cup \overline{M}_1$, $M_{\bar{\rho}_2}=M_1 \cup \overline{M}_2$, and $M^{\bar{\rho}_2}=M_1 \cup \overline{M}_3$. Hence, if $f$ is a minimum cut representable function whose domain are the stable matchings of the instance, by Definition~\ref{def:sod} and Theorem~\ref{thm:main1} we have: $$\partial^2 f_{\rho_1,\bar{\rho}_2}=f(M_2\cup\overline{M}_2)+f(M_1\cup\overline{M}_3)-f(M_1\cup\overline{M}_2)-f(M_2\cup\overline{M}_3)\geq 0.$$
Observing that $M_2\cup\overline{M}_2$ and $M_1\cup\overline{M}_3$ are activity-stable while $M_1\cup\overline{M}_2$ is not, we obtain that the set of minima of $f$ does not coincide with the activity-stable matchings, as required.
\fi 

\section{Two-Stage Stochastic Stable Matchings: Omitted Proofs}\label{ap:2-stage}

\subsection{Explicit Second-Stage Distribution: Omitted Proofs}\label{ap:2-stage-exp}

Recall that $I^{\sqcup} = I \sqcup \bigsqcup_{k=1}^K J_k = (A^{\sqcup}, B^{\sqcup}, \succ^{\sqcup}, q^{\sqcup})$ denotes the disjoint union of the first-stage instance $I$ and the second-stage instances $\{J_k\}_{k\in [K]}$. We begin this section with the following definitions.   \begin{definition}[Projection of a matching] For $M \in \cS(I^{\sqcup})$ and $L \in \{I\} \cup \{J_k\}_{k \in [K]}$, define 
        $$
        {\cal P}_{L}(M) = \{ab\}_{\substack{ab \in A^{L} \times B^{L} \\ M((a, A^{L})) = (b, B^{L})}} \cup \{a\emptyset\}_{\substack{ a \in A^{L} \\ M((a, A^{L})) = \emptyset}} \cup \{\emptyset b\}_{\substack{b \in B^{L}\\ M((b, B^{L})) = \emptyset}}.
        $$
        We refer to ${\cal P}_{L}(M)$ as the \emph{projection} of the stable matching $M$ on instance $L$.
\end{definition}
\begin{definition}[Lifting of a matching]
For $M^{I} \in \cS(I)$, $M^{J_k} \in \cS(J_k)$ for $k \in [K]$, define,
        $$
        {\cal L}(M^{I}, M^{J_1}, \dots, M^{J_K})=\bigcup_{L \in \{I\} \cup \{J_k\}_{k \in [K]}}\{(a, A^{L})(b, B^{L})\}_{\substack{ab \in A^{L} \times B^{L} \\ M^{L}(a)=b}} \cup \{(a, A^{L})\emptyset\}_{\substack{a \in A^{L} \\ M^{L}(a)=\emptyset}} \cup \{\emptyset(b, B^{L})\}_{\substack{b \in B^{L} \\ M^{L}(b)=\emptyset}}.
        $$ We refer to ${\cal L}(M^{I}, M^{J_1}, \dots, M^{J_K})$ as the \emph{lifting} of the stable matchings $M^{I}, M^{J_1}, \dots, M^{J_K}$ to instance $I^{\sqcup}$.\end{definition}
       
        The next lemma justifies the names given to the objects above.

\begin{lemma}
\label{lem:projection-lifting} The following holds,
    \begin{itemize}
        \item Let $M \in \cS(I^{\sqcup})$ and $L \in \{I\} \cup \{J_k\}_{k \in [K]}$. Then, ${\cal P}_{L}(M) \in \cS(L)$. 
        \item Let $M^{I} \in \cS(I)$ and $M^{J_k} \in \cS(J_k)$ for every $k \in [K]$. Then ${\cal L}(M^{I}, M^{J_1}, \dots, M^{J_K}) \in \cS(I^{\sqcup})$.
    \end{itemize}
\end{lemma}
% \begin{proof}

{\noindent \em Proof.}
    Given $M \in \cS(I^\sqcup)$, ${\cal P}_{L}(M)$ for $L \in \{I\} \cup \{J_k\}_{k \in [K]}$ is a matching by construction and it is stable because if $ab$ is a blocking pair of ${\cal P}_{L}(M)$ then, by construction, $(a, A^L)(b,B^L)$ would be a blocking pair for $M$. Similarly, if $a$ (resp. $b$) is a blocking student (resp. school) of ${\cal P}_{L}(M)$ then, by construction, $(a,A^L)\emptyset$ (resp. $\emptyset(b,AB^L)$) would be a blocking student (resp. school) of $M$.

    Given $M^{I} \in \cS(I)$ and $M^{J_k} \in \cS(J_k)$ for every $k \in [K]$, ${\cal L}(M^{I}, M^{J_1}, \dots, M^{J_K})$ is a matching by construction and is stable because any blocking pair of ${\cal L}(M^{I}, M^{J_1}, \dots, M^{J_K})$ must involve copies of a pair of agents $(a, A^L)$ and $(b, B^L)$ belonging to the same instance $L \in \{I\} \cup \{J_k\}_{k \in [K]}$ as copies of agents from different instances list each other after the outside option and cannot form a blocking pair. This same pair of agents $ab$ is therefore a blocking pair for the stable matching $M^L$ as $\succ^{\sqcup}$ and $\succ^{L}$ agree over the agents of the instance $L$. Similarly, for any blocking student $(a, A^L)$ (resp. school $(b, B^L)$) of ${\cal L}(M^{I}, M^{J_1}, \dots, M^{J_K})$, $a$ (resp. $b$) would be a blocking student (resp. school) for $M^L$.
% \end{proof}
\hfill $\square$

Next, we state the following identity whose proof is by immediate case study and hence is omitted here.

\begin{lemma}\label{lem:interval}
Let $x,y,z,w\in\mathbb{R}$ such that $x\le y$ and $z\le w$. Then,
$|[x,y] \cap [z,w]|=[z-y]^++[w-x]^+-[z-x]^+-[w-y]^+$, where $|[a,b]|=b-a$ denotes the length of the interval $[a,b]$.
\end{lemma}
We are now ready to prove our results.

\subsubsection{Proof of Lemma~\ref{lem:reformulation}}
Let $M \in \cF$. By Lemma~\ref{lem:projection-lifting}, $M^{I} = {\cal P}_I(M)$  and  $M^{J_k} = {\cal P}_{J_k}(M)$ for every $k \in [K]$ are stable matchings. Hence  $(M^I, \{M^{J_k}\}_{k \in [K]})$ is a feasible solution of~\eqref{eq:exp2sto}. The value of this solution is given by
\begin{align*}
    &c_1(M^I) + \sum_{k\in [K]} p^{J_k}\left(c_2(M^{J_k}) + d(M^I, M^{J_k})\right)\\
    &= c_1^{(I)}(M^{(I)}) + \sum_{k \in [K]} p^{J_k}c^{(J_k)}_2(M^{(J_k)}) + \lambda \sum_{k \in [K]} p^{J_k}\sum_{a \in A^I \cap A^{J_k}} [R_a(M^{(J_k)}(a))-R_a(M^{(I)}(a))]^+\\
    &= c_1^{(I)}(M) + \sum_{k \in [K]} p^{J_k}c^{(J_k)}_2(M) + \lambda \sum_{k \in [K]} p^{J_k}\sum_{a \in A^I \cap A^{J_k}} [R_a(M^{(J_k)}(a))-R_a(M^{(I)}(a))]^+ = f(M)
\end{align*}

Conversely, let $M^I \in \cS(I)$ and $M^{J_k} \in \cS(J_k)$, The lifting $M = {\cal L}(M^{I}, M^{J_1}, \dots, M^{J_K})$ of these stable matchings to $I^{\sqcup}$ is a stable matching, and hence a feasible solution of~\eqref{problem}. Its value is,
\begin{align*}
    f(M)
    &= c_1^{(I)}(M) + \sum_{k \in [K]} p^{J_k}c^{(J_k)}_2(M) + \lambda \sum_{k \in [K]} p^{J_k}\sum_{a \in A^I \cap A^{J_k}} [R_a(M^{(J_k)}(a))-R_a(M^{(I)}(a))]^+\\
    &= c_1^{(I)}(M^{(I)}) + \sum_{k \in [K]} p^{J_k}c^{(J_k)}_2(M^{(J_k)}) + \lambda \sum_{k \in [K]} p^{J_k}\sum_{a \in A^I \cap A^{J_k}} [R_a(M^{(J_k)}(a))-R_a(M^{(I)}(a))]^+\\
    &= c_1(M^I) + \sum_{k\in [K]} p^{J_k}\left(c_2(M^{J_k}) + d(M^I, M^{J_k})\right),
\end{align*}
which is the same as the objective value of $(M^I, \{M^{J_k}\}_{k \in [K]})$ in~\eqref{eq:exp2sto}.
Hence,~\eqref{eq:exp2sto} and~\eqref{problem} are equivalent.
\hfill$\square$

\subsubsection{Proof of Theorem~\ref{thm:two-stage-main1}}
We show that for every $k \in [K]$ and $a \in A^I \cap A^{J_k}$, $(\Pi(f^3_{a,k}, \cF))$ is minimum-cut representable. 

First, recall that in a school matching instance, if a student is matched to the outside option in one stable matching then it is matched to the outside option in all stable matching (see the rural hospital theorem, e.g., \cite[Theorem 1.6.3]{gusfieldbook}). Now, if $a$ is matched to the outside option in the first-stage instance $I$ or in the second-stage $J_k$ or in both, then by Lemma~\ref{lem:projection-lifting}, it must be the case that $M((a, A^I))=\emptyset$ or $M((a, A^{J_k}))=\emptyset$ or both, for every stable matching $M \in \cS(I^\sqcup)$, otherwise the projection of such stable matching to $I$ or $J_k$ or both would have $a$ matched to a partner different from the outside option which is absurd. Hence, $f_{a,k}^3(M)= [R_a(M^{(J_k)}(a))-R_a(\emptyset)]^+ = 0$ or $f_{a,k}^3(M)= [R_a(\emptyset)-R_a(M^{(I)}(a))]^+ = R_a(\emptyset)-R_a(M^{(I)}(a))$ or $f_{a,k}^3(M)= [R_a(\emptyset)-R_a(\emptyset)]^+ = 0$ for such $a$ and $k$. In all of these cases $f^3_{a,k}$ is a linear function of $M$ and because $\cF$ is the lattice of all stable matchings, the corresponding problem $(\Pi(f^3_{a,k},\cF))$ is minimum-cut representable.

Now fix $k \in [K]$ and $a \in A^I \cap A^{J_k}$, such that $a$ is not matched to the outside option in either $I$ nor $J_k$. We show that $f^3_{a,k}$ and $\cF$ verify the conditions of Theorem~\ref{thm:main1}. 
Condition (i) of Theorem~\ref{thm:main1} is trivial as $\cF$ is the lattice of all stable matchings of $I^{\sqcup}$. We now show that (ii) and (iii) hold.

For every instance $L \in \{I\} \cup \{J_k\}_{k \in [K]}$, let $b^{L}_0 \succ^I_{a} b^{L}_1 \succ^I_{a} \dots \succ^I_{a} b^{L}_{m^{L}}$ denote the ordered sequence of stable schools (schools to which $a$ is matched in some stable matching) of $a$ in instance $L$ . Let $\cR_a(I^{\sqcup})$ denote the subset of rotations that change the partner of $(a, A^L)$ for some instance $L \in \{I\} \cup \{J_k\}_{k \in [K]}$. The following claim characterizes the stable matchings $M^{\rho}$ and $M_{\rho}$ for $\rho \in \cR(I^{\sqcup})$.

\begin{claim}\label{cl:rotation-per-instance}
There is a one-to-one correspondence $\Psi$ between the rotations $\cR_a(I^{\sqcup})$ and the 
set of schools $\{(b^L_{j}, B^L)\}_{L \in \{I\} \cup \{J_k\}_{k \in [K]}, j \in [m^{L}]}$, such that for every $\rho \in \cR(I^{\sqcup})$, it holds that
\begin{itemize}
    \item If $\rho \in \cR(I^{\sqcup}) \setminus \cR_a(I^{\sqcup})$ then $M^\rho((a, A^{L'})) = M_\rho((a, A^{L'})) = (b^{L'}_0, B^{L'})$ for every $L' \in \{I\} \cup\{J_k\}_{k \in [K]}$.
    \item If $\rho \in \cR_a(I^\sqcup)$ with $\Psi(\rho) = (b^L_j, B^L)$ then,
    \begin{itemize}
    \item[*] $M^\rho((a, A^{L'})) = M_\rho((a, A^{L'})) = (b^{L'}_0, B^{L'})$ for every $L' \in \{I\} \cup\{J_k\}_{k \in [K]} \setminus \{L\}$.
    \item[*] $M^\rho((a, A^L)) = (b^L_j, B^L)$ and $M_\rho((a, A^L)) = (b^L_{j-1}, B^L)$.
\end{itemize}
\end{itemize}
\end{claim}

{\noindent \em Proof.} We break the proof into three parts.

{\noindent \bf Correspondence between rotations of $I^{\sqcup}$ and rotations of $I, J_1, \dots, J_K$.}
Let $\rho \in \cR(I^{\sqcup})$. For every $L' \in \{I\} \cup \{J_k\}_{k \in [K]}$, let $M_0^{L'}$ denote the student optimal matching of instance $L'$. We show that there exists $L \in \{I\} \cup \{J_k\}_{k \in [K]}$ and $\hat{\rho} \in \cR(L)$ such that $\rho = \rho(M, M')$ where $M = {\cal L}(M_0^{I}, M_0^{J_1}, \dots, M_{\hat{\rho}}, \dots, M_0^{J_K})$ and $M' = {\cal L}(M_0^{I}, M_0^{J_1}, \dots, M^{\hat{\rho}}, \dots, M_0^{J_K})$ where $M_{\hat{\rho}}$ and $M^{\hat{\rho}}$ are at the position corresponding to instance $L$.

Consider $M, M' \in \cS(I^{\sqcup})$ such that $M'$ is an immediate predecessor of $M$ and $\rho = \rho(M, M')$.

First, we claim that eliminating $\rho$ changes the partner of students $(a', A^L)$ for only a unique instance $L$. Indeed, assume there exists two instance $L \neq L'$ such that eliminating $\rho$ changes the partners of both $(a', A^L), (a'', A^{L'}) \in A^\sqcup$. Then $M((a', A^L)) \neq M'((a', A^L))$ and $M((a'', A^{L'})) \neq M'((a'', A^{L'}))$. Consider the matching 
    $$
    M'' = {\cal L}({\cal P}_{I}(M'), {\cal P}_{J_1}(M'), \dots, {\cal P}_{L}(M) ,\dots, {\cal P}_{J_K}(M'))
    $$
    then, by Lemma~\ref{lem:projection-lifting}, $M''$ is a stable matching. Moreover, by construction, $M' <^\sqcup M'' <^\sqcup M$ which contradicts the fact that $M'$ is an immediate predecessor of $M$.

Second, let $L$ be the unique instance where (a subset of) students change partners when $\rho$ is eliminated. We show that ${\cal P}(M'^L)$ is an immediate predecessor of ${\cal P}(M^L)$ in $L$. Indeed, if that is not the case then there would exists $M''^L \in \cS(L)$ such that ${\cal P}(M'^L) <^L M''^L <^L {\cal P}(M^L)$. Then, the lifting 
    $
    M'' = {\cal L}({\cal P}_{I}(M'), {\cal P}_{J_1}(M'), \dots, {\cal P}_{L}(M'') ,\dots, {\cal P}_{J_K}(M'))
    $ is a stable matching such that $M' <^\sqcup M'' <^\sqcup M$ which contradicts the fact that $M'$ is an immediate predecessor of $M$.

Finally, let $\hat{\rho} \in \cR(L)$ such that $\hat{\rho} = \rho({\cal P}_{L}(M), {\cal P}_{L}(M'))$ in $L$. Let
    $$
    \tilde{M} = {\cal L}(M^I_0, M_0^{J_1}, \dots, M_{\hat{\rho}} ,\dots, M^{J_K}_0)
    \quad \text{ and } \quad 
    \tilde{M}' = {\cal L}(M^I_0, M_0^{J_1}, \dots, M^{\hat{\rho}} ,\dots, M^{J_K}_0).
    $$ Then by the first two points, the difference $M \setminus M'$ (resp. $M' \setminus M$) is the same as $\tilde{M} \setminus \tilde{M}'$ (resp. $\tilde{M}' \setminus \tilde{M}$). Hence, $\rho = \rho(\tilde{M}, \tilde{M}')$.

{\noindent \bf The one-to-one correspondence $\Psi$.} Consider $\rho \in \cR_a(I^{\sqcup})$
such that $\rho = \rho(M, M')$ where
    $$
    M = {\cal L}(M^I_0, M_0^{J_1}, \dots, M_{\hat{\rho}} ,\dots, M^{J_K}_0)
    \quad \text{ and } \quad
    M' = {\cal L}(M^I_0, M_0^{J_1}, \dots, M^{\hat{\rho}} ,\dots, M^{J_K}_0),
    $$
for some $L \in \{I\} \cup \{J_k\}_{k \in [K]}$ and $\hat{\rho} \in \cR(L)$.
% such that $\rho = \rho(M, M')$ where $M$ is the lifting of the stable matchings $\{M_0^{L'}\}_{\substack{L' \in \{I\} \cup \{J_k\}_{k \in [K]}\\L' \neq L}}$ and $M_{\hat{\rho}}$ and $M'$ is the lifting of the stable matchings $\{M_0^{L'}\}_{\substack{L' \in \{I\} \cup \{J_k\}_{k \in [K]}\\L' \neq L}}$ and $M^{\hat{\rho}}$.
Note that because $\rho \in \cR_a(I^\sqcup)$, the partner of $(a, A^L)$ is changed when $\rho$ is eliminated, this implies that $M_{\hat{\rho}}(a) \neq M^{\hat{\rho}}(a)$ and the elimination of $\hat{\rho}$ in $L$ changes the partner of $a$ from $b_{j-1}^L$ to $b_{j}^L$ for some $j \in [m^L]$. We define $\Psi(\rho) = (b_j^L, B^L)$. 

We next show that $\Psi$ is indeed a one-to-one correspondence. Fix $(b_j^L, B^L)$ where $L \in \{I\} \cup \{J_k\}_{k \in [K]}$ and $j \in [m^L]$, there exists a rotation $\hat{\rho} \in \cR(L)$ such that eliminating $\hat{\rho}$ changes the partner of $a$ from $b_{j-1}^L$ to $b_j^L$. The lifting $
    M' = {\cal L}(M^I_0, M_0^{J_1}, \dots, M^{\hat{\rho}} ,\dots, M^{J_K}_0),
    $ is an immediate predecessor of the lifting 
$M = {\cal L}(M^I_0, M_0^{J_1}, \dots, M_{\hat{\rho}} ,\dots, M^{J_K}_0)$ in $I^\sqcup$ (otherwise if there exists $M'' \in \cS(I^\sqcup)$ such that $M' <^\sqcup M'' <^\sqcup M$ then $\cP_L(M'')$ would be a stable matching of $L$ such that $M^{\hat{\rho}} <^L \cP_L(M'') <^L M_{\hat{\rho}}$ which is absurd). Hence there exists a rotation $\rho = \rho(M, M') \in \cR(I^\sqcup)$ such that $\Psi(\rho) = (b_j^L, B^L)$. The rotation $\rho$ is moreover unique. Indeed, consider two rotations $\rho_1 = \rho(M_1, M_1'), \rho_2 = \rho(M_2, M_2')$ such that
$$
M_1 = {\cal L}(M^I_0, M_0^{J_1}, \dots, M_{\hat{\rho}_1} ,\dots, M^{J_K}_0)
    \quad \text{ and } \quad
M_1' = {\cal L}(M^I_0, M_0^{J_1}, \dots, M^{\hat{\rho}_1} ,\dots, M^{J_K}_0),
    $$
for some $L_1 \in \{I\} \cup \{J_k\}_{k \in [K]}$ and $\hat{\rho}_1 \in \cR(L_1)$
and that
$$
    M_2 = {\cal L}(M^I_0, M_0^{J_1}, \dots, M_{\hat{\rho}_2} ,\dots, M^{J_K}_0)
    \quad \text{ and } \quad
    M'_2 = {\cal L}(M^I_0, M_0^{J_1}, \dots, M^{\hat{\rho}_2} ,\dots, M^{J_K}_0),
    $$
for some $L_2 \in \{I\} \cup \{J_k\}_{k \in [K]}$ and $\hat{\rho}_2 \in \cR(L_2)$ and assume that $\Psi(\rho_1) = \Psi(\rho_2) = (b_j^{L}, B^{L})$ then $L=L_1=L_2$ and, as we argued above, $\hat{\rho}_1$ and $\hat{\rho}_2$ change the partner of $a$ to $b_{j}^L$ when eliminated. Since there is only a unique rotation in $\cR(L)$ that can change the partner of $a$ to $b_{j}^L$ when eliminated (see Lemma~\ref{lem:rotation-all-or-nothing}) $\hat{\rho}_1 = \hat{\rho}_2$ implying that $\rho_1 = \rho_2$.

{\noindent \bf Image of $(a, A^L)$ by $M^\rho$ and $M_{\rho}$.} Consider $\rho \in \cR(I^{\sqcup})$
such that $\rho = \rho(M, M')$ where
    $$
    M = {\cal L}(M^I_0, M_0^{J_1}, \dots, M_{\hat{\rho}} ,\dots, M^{J_K}_0)
    \quad \text{ and } \quad
    M' = {\cal L}(M^I_0, M_0^{J_1}, \dots, M^{\hat{\rho}} ,\dots, M^{J_K}_0),
    $$
for some $L \in \{I\} \cup \{J_k\}_{k \in [K]}$ and $\hat{\rho} \in \cR(L)$.

If $\rho \in \cR_a(I^\sqcup)$ such that $\Psi(\rho) = (b_j^L, B^L)$, then $(a, A^L) (b_j^L, B^L) \in \rho^+$ and $(a, A^L) (b_{j-1}^L, B^L) \in \rho^-$ and hence $M^\rho((a, A^L)) = (b_j^L, B^L)$ and $M_\rho((a, A^L)) = (b_{j-1}^L, B^L)$.

If $\rho \in \cR(I^{\sqcup}) \setminus \cR_a(I^{\sqcup})$ then $M((a, A^L)) = M'((a, A^L))$ implying that the image of $a$ by $M^\rho$ and $M_{\rho}$ is the same and is equal to $b^L_0$ by definition of $M^\rho$ and $M_\rho$. This implies in turn that $M((a, A^L)) = M'((a, A^L)) = (b^L_0, A^L)$.

Next, if $\rho \in \cR_a(I^{\sqcup})$ and $L' \in \{I\} \cup \{J_k\}_{k \in [K]} \setminus \{L\}$ or $\rho \in \cR(I^{\sqcup}) \setminus \cR_a(I^{\sqcup})$ and $L' \in \{I\} \cup \{J_k\}_{k \in [K]}$. Recall that $M^\rho$ is the stable matching corresponding to the upper-closure of $\{\rho\}$ in the poset $\cR(I^{\sqcup})$ and hence $M^\rho$ is the largest stable matching (with respect to the order on stable matchings $>^\sqcup$) where $\rho$ is eliminated. Hence, $M^\rho((a, A^{L'})) \succeq^\sqcup_{(a, A^{L'})} M((a, A^{L'}))=(b_0^{L'}, B^{L'})$. But since the lifting 
$M = {\cal L}(M^I_0, M_0^{J_1}, ,\dots, M^{J_K}_0)$
is the student optimal stable matching of $I^\sqcup$ (otherwise one can find a projection onto some instance $L'$ that is a stable matching and where a student has strictly better partner than in $M_0^{L'}$). Then,  $M^\rho((a, A^{L'})) \preceq_{(a, A^{L'})} (b_0^{L'}, B^{L'})$ and hence $M^\rho((a, A^{L'})) = (b_0^{L'}, B^{L'})$. Finally, since $M^\rho((a, A^{L'}))$ is the stable matching we get by eliminating $\rho$ from $M_{\rho}$, and $\rho$ does not change the 
partner of $(a, A^{L'})$, we have $M^\rho((a, A^{L'})) =M_\rho((a, A^{L'}))=(b_0^{L'}, B^{L'})$. \hfill $\square$

The following claim describes the second order differentials of $f_{a,k}^3$.

\begin{claim}
\label{clm:secondorderf3}
Let $\rho \neq \rho' \in \cR(I^{\sqcup})$. If $\rho, \rho' \in \cR_a(I^{\sqcup})$ and $(\Psi(\rho), \Psi(\rho')) = ((b^I_j,B^I), (b^{J_k}_{l}, B^{J_k}))$ or $((b^{J_k}_{l}, B^{J_k}), (b^I_j,B^I))$ for some $j \in [m^I]$ and $l \in [m^{J_k}]$, then, 
$$
\partial^2 (f^3_{a,k})_{\rho, \rho'} = \left|[R_a(b^I_{j-1}), R_a(b^I_{j})] \cap [R_a(b^{J_k}_{l-1}), R_a(b^{J_k}_{l})]\right|
$$
otherwise $\partial^2 (f^3_{a,k})_{\rho, \rho'} = 0$.
In particular for every $M \in \cS(I^{\sqcup})$,
$$
    \frac{1}{2}\sum_{\rho \neq \rho' \in R_M} \partial^2 (f^3_{a,k})_{\rho, \rho'} = \left|[R_a(M_0^{(I)}(a)), R_a(M^{(I)}(a))] \cap [R_a(M_0^{(J_k)}(a)), R_a(M^{(J_k)}(a))]\right|
    $$
\end{claim}
{\noindent \em Proof. }Recall that $f^3_{a,k}(M)=[R_a(M^{(J_k)}(a))-R_a(M^{(I)}(a))]^+$. First, if we have $(\Psi(\rho), \Psi(\rho')) = ((b^I_j,B^I), (b^{J_k}_{l}, B^{J_k}))$, then $[M^\rho\wedge M_{\rho'}]^{(I)}(a)=\min_{\succ_a}\{b_j^I,b_0^I\}=b_j^I$ and $[M^\rho\wedge M_{\rho'}]^{(J_k)}(a)=\min_{\succ_a}\{b_0^{J_k},b_{l-1}^{J_k}\}=b_{l-1}^{J_k}$, therefore we have $f^3_{a,k}(M^\rho\wedge M_{\rho'})=[R_a(b_{l-1}^{J_k})-R_a(b_j^I)]^+$. Similarly, we have $f^3_{a,k}(M_\rho\wedge M^{\rho'})=[R_a(b_{l}^{J_k})-R_a(b_{j-1}^I)]^+$, $f^3_{a,k}(M_\rho\wedge M_{\rho'})=[R_a(b_{l-1}^{J_k})-R_a(b_{j-1}^I)]^+$ and $f^3_{a,k}(M^\rho\wedge M^{\rho'})=[R_a(b_{l}^{J_k})-R_a(b_{j}^I)]^+$. Therefore, we obtain
\begin{align*}
    \partial^2 (f^3_{a,k})_{\rho,\rho'}&=f^3_{a,k}(M^\rho\wedge M_{\rho'})+f^3_{a,k}(M_\rho\wedge M^{\rho'})-f^3_{a,k}(M_\rho\wedge M_{\rho'})-f^3_{a,k}(M^\rho\wedge M^{\rho'})\\
    &=[R_a(b_{l-1}^{J_k})-R_a(b_j^I)]^++[R_a(b_{l}^{J_k})-R_a(b_{j-1}^I)]^+ -[R_a(b_{l-1}^{J_k})-R_a(b_{j-1}^I)]^+-[R_a(b_{l}^{J_k})-R_a(b_{j}^I)]^+\\
    &=\left|[R_a(b^I_{j-1}), R_a(b^I_{j})] \cap [R_a(b^{J_k}_{l-1}), R_a(b^{J_k}_{l})]\right|
\end{align*}
where the last equality follows from Lemma~\ref{lem:interval}. If $(\Psi(\rho), \Psi(\rho')) = ( (b^{J_k}_{l}, B^{J_k}), (b^I_j,B^I))$, we have the same formula by symmetry of $\rho,\rho'$. 

Second, suppose that $(\Psi(\rho), \Psi(\rho')) = ((b^L_j, B^L), (b^{L'}_{l}, B^{L'}))$ such that $L\notin\{I,J_k\}$, then $[M^\rho\wedge M_{\rho'}]^{(I)}(a)=[M_{\rho'}]^{(I)}(a)$ since $[M^\rho]^{(I)}(a)=M_0^{(I)}(a)=b_0^I$. Similarly, $[M^\rho\wedge M_{\rho'}]^{(J_k)}(a)=[M_{\rho'}]^{(J_k)}(a)$, and therefore we have $f^3_{a,k}(M^\rho\wedge M_{\rho'})=f^3_{a,k}(M_{\rho'})$. A similar argument implies that $f^3_{a,k}(M_\rho\wedge M^{\rho'})=f^3_{a,k}(M^{\rho'})$, $f^3_{a,k}(M_\rho\wedge M_{\rho'})=f^3_{a,k}(M_{\rho'})$, and $f^3_{a,k}(M^\rho\wedge M^{\rho'})=f^3_{a,k}(M^{\rho'})$. These facts imply that $\partial^2 (f^3_{a,k})_{\rho,\rho'}=0$. Similarly, if $L'\notin\{I,J_k\}$, we also have $\partial^2 (f^3_{a,k})_{\rho,\rho'}=0$.

Third, if $(\Psi(\rho), \Psi(\rho')) = ((b^I_j, B^I), (b^I_{l}, B^I))$, then let us assume without loss of generality that $j<l$ (note that $j \neq l$ since $\rho \neq \rho'$). Then, $[M^\rho\wedge M_{\rho'}]^{(I)}(a)=\min_{\succ_a}\{b_j^I,b_{l-1}^I\}=b_{l-1}^I$ and $[M^\rho\wedge M_{\rho'}]^{(J_k)}(a)=\min_{\succ_a}\{b_0^{J_k},b_0^{J_k}\}=b_0^{J_k}$. Therefore, $f^3_{a,k}(M^\rho\wedge M_{\rho'})=[R_a(b_0^{J_k})-R_a(b_{l-1}^I)]^+$. Similarly, we can deduce that $f^3_{a,k}(M_\rho\wedge M^{\rho'})=[R_a(b_0^{J_k})-R_a(b_l^I)]^+$, $f^3_{a,k}(M_\rho\wedge M_{\rho'})=[R_a(b_0^{J_k})-R_a(b_{l-1}^I)]^+$ and $f^3_{a,k}(M^\rho\wedge M^{\rho'})=[R_a(b_0^{J_k})-R_a(b_l^I)]^+$. Combining these, we obtain $\partial^2 (f^3_{a,k})_{\rho,\rho'}=0$. The case when $(\Psi(\rho), \Psi(\rho')) = ((b^{J_k}_j, B^{J_k}), (b^{J_k}_{l}, B^{J_k}))$ can be verified similarly. 

Lastly, if $\rho\notin \cR_a(I^\sqcup)$, then $M^\rho((a,A^L))=(M_0^{(L)}(a),B^L)$ and $M_\rho((a,A^L))=(M_0^{(L)}(a),B^L)$ for every $L\in\{I\}\cup\{J_k\}_{k\in K}$. Therefore, $M^\rho\wedge M_{\rho'}((a,A^L))=M_{\rho'}((a,A^L))$ and similarly, $M_\rho\wedge M^{\rho'}((a,A^L))=M^{\rho'}((a,A^L))$, $M_\rho\wedge M_{\rho'}((a,A^L))=M_{\rho'}((a,A^L))$ and $M^\rho\wedge M^{\rho'}((a,A^L))=M^{\rho'}((a,A^L))$, implying that $\partial^2 (f^3_{a,k})_{\rho, \rho'}=0$. The same result holds for the case when $\rho'\notin \cR_a(I^\sqcup)$ by symmetry.

Now, using the formula of $\partial^2 (f^3_{a,k})_{\rho,\rho'}$ above, we have
\begin{align*}
      &\frac{1}{2}\sum_{\rho \neq \rho' \in R_M} \partial^2 (f^3_{a,k})_{\rho, \rho'} \\
      &=\frac{1}{2}\left(\sum_{\substack{\rho\in R_M:\exists j\, s.t. \Psi(\rho) = (b^I_j, B^I) \\
     \rho'\in R_M:\exists l\, s.t. \Psi(\rho') = (b^{J_k}_l, B^{J_k})}}\partial^2 (f^3_{a,k})_{\rho, \rho'} + \sum_{\substack{\rho'\in R_M:\exists j\, s.t. \Psi(\rho') = (b^I_j, B^I) \\
     \rho\in R_M:\exists l\, s.t. \Psi(\rho) = (b^{J_k}_l,B^{J_k})}}\partial^2 (f^3_{a,k})_{\rho, \rho'}\right)\\
     &=\sum_{\rho\in R_M:\exists j\, s.t. \Psi(\rho) = (b^I_j, B^I)}\left(\sum_{\rho'\in R_M:\exists l\, s.t. \Psi(\rho') = (b^{J_k}_l, B^{J_k})}\left|[R_a(b^I_{j-1}), R_a(b^I_{j})] \cap [R_a(b^{J_k}_{l-1}), R_a(b^{J_k}_{l})]\right|\right)\\
     &=\sum_{\rho\in R_M:\exists j\, s.t. \Psi(\rho) = (b^I_j, B^{I})}\left|[R_a(b^I_{j-1}), R_a(b^I_{j})] \cap [R_a(M_0^{(J_k)}(a)), R_a(M^{(J_k)}(a))]\right|\\
     &= \left|[R_a(M_0^{(I)}(a)), R_a(M^{(I)}(a))] \cap [R_a(M_0^{(J_k)}(a)), R_a(M^{(J_k)}(a))]\right|,
\end{align*}
implying the desired formula.\hfill $\square$

Note that all of the second order differentials are non-negative and hence (ii) holds. The next claim gives the first order differentials of $f^3_{a,k}$.
\begin{claim}
\label{clm:firstorderf3}
For every $\rho \in \cR(I^{\sqcup})$,  if $\rho \in \cR_a(I^{\sqcup})$ and,
\begin{itemize}
    \item $\Psi(\rho) = (b_j^I, B^I)$ then,
$
\partial (f^3_{a,k})_{\rho} = 
        [R_a(M_0^{(J_k)}(a)) - R_a(b_j^I)]^+ - [R_a(M_0^{(J_k)}(a)) - R_a(b_{j-1}^I)]^+.
$
\item $\Psi(\rho) = (b_l^{J_k}, B^{J_k})$ then, $
\partial (f^3_{a,k})_{\rho} = [R_a(b_l^{J_k}) - R_a(M_0^{(I)}(a))]^+ - [R_a(b_{l-1}^{J_k}) - R_a(M_0^{(I)}(a))]^+.
$
\end{itemize}
Otherwise $\partial (f^3_{a,k})_{\rho} = 0$.

In particular for every $M \in \cS(I^{\sqcup})$, 
\begin{align*}
    \sum_{\rho \in R_M} \partial (f^3_{a,k})_\rho &= [R_a(M_0^{(J_k)}(a)) - R_a(M^{(I)}(a))]^+ + [R_a(M^{(J_k)}(a)) - R_a(M_0^{(I)}(a))]^+ 
    \\&\quad - 2[R_a(M_0^{(J_k)}(a)) - R_a(M_0^{(I)}(a))]^+
\end{align*}
\end{claim}

{\noindent \em Proof. } If $\Psi(\rho) = (b_j^I, B^I)$, then by Claim~\ref{cl:rotation-per-instance}, $[M^\rho]^{(I)}(a)=b_j^I$, $[M^\rho]^{(J_k)}(a)=b_0^{J_k}=M_0^{(J_k)}(a)$ and $[M_\rho]^{(I)}(a)=b_{j-1}^I$, $[M_\rho]^{(J_k)}(a)=b_0^{J_k}=M_0^{(J_k)}(a)$. Therefore, by definition we have $\partial (f^3_{a,k})_{\rho} = f^3_{a,k}(M^\rho)-f^3_{a,k}(M_\rho)=[R_a(M_0^{(J_k)}(a)) - R_a(b_j^I)]^+ - [R_a(M_0^{(J_k)}(a)) - R_a(b_{j-1}^I)]^+$. The case of $\Psi(\rho) = (b_l^{J_k}, B^{J_k})$ follows by a similar argument. On the other hand, if $\rho\notin \cR_a(I^\sqcup)$, then $M^\rho((a,A^L))=(M_0^{(L)}(a),B^L)$ and $M_\rho((a,A^L))=(M_0^{(L)}(a),B^L)$ for every $L\in\{I\}\cup\{J_k\}_{k\in K}$. Therefore, $\partial (f^3_{a,k})_{\rho} = f^3_{a,k}(M^\rho)-f^3_{a,k}(M_\rho)=[R_a(M_0^{(J_k)}(a)) - R_a(M_0^{(I)}(a))]^+ - [R_a(M_0^{(J_k)}(a)) - R_a(M_0^{(I)}(a))]^+=0$.

Now, using the first order differential computed above, we have
\begin{align*}
     \sum_{\rho \in R_M} \partial (f^3_{a,k})_\rho&=\sum_{\rho\in R_M:\exists j\, s.t. \Psi(\rho) = (b^I_j, B^I)} [R_a(M_0^{(J_k)}(a)) - R_a(b_j^I)]^+ - [R_a(M_0^{(J_k)}(a)) - R_a(b_{j-1}^I)]^+\\
     &\quad +   \sum_{\rho\in R_M:\exists l\, s.t. \Psi(\rho) = (b^{J_k}_l, B^{J_k})}[R_a(b_l^{J_k}) - R_a(M_0^{(I)}(a))]^+ - [R_a(b_{l-1}^{J_k}) - R_a(M_0^{(I)}(a))]^+  \\
     &= [R_a(M_0^{(J_k)}(a)) - R_a(M^{(I)}(a))]^+ - [R_a(M_0^{(J_k)}(a)) - R_a(M_0^{(I)}(a))]^+\\
     &\quad +[R_a(M^{(J_k)}(a)) - R_a(M_0^{(I)}(a))]^+ - [R_a(M_0^{(J_k)}(a)) - R_a(M_0^{(I)}(a))]^+\\
     &= [R_a(M_0^{(J_k)}(a)) - R_a(M^{(I)}(a))]^+ + [R_a(M^{(J_k)}(a)) - R_a(M_0^{(I)}(a))]^+\\
     &\quad- 2[R_a(M_0^{(J_k)}(a)) - R_a(M_0^{(I)}(a))]^+.
\end{align*}
\hfill $\square$

Finally, from Claims~\ref{clm:secondorderf3} and~\ref{clm:firstorderf3} we conclude that for every $M \in \cS(I^{\sqcup})$,
\begin{align*}
    &f^3_{a,k}(M_0) + \sum_{\rho \in R_M} \partial (f^3_{a,k})_\rho - \frac{1}{2} \sum_{\rho \neq \rho'} \partial (f^3_{a,k})_{\rho, \rho' \in R_M}  
    \\
    &= [R_a(M_0^{(J_k)}(a)) - R_a(M^{(I)}(a))]^+ + [R_a(M^{(J_k)}(a)) - R_a(M_0^{(I)}(a))]^+ - [R_a(M_0^{(J_k)}(a)) - R_a(M_0^{(I)}(a))]^+ 
    \\
    &\quad - \left|[R_a(M_0^{(I)}(a)), R_a(M^{(I)}(a))] \cap [R_a(M_0^{(J_k)}(a)), R_a(M^{(J_k)}(a))]\right|\\
    &= [R_a(M^{(J_k)}(a)) - R_a(M^{(I)}(a))]^+ = f^3_{a,k}(M),
\end{align*}
where the second to last equality follows from Lemma~\ref{lem:interval}. This proves (iii).

\subsection{Implicit Second-Stage Distribution: Omitted Proofs}

\subsubsection{Proof of Theorem~\ref{thm:imp-hard}}

Suppose there exists an algorithm $\cA$ such that given any instance $\cI$ of \eqref{eq:2stagesto}, where the second-stage distribution is given by a sampling oracle, solves the two-stage problem over $\cI$ (in the sense that it computes at least one of the two, the optimal value or the optimal solution) in time and number of calls to the sampling oracle that is polynomial in the input size ${\sf N}^{\sf imp}$. Consider an undirected graph $G(V,E)$. 

Let $\alpha > 0$. Consider the following instance $\cI_G(\alpha)$ of the two-stage problem that we will use to count the number of vertex covers of $G$.

\medskip
{\noindent \bf Construction of $\cI_G(\alpha)$.}

\begin{itemize}
    \item {\em Aggregate Market:} The set of students $A$ is: For each vertex $v \in V$, let $\{v,w_1\}, \dots, \{v, w_{\degree{v}}\}$ be the edges adjacent to $v$, where $\degree{v}$ denotes the degree of vertex $v$. Add a student for each edge $\{v,w_i\}$ for $i \in [\degree{v}]$. We denote this student by $(v,w_i)$ and refer to such student as a $G$-student. Add three more students that we denote by $a$, $a'$ and $a''$ and refer to these students as dummy students. In total we have $\sum_{v \in V} \degree{v} + 3 = 2|E|+3$ students. For the set of schools $B$, add a school for each edge $\{v,w\} \in E$. We denote such school by $\{v,w\}$ and refer to such school as a $G$-school. Add three more schools that we denote by $b$, $b'$ and $b''$ and refer to such schools by as dummy schools. Hence, we have a total of $|E|+3$ schools. Now the orders $\AgentGreater{}{}$ are as follows:
    
    Student $(v,w)$ has the order $\{v,w\} \AgentGreater{(v,w)}{} \emptyset$. Student $a$ has the order $G\text{-schools} \AgentGreater{a}{} b  \AgentGreater{a}{} b' \AgentGreater{a}{} \emptyset$. Student $a'$ has the order $b' \AgentGreater{a'}{} b''  \AgentGreater{a'}{} \emptyset$. Student $a''$ has the order $b'' \AgentGreater{a''}{} b'  \AgentGreater{a''}{} \emptyset$. School $\{v,w\}$ has the order $(v,w) \AgentGreater{\{v,w\}}{} (w,v) \AgentGreater{\{v,w\}}{} a \AgentGreater{\{v,w\}}{} \emptyset,$
    where we note that the order of the first two students $(v,w)$ and $(w,v)$ does not matter and can be switched. School $b$ has the order $a \AgentGreater{b}{} \emptyset$. School $b'$ has the order $a \AgentGreater{b'}{} a'' \AgentGreater{b'}{} a' \AgentGreater{b'}{} \emptyset$. Finally, school $b''$ has the order $a' \AgentGreater{b''}{} a'' \AgentGreater{b''}{} \emptyset$.
    \smallskip
    \item {\em First-stage instance:}
    $I_1 = (A, B, \succ, q)$, where all the school quotas are $1$.
    
    \smallskip
    \item {\em Second-stage instance:} The second-stage distribution $\cD$ is such that sampling $J \sim \cD$ consists of the following: Sample a subset $\hat{S}$ of vertices of $V$ uniformly randomly (by including each element $v \in V$ independently with probability $1/2$). Then for each $v \in \hat{S}$, let $\{v,w_1\}, \dots, \{v, w_{\degree{v}}\}$ be the edges adjacent to $v$, all the students $(v,w_1), \dots, (v, w_{\degree{v}})$ will be present in the second-stage. Finally, the set of all $G$-schools, $b'$, $b''$, $a$, $a'$ and $a''$ will also be present in the second-stage. The quotas are still $1$ for all schools.

    \smallskip
    \item {\em Costs:} $c_1$ is such that $c_1(a'b')=-\alpha$ and $c_1(e)=0$ for every other pair $e \in A^+ \times B^+$, $\lambda = 1$, and $c_2 \equiv 0$.

\end{itemize}

\medskip
{\noindent \bf Relation to the number of vertex covers of $G$.} We now show that the instances $\cI_G(\alpha)$ can be used to compute the number of vertex covers of $G$. Let us begin by relating the optimal value of $\cI_G(\alpha)$ to the number of vertex covers of $G$.

First of all, note that in any first-stage stable matching, it must be the case that $G$-students and $G$-schools are matched only to each other with a subset of $G$-students matched to the outside option (note that $G$-schools are the only acceptable--better than outside option--schools for the $G$-students, that each $G$-school has its corresponding $G$-students as top choices, and that there are more $G$-students then $G$-schools), $a$ is matched to $b$ as they are top choice of each other (when the $G$-schools are taken), and hence $a'$ and $a''$ are matched to $b'$ and $b''$. Note that $a'$ prefers $b'$ to $b''$, $b'$ prefers $a''$ to $a'$, $a''$ prefers $b''$ to $b'$ and finally $b''$ prefers $a'$ to $a''$ and hence, switching the partners of $a'$ and $a''$ in any first-stage stable matching will still give a stable matching. When $a'$ and $a''$ are matched to $b'$ and $b''$ respectively, we say that the matching is of type $1$, such matching has a first-stage cost of $-\alpha$. When the reverse is true, we say that the matching is of type $0$ and such matching has cost $0$.

Let $N^{\sf vc}_G$ denote the number of vertex covers of $G$. We now show that if the first-stage matching is of type $i$ then the total cost paid across both stages is given by,
$
\left(\frac{N^{\sf vc}_G}{2^{|V|-1}}-\alpha\right) \cdot i + \frac{N^{\sf vc}_G}{2^{|V|-1}}.
$

First of all, note that whatever instance was sampled in the second-stage, all $G$-students will always improve their matching in the second-stage and hence their rank increase is always $0$. To show this, it is sufficient to show that any $G$-student who was matched to its corresponding $G$-school in the first-stage and was sampled in the second-stage, will still be matched to its corresponding $G$-school in the second-stage. Let $(v,w)$ be one such student, then if $(v,w)$ is matched to the outside option in the second-stage, it means that $\{v,w\}$ is matched to a student it prefers more, the only possible such student is $(w,v)$. However, this implies that $\{v,w\}$ prefers $(w,v)$ to $(v,w)$, and $(w,v)$ was matched to the outside option in the first-stage (as its only better choice $\{v,w\}$ was taken by $(v,w)$) implying that  $\{v,w\}$ and $(w,v)$ is a blocking pair for the first-stage matching which is absurd.

Suppose now that the sampled $\hat{S}$ is a vertex cover of $G$. We claim that the rank increase of $a$ is exactly $1$, the rank increase of $a'$ is $i$ and finally the rank increase of $a''$ is $i+1$. In fact, if $\hat{S}$ is a vertex cover of $G$, every $G$-school $\{v,w\}$ will have at least one of $(v,w)$ or $(w,v)$ present in the second-stage that it prefers to every other student and hence $G$-schools will be entirely matched to $G$-students (the rest of $G$-students will be matched to the outside option as there are more $G$-students then $G$-schools and $G$-students prefer the outside option to any school other than their corresponding $G$-school). Now, since $b$ is not included in the second-stage, $a$ must be matched to $b'$ as they are top choice of each other (when $b$ and $G$-schools are not available for $a$). This implies a rank increase of $a$ of exactly $1$. Next, if the first-stage matching is of type $1$, i.e., $a'$ and $a''$ are matched to $b'$ and $b''$ respectively, $a'$ will have to change its partner taken by $a$ to $b''$ (note that $a'$ is the top choice of $b''$), this implies a rank increase of $a'$ of exactly $1$. At the same time, $a''$ will be matched to the outside option, and hence its rank increase is $2$. If on the other hand the first-stage matching is of type $0$, $a'$ will not change partner and hence have a rank increase $0$, and $a''$ will become matched to the outside option as its partner is taken by $a$, hence a rank increase of $a''$ of $1$. In conclusion, when $\hat{S}$ is a vertex cover, the total cost paid in the second-stage is given by $C_{\hat{S}}= 2 i + 2$.

Now when $\hat{S}$ is not a vertex cover of $G$, we claim that the rank increase of all the dummy vertices $a$, $a'$ and $a''$ is $0$. In fact, if $\hat{S}$ is not a vertex cover of $G$, this implies that there exists a $G$-school $\{v,w\}$ for which both associated $G$-students $(v,w)$ and $(w,v)$ are not in $S$. This implies that $a$ will be matched to one such school in the second-stage as $G$-schools are $a$'s top choice and $a$ is every $G$-school's top choice (in the absence of its corresponding $G$-students). In terms of rank increase, $a$ improves its matchings and hence has rank increase $0$. Now because we are minimizing the rank increase between the two-stages, and because any initial matching of $a'$ and $a''$ can be kept without altering stability (as none of their partners are taken by some other student), the total rank increase of $a'$ and $a''$ in an optimal second-stage solution is also $0$. In conclusion, the total cost paid in the second-stage in this case is $C_{\hat{S}} = 0$.

Given the above, we conclude that the total cost paid across both stages is
\begin{align*}
-\alpha \cdot i + \frac{1}{2^{|V|}}\sum_{\hat{S}} C_{\hat{S}}
= -\alpha \cdot i + \frac{1}{2^{|V|}}\sum_{\hat{S}:\text{is vertex cover}} C_{\hat{S}} = \left(\frac{N^{\sf vc}_G}{2^{|V|-1}}-\alpha\right) \cdot i + \frac{N^{\sf vc}_G}{2^{|V|-1}}
\end{align*}
Hence, when $\alpha > \frac{N^{\sf vc}_G}{2^{|V|-1}}$, the optimal first-stage stable matching must be of type $1$, if $\alpha < \frac{N^{\sf vc}_G}{2^{|V|-1}}$, the optimal first-stage stable matching must be of type $0$. Finally, since $\frac{N^{\sf vc}_G}{2^{|V|-1}} \in [2^{1-|V|},2]$, this implies that: (i) computing the optimal value of the two stage problem for the instance $\cI_G(0)$ would give the number of vertex covers given the above expression. (ii) computing the optimal solution (therefore knowing whether its of type $0$ or type $1$ and hence knowing whether $\alpha > \frac{N^{\sf vc}_G}{2^{|V|-1}}$ or $\alpha < \frac{N^{\sf vc}_G}{2^{|V|-1}}$) of at most $\log(\frac{2}{2^{1-|V|}})=|V|$ many instances $\cI(\alpha)$ gives the value of $\frac{N^{\sf vc}_G}{2^{|V|-1}}$ and hence that of $N^{\sf vc}_G$ by binary search. Note that all of our constructions are polynomial in $|V|$, and hence this gives a polynomial time algorithm to compute the number of vertex covers of $G$.\hfill\Halmos

\subsubsection{Proof of Lemma~\ref{lem:samplecomplexity}}
For every $k \in [K]$ and every first-stage stable matching $M^I \in \cS(I)$, let $G_k(M^I) = \min_{
    M^{J_k} \in \cS(J_k)}  c_2(M^{J_k}) + d(M^I, M^{J_k})$
and let,
$G(M^I) = \mathbb{E}_{J \sim \cD}\left[\min_{
    M^J \in \cS(J)}  c_2(M^J) + d(M^I, M^J)\right]$.
Note that $G(M^I) 
= \mathbb{E}_{J_k \sim \cD}\left[G_k(M^I)\right]$. The two-stage problem \eqref{eq:2stagesto} is given by,
\begin{align*}
    \min_{M^I \in \cS(I)} c_1(M^I) + G(M^I)
\end{align*}
while the \eqref{eq:saa} problem is given by,
\begin{align*}
    \min_{M^I \in \cS(I)} c_1(M^I) + \frac{1}{K} \sum_{k=1}^K G_k(M^I)
\end{align*}

Let $\epsilon = 4|A|(\max_{ab} |c_2(ab)| + \lambda |B|)\sqrt{\frac{\max\{|A|, |B|\}\log(3.88) - \log(\alpha)}{K}}$ and Let $\cF^\epsilon = \{M^I \in \cS(I) \;|\; \val{}{M^I} \leq \val{}{M^{I, *}} + \epsilon\}$ denote the set of $\epsilon$ approximate solutions of \eqref{eq:2stagesto}. We have,
\begin{align}
\label{eq:unionbound}
    \mathbb{P}(\hat{M}^{I, *} \notin \textstyle{\cF}^\epsilon) &\leq \sum_{M^I \notin \textstyle{\cF}^\epsilon} \mathbb{P}(M^I\text{ is an optimal solution of \eqref{eq:saa}})\notag
    \\
    &\leq \sum_{M^I \notin \textstyle{\cF}^\epsilon} \mathbb{P}\left(c_1(M^I) + \frac{1}{K} \sum_{k=1}^K G_k(M^I) \leq c_1(M^{I, *}) + \frac{1}{K} \sum_{k=1}^K G_k(M^{I, *})\right)\\
    &= \sum_{M^I \notin \textstyle{\cF}^\epsilon} \mathbb{P}\left(\frac{1}{K} \sum_{k=1}^K \left[c_1(M^I) + G_k(M^I) - c_1(M^{I, *}) - G_k(M^{I, *})\right] \leq 0\right),\notag
\end{align}
where the probability is taken over the samples $J_1, \dots, J_K \sim \cD$. 

Fix $M^I \notin \cF^\epsilon$ and let $
X_k = c_1(M^I) + G_k(M^I) - c_1(M^{I, *}) - G_k(M^{I, *})
,$ for all $k \in [K]$. Let us bound the probability,
$
\mathbb{P}\left(\frac{1}{K} \sum_{k=1}^K X_k \leq 0\right)
$.
We will use the classical Hoeffding's inequality stated below. 
\begin{lemma}[\cite{hoeffding1994probability}]
\label{lem:hoeffding}
    Let $X_1, \dots, X_K$ be independent random variables such that $a_k \leq X_k \leq b_k$ a.s. Let $S_K = X_1 + \dots + X_K$ and $t > 0$. Then, 
    $
    \mathbb{P}\left( \mathbb{E}(S_K) - S_K \geq t \right) \leq \exp \left(-\frac{2t^2}{\sum_{k=1}^{K}(b_k - a_k)^2}\right)
    $.
\end{lemma}

In particular, let $k \in [K]$, we have for every $M^{J_k} \in \cS(J_k)$,
$$
- |A|\max_{ab} |c_2(ab)| \leq  c_2(M^{J_k}) + d(M^I, M^{J_k}) \leq |A|\left(\max_{ab} |c_2(ab)| + \lambda |B|\right)
.$$
Hence,
$$
- |A|\max_{ab} |c_2(ab)| \leq G_k(M^I) \leq |A|\left(\max_{ab} |c_2(ab)| + \lambda |B|\right).
$$
Let
$$
a_k = c_1(M^I) - c_1(M^{I, *}) - |A|\left(2\max_{ab} |c_2(ab)| + \lambda |B|\right)
$$
and
$$
b_k = c_1(M^I) - c_1(M^{I, *}) + |A|\left(2\max_{ab} |c_2(ab)| + \lambda |B|\right)
.$$
Then the random variable $X_k$ is bounded such that $X_k \in [a_k, b_k]$. Now since, $M^I \notin \cF^{\epsilon}$ it holds that $\mathbb{E}(X_k) \geq \epsilon$. Let $S_K = \sum_{k=1}^K X_k$. We have,

\begin{align*}
     \mathbb{P}\left(\frac{1}{K} \sum_{k=1}^K X_k \leq 0\right) 
    & =
    \mathbb{P}\left(S_K \leq 0\right)=
    \mathbb{P}\left(\mathbb{E}(S_K) - S_K \geq \mathbb{E}(S_K)\right) \\
    &\leq 
    \mathbb{P}\left(\mathbb{E}(S_K) - S_K \geq K\epsilon\right)
    \leq \exp \left(-\frac{2K^2\epsilon^2}{\sum_{k=1}^{K}(b_k - a_k)^2}\right) \leq \frac{\alpha}{3.88^{\max\{|A|, |B|\}}},
\end{align*}
where the second inequality follows from Hoeffding's inequality Lemma~\eqref{lem:hoeffding}. Hence, by summing over all $M_1 \notin \cF^\epsilon$, equation~\eqref{eq:unionbound} implies that $\mathbb{P}(\hat{M}_1 \notin \textstyle{\cF}^\epsilon) \leq \frac{\alpha |\cF^\epsilon|}{3.88^{\max\{|A|, |B|\}}} \leq \alpha$,
where the last inequality follows from the fact that the number of stable matchings (and hence the cardinality of $\textstyle{\cF}^\epsilon$) of an instance between $n$ students and $n$ schools is at most $3.88^{n}$, \cite{palmer2022most}, and by the fact that if an instance is unbalanced, for example, has more students then schools, one can add dummy schools that are the least preferred to all students to get an instance with same number of students and schools and with a number of stable matchings that is at least the number of stable matchings of the original instance.\hfill\Halmos

\end{document}